\documentclass[11pt,reqno]{amsart}    

\usepackage[margin=1in]{geometry}
\usepackage{amsfonts,amsmath,amssymb,enumerate}    
\usepackage{amsthm, bm}    
\usepackage{mathrsfs}
\usepackage{setspace}
\onehalfspace
\usepackage{comment}
\usepackage[mathcal]{euscript}
\usepackage{enumerate}
\usepackage{amscd}
\usepackage{bm} 
\usepackage{cancel}
\usepackage{xcolor}
\usepackage{tikz}
\usetikzlibrary{arrows,matrix,snakes}
\usetikzlibrary{decorations.markings,shapes.geometric,shapes.misc}

\usepackage[latin1]{inputenc} 
\usepackage[T1]{fontenc} 
\usepackage[english]{babel}

\definecolor{lightred}{RGB}{255,127,127}
\definecolor{lightgreen}{RGB}{127,255,127}
\definecolor{lightblue}{RGB}{127,127,255}
\definecolor{linkcolor}{rgb}{0,0,0.6}
\usepackage[ pdftex,colorlinks=true,
pdfstartview=FitV,
linkcolor= linkcolor,
citecolor= linkcolor,
urlcolor= linkcolor,
hyperindex=true,
hyperfigures=false]
{hyperref}

\usepackage{tabularx} 

\theoremstyle{plain}
\newtheorem{theorem}{Theorem}[section]
\newtheorem{lemma}[theorem]{Lemma}
\newtheorem{proposition}[theorem]{Proposition}
\newtheorem{corollary}[theorem]{Corollary}

\newcommand{\nocontentsline}[3]{}
\newcommand{\tocless}[2]{\bgroup\let\addcontentsline=\nocontentsline#1{#2}\egroup}

\newcommand{\p}{\mathfrak{p}}
\newcommand{\g}{\mathfrak{g}}
\newcommand{\Lg}{{\null^L \g}}

\renewcommand{\sl}{\mathfrak{sl}}

\newcommand{\sO}{{\mathscr O}}

\newcommand{\C}{\mathbb{C}}
\newcommand{\CP}{\mathbb{P}^1}
\newcommand{\s}{\sigma}
\newcommand{\vs}{\varsigma}

\newcommand{\Id}{\textup{Id}}
\newcommand{\ms}[1]{{\mathsf #1}}
\newcommand{\Vc}[1]{\a_{#1}}
\newcommand{\M}{\mathcal{M}}

\newcommand{\Z}{\mathbb{Z}}

\numberwithin{equation}{section}
\usepackage{mathtools}
\theoremstyle{plain}
\newtheorem{thm}{Theorem}[section]

\newtheorem*{prop*}{Proposition}
\newtheorem{conj}[thm]{Conjecture}
\theoremstyle{definition}

\theoremstyle{remark}
\newtheorem{remarkx}{Remark}
\newenvironment{rem}
  {\pushQED{\qed}\remarkx}
  {\popQED\endremarkx}

\newtheorem{examplex}{Example}
\newenvironment{exmp}
  {\pushQED{\qed}\examplex}
  {\popQED\endexamplex}

\newtheorem*{rem*}{Remark}

\newcommand{\be}{\begin{equation}}    
\newcommand{\ee}{\end{equation}}    
\newcommand{\beu}{\begin{equation*}}    
\newcommand{\eeu}{\end{equation*}}    
\newcommand{\bea}{\begin{eqnarray}}    
\newcommand{\eea}{\end{eqnarray}}    
\newcommand{\beaa}{\begin{eqnarray*}}    
\newcommand{\eeaa}{\end{eqnarray*}}    
\newcommand{\bmx}{\begin{pmatrix}}    
\newcommand{\emx}{\end{pmatrix}}    
    
\newcommand{\ha}{\mbox{\small $\frac{1}{2}$}}

\newcommand{\n}{{\mathfrak n}}
\newcommand{\h}{{\mathfrak h}}    
\renewcommand{\b}{{\mathfrak b}}

\renewcommand{\a}{{\mathfrak a}}

\newcommand{\co}[1]{\check{ #1}}

\newcommand{\SimTo}{%
  \xrightarrow{\raisebox{-0.35 em}{\smash{\ensuremath{\sim}}}}%
}

\newcommand{\btp}{\begin{tikzpicture}[baseline=0pt,scale=0.9,line width=0.25pt]}    
\newcommand{\etp}{\end{tikzpicture}}

\newcommand{\atp}[1]{}

\DeclareMathOperator{\hgt}{ht}

\DeclareMathOperator{\res}{res}
\DeclareMathOperator{\tr}{tr}
\DeclareMathOperator{\rk}{rk}

\DeclareMathOperator{\Ad}{Ad}
\DeclareMathOperator{\ad}{ad}
\DeclareMathOperator{\Aut}{Aut}
\DeclareMathOperator{\Inn}{Inn}
\DeclareMathOperator{\Conn}{Conn}
\DeclareMathOperator{\Op}{Op}
\DeclareMathOperator{\op}{op}
\DeclareMathOperator{\MOp}{MOp}

\DeclareMathOperator{\Spec}{Spec}

\newcommand*{\longhookrightarrow}{\ensuremath{\lhook\joinrel\relbar\joinrel\rightarrow}}
\newcommand*{\longtwoheadrightarrow}{\ensuremath{\relbar\joinrel\twoheadrightarrow}}

\newcommand{\up}{\upsilon}

\newcommand{\vac}{v_0}

\DeclareMathOperator{\End}{End}

\DeclareMathOperator{\Der}{Der}
\DeclareMathOperator{\Lie}{Lie}
\DeclareMathOperator{\GL}{\mathsf{GL}}

\begin{document}

\title[Cyclotomic Gaudin models, Miura opers and flag varieties]{Cyclotomic Gaudin models,\\
Miura opers and flag varieties}

\author{Sylvain Lacroix}
\address{Univ Lyon, Ens de Lyon, Univ Claude Bernard, CNRS, Laboratoire de Physique, F-69342 Lyon, France.} \email{sylvain.lacroix@ens-lyon.fr}
\author{Beno\^{\i}t Vicedo}
\address{School of Physics, Astronomy and Mathematics, University of Hertfordshire, College Lane, Hatfield AL10 9AB, UK.} \email{benoit.vicedo@gmail.com}

\begin{abstract}
Let $\g$ be a semisimple Lie algebra over $\C$. Let $\nu \in \Aut \g$ be a diagram automorphism whose order divides $T \in \Z_{\geq 1}$.
We define \emph{cyclotomic $\g$-opers} over the Riemann sphere $\CP$ as gauge equivalence classes of $\g$-valued connections of a certain form, equivariant under actions of the cyclic group $\Z/ T\Z$ on $\g$ and $\CP$. It reduces to the usual notion of $\g$-opers when $T = 1$. 

We also extend the notion of Miura $\g$-opers to the cyclotomic setting. To any cyclotomic Miura $\g$-oper $\nabla$ we associate a corresponding cyclotomic $\g$-oper.
Let $\nabla$ have residue at the origin given by a $\nu$-invariant rational dominant coweight $\co\lambda_0$ and be monodromy-free on a cover of $\CP$. We prove that the subset of all cyclotomic Miura $\g$-opers associated with the same cyclotomic $\g$-oper as $\nabla$ is isomorphic to the $\vartheta$-invariant subset of the full flag variety of the adjoint group $G$ of $\g$, where the automorphism $\vartheta$ depends on $\nu$, $T$ and $\co\lambda_0$. The big cell of the latter is isomorphic to $N^\vartheta$, the $\vartheta$-invariant subgroup of the unipotent subgroup $N \subset G$, which we identify with those cyclotomic Miura $\g$-opers whose residue at the origin is the same as that of $\nabla$. In particular, the cyclotomic generation procedure recently introduced in \cite{CharlesVar} is interpreted as taking $\nabla$ to other cyclotomic Miura $\g$-opers corresponding to elements of $N^\vartheta$ associated with simple root generators.

We motivate the introduction of cyclotomic $\g$-opers by formulating two conjectures which relate them to the cyclotomic Gaudin model of \cite{VY1}. 
\end{abstract}

\input{epsf}

\maketitle
\setcounter{tocdepth}{1}

\vspace{10mm}

\tableofcontents


\newpage

\section{Introduction and motivation}

The Gaudin model \cite{Gaudin} is a quantum integrable long-range
spin chain of any length $N \in \Z_{\geq 1}$ which can be associated with any semisimple Lie algebra $\g$ over $\C$. Among its many different possible generalisations, we shall be interested in the so called \emph{cyclotomic Gaudin model}. It can be regarded as a particular example of a general family of Gaudin models associated with \emph{non}-skew-symmetric solutions of the classical Yang-Baxter equation, introduced in \cite{Skryp}.

The algebra of observables of the cyclotomic Gaudin model, as introduced by C. Young and one of the present authors in \cite{VY1}, is the $N$-fold tensor product $U(\g)^{\otimes N}$ of the universal enveloping algebra $U(\g)$. Given any $\g$-modules $M_i$ for $i = 1, \ldots, N$, the Hilbert space, or spin chain, is taken to be the $N$-fold tensor product $\bigotimes_{i=1}^N M_i$. To define the Hamiltonians, let $\{ I^a \}_{a=1}^{\dim \g}$ and $\{ I_a \}_{a=1}^{\dim \g}$ be dual bases of $\g$ with respect to a chosen non-degenerate bilinear form on $\g$. Let $T \in \Z_{\geq 1}$, pick a primitive $T^{\rm th}$-root of unity $\omega$ and consider the cyclic group $\Gamma \coloneqq \langle \omega \rangle \cong \Z/ T \Z$. Let $\s \in \Aut \g$ be an automorphism of $\g$ such that $\sigma^T = \Id$. Fix a collection of $N$ distinct complex numbers $z_i \in \C^\times$, $i=1,\ldots, N$ with disjoint $\Gamma$-orbits, \emph{i.e.} such that $z_i \neq \omega^k z_j$ for all distinct $i, j = 1, \ldots, N$ and $k = 0, \ldots, T-1$. The quadratic cyclotomic Gaudin Hamiltonians are defined as (see \cite{Skryp2})
\begin{equation}\label{CycloH}
\mathcal H_i \coloneqq \sum_{k=0}^{T-1} \sum_{\substack{j=1 \\ j\neq i}}^N \frac{I^{a (i)} \s^k I_a^{(j)}}{z_i-\omega^{-k} z_j} + \sum_{k=1}^{T-1} \frac{I^{a (i)} \s^k I_a^{(i)}}{(1 - \omega^k) z_i} \in U(\g)^{\otimes N}, \qquad i = 1, \ldots, N
\end{equation}
where for any $A \in U(\g)$ we let $A^{(i)}$ denote the element of $U(\g)^{\otimes N}$ with $A$ in the $i^{\rm th}$ tensor factor and the identity in every other factor. One checks directly that these Hamiltonians mutually commute, \emph{i.e.} $[\mathcal H_i, \mathcal H_j] = 0$ for all $i, j = 1, \ldots, N$, and commute with the diagonal action of the $\sigma$-invariant subalgebra $\g^\sigma$. The Hamiltonians \eqref{CycloH} belong to a large commutative subalgebra of $U(\g)^{\otimes N}$, the so called \emph{cyclotomic Gaudin algebra}, whose definition we briefly recall below. For $\rk \g \geq 2$ the latter contains also ``higher'' Gaudin Hamiltonians of degrees equal to the exponents of $\g$ plus one. Taking $T=1$, \emph{i.e.} $\Gamma = \{ 1 \}$, so that $\s = \Id$, we recover the usual quadratic Gaudin Hamiltonians
\begin{equation}
\label{GaudinH}
\mathcal H^{\rm G}_i \coloneqq \sum_{\substack{j=1 \\ j\neq i}}^N \frac{I^{a (i)} I_a^{(j)}}{z_i - z_j} \in U(\g)^{\otimes N}, \qquad i = 1, \ldots, N.
\end{equation}
The cyclotomic Gaudin algebra reduces in this case to the \emph{Gaudin algebra} \cite{Frenkel2}.

The construction of the Gaudin algebra by Feigin, Frenkel and Reshetikhin \cite{FFR} exploits
the commutative algebra structure of the subspace of singular vectors in the vacuum Verma module at the critical level over the untwisted affine Kac-Moody algebra $\widehat \g$ associated with $\g$. This approach was recently generalised to construct the cyclotomic Gaudin algebra by C. Young and one of the present authors in \cite{VY1}. The central ingredient in the construction of \cite{FFR} is the notion of coinvariant of an $N$-fold tensor product of $\widehat \g$-modules with respect to the algebra of rational functions $\CP \to \g$ vanishing at infinity and with poles at most at the marked points $z_i$, $i = 1, \ldots, N$. In the cyclotomic setting this gets replaced by the notion of \emph{cyclotomic} coinvariants \cite{FreSzcz, VY2}, \emph{i.e.} coinvariants with respect to the algebra of $\Gamma$-equivariant rational functions $\CP \to \g$ vanishing at infinity and regular away from $\omega^k z_i$, $i = 1, \ldots, N$, $k = 0, \ldots, T-1$, where $\omega \in \Gamma$ acts on $\CP$ by multiplication and on $\g$ as $\sigma$. Let $u \in \C^\times \setminus \{ z_1, \ldots, z_N\}$ be such that its $\Gamma$-orbit is disjoint from those of the $z_i$, and denote by $\mathscr O_u$ and $\mathscr K_u$ the local ring and local field at $u$, respectively. Let $\widehat \g_u \coloneqq \g (\mathscr K_u) \oplus \C K$ be the affine Kac-Moody algebra, where for any $\C$-algebra $R$ we let $\g(R) \coloneqq \g \otimes R$, and consider its subalgebra $\widehat \g_u^+ \coloneqq \g (\mathscr O_u) \oplus \C K$. The vacuum Verma module $\mathbb V_{0,u}^{\rm crit}(\g) \coloneqq U(\widehat \g_u) \otimes_{U(\widehat \g^+_u)} \C \vac$ is the $\widehat \g_u$-module induced from the one-dimensional $\widehat \g^+_u$-module $\C \vac$ on which $\g(\mathscr O_u)$ acts trivially and $K$ acts by the critical level (whose specific value depends on the choice of normalisation of the bilinear form on $\g$).
The subspace $\mathfrak z\big( \mathbb V_{0,u}^{\rm crit}(\g) \big) \coloneqq \{ X \in \mathbb V_{0,u}^{\rm crit}(\g) \,|\, \g(\mathscr O_u) X = 0 \}$ of singular vectors in $\mathbb V_{0,u}^{\rm crit}(\g)$ comes naturally equipped with the structure of a commutative algebra and the upshot of the construction of \cite{VY1}, generalising that of \cite{FFR} to the cyclotomic setting, is an algebra homomorphism
\begin{equation*}
\Psi_{(z_i), u}^\Gamma : \mathfrak z \big( \mathbb V_{0,u}^{\rm crit}(\g) \big) \longrightarrow U(\g)^{\otimes N}.
\end{equation*}
Working at the critical level ensures that this homomorphism is non-trivial since it is only then that the vacuum Verma module admits singular vectors not proportional to the vacuum $\vac$. In particular, $\mathfrak z \big( \mathbb V_{0,u}^{\rm crit}(\g) \big)$ always contains the quadratic vector $S = \ha I^a(-1) I_a(-1) \vac$. Here we fix once and for all a global coordinate $t$ on $\C \subset \CP$ and let $X(n) \coloneqq X \otimes (t - u)^n \in \g(\mathscr K_u)$ for any $X \in \g$ and $n \in \Z$. The image of $S$ under $\Psi_{(z_i), u}^\Gamma$ gives rise to the quadratic cyclotomic Gaudin Hamiltonians \eqref{CycloH} as $\mathcal H_j = \res_{z_j} \Psi_{(z_i), u}^\Gamma(S) du$ for each $j = 1, \ldots, N$.
The \emph{cyclotomic Gaudin algebra} $\mathscr Z_{(z_i)}^\Gamma(\g)$ is defined as the image of $\Psi_{(z_i), u}^\Gamma$. We thus obtain a surjective homomorphism of commutative $\C$-algebras
\begin{equation} \label{hom V0 to ZG}
\Psi^\Gamma_{(z_i), u} : \mathfrak z \big( \mathbb V_{0,u}^{\rm crit}(\g) \big) \longtwoheadrightarrow \mathscr Z^\Gamma_{(z_i)}(\g).
\end{equation}

Given any $\g$-modules $M_i$, $i = 1, \ldots, N$, one of the main problems in the study of the (cyclotomic) Gaudin model is to simultaneously diagonalise the quadratic (cyclotomic) Gaudin Hamiltonians \eqref{CycloH}, or indeed the entire (cyclotomic) Gaudin algebra $\mathscr Z^\Gamma_{(z_i)}(\g)$, on the tensor product $\bigotimes_{i=1}^N M_i$. If $\psi \in \bigotimes_{i=1}^N M_i$ is a joint eigenvector of the cyclotomic Gaudin algebra $\mathscr Z_{(z_i)}^\Gamma(\g)$, then the common eigenvalues of $\mathscr Z_{(z_i)}^\Gamma(\g)$ on $\psi$ are encoded in a $\C$-algebra homomorphism $\eta_\psi : \mathscr Z_{(z_i)}^\Gamma(\g) \to \C$, defined by $z \, \psi = \eta_\psi(z) \psi$ for every $z \in \mathscr Z_{(z_i)}^\Gamma(\g)$. In other words, the joint spectrum of $\mathscr Z_{(z_i)}^\Gamma(\g)$ on any spin chain $\bigotimes_{i=1}^N M_i$ forms a subset of the maximal spectrum of $\mathscr Z_{(z_i)}^\Gamma(\g)$. It is therefore of interest to first obtain a description of the spectrum $\Spec \mathscr Z_{(z_i)}^\Gamma(\g)$ before attempting to diagonalise the cyclotomic Gaudin Hamiltonians on any tensor product of $\g$-modules. In view of the homomorphism \eqref{hom V0 to ZG}, let us begin by recalling the description of the spectrum of $\mathfrak z \big( \mathbb V_{0,u}^{\rm crit}(\g) \big)$.

\medskip

The commutative algebra $\mathfrak z\big( \mathbb V_{0,u}^{\rm crit}(\g) \big)$ is naturally isomorphic, by a theorem of Feigin and Frenkel \cite{FF1} (see also the book \cite{Langlands}), to 
the \emph{classical $W$-algebra} $\mathcal W_u(\Lg)$ for the Langlands dual Lie algebra $\Lg$ of $\g$, whose Cartan matrix is the transpose of that of $\g$.
The classical $W$-algebra $\mathcal W_u(\g)$ is obtained by Drinfel'd-Sokolov reduction of the algebra of functions on the dual of the affine Kac-Moody algebra $\widehat \g_u$ \cite{Drin}. Its spectrum is isomorphic to the space of so called \emph{$\g$-opers} on the formal disc $D_u = \Spec \mathscr O_u$ around the point $u$. In order to introduce the notion of $\g$-oper it is convenient to first recall the description, due to Kostant \cite{KostantPrincipal, KostantPolynomial}, of the classical finite $W$-algebra $\mathcal W^{\rm fin}(\g, p_{-1})$ associated with a principal nilpotent element $p_{-1}$. Embedding the latter into an $\sl_2$-triple $\{ p_{-1}, 2 \co\rho, p_1 \}$, we set $\n \coloneqq \bigoplus_{i > 0} \g_i$ and $\b \coloneqq \bigoplus_{i \geq 0} \g_i$ using the $\Z$-grading $\g = \bigoplus_{i \in \Z} \g_i$ defined by $\ad_{\co\rho}$ and denote by $N$ the unipotent subgroup of the adjoint group of $\g$ with Lie algebra $\n$. The algebra $\mathcal W^{\rm fin}(\g, p_{-1})$ is isomorphic to the algebra of $N$-invariant polynomial functions on $p_{-1} + \b$, or equivalently to the algebra of polynomial functions on the Slodowy slice $p_{-1} + \a$, where $\a$ is the centraliser of $p_1$. In other words, we have isomorphisms $\Spec \mathcal W^{\rm fin}(\g, p_{-1}) \simeq (p_{-1} + \b) / N \simeq p_{-1} + \a$. To emphasise the parallel with $\g$-opers, we will call $\Op_\g^{\rm fin} \coloneqq (p_{-1} + \b) / N$ the space of \emph{finite $\g$-opers}. We denote by $[X]_\g$ the class of $X \in p_{-1} + \b$ in $\Op_\g^{\rm fin}$.

A $\g$-oper on $D_u$ is then defined as an equivalence class of connections of the form $d + p_{-1} dt + v \, dt$, with $v \in \b(\mathscr O_u)$, modulo the gauge action of the group $N(\mathscr O_u)$ whose definition we recall in \S\ref{TwistedAlgebra}. Each class admits a canonical representative of the form $d + p_{-1} dt + c\, dt$, with $c \in \a(\mathscr O_u)$, which can be regarded as an affine analog of the Slodowy slice. Denoting the space of $\g$-opers on $D_u$ by $\Op_\g(D_u)$, we have an isomorphism of varieties \cite{FF1}
\begin{equation} \label{FF isomorphism}
\Spec \mathfrak z\big( \mathbb V_{0,u}^{\rm crit}(\g) \big) \simeq \Op_{\Lg}(D_u).
\end{equation}
Now the surjective homomorphism \eqref{hom V0 to ZG} induces an injective map
\begin{equation} \label{hom V0 to ZG Spec}
\Spec \mathscr Z^\Gamma_{(z_i)}(\g) \longhookrightarrow \Spec \mathfrak z \big( \mathbb V_{0,u}^{\rm crit}(\g) \big)
\end{equation}
between the corresponding spectra. We may thus regard $\Spec \mathscr Z^\Gamma_{(z_i)}(\g)$ as a subvariety of $\Op_{\Lg}(D_u)$.
In the non-cyclotomic case, $\Gamma = \{ 1 \}$, it was shown by Frenkel \cite{Frenkel} that the usual Gaudin algebra $\mathscr Z_{(z_i)}(\g) \coloneqq \mathscr Z^{\{ 1 \}}_{(z_i)}(\g)$ is isomorphic to the algebra of functions on the space of \emph{global} $\Lg$-opers on the complex projective line $\CP$ with regular singularities at the points $z_i$, $i = 1, \ldots, N$ and at infinity. In other words, denoting the space of such $\Lg$-opers by $\Op_{\Lg}(\CP)^{\text{RS}}_{(z_i), \infty}$ which we recall the definition of in \S\ref{DS}, we have the isomorphism
\begin{equation} \label{Op CP iso}
\Spec \mathscr Z_{(z_i)}(\g) \simeq \Op_{\Lg}(\CP)^{\text{RS}}_{(z_i), \infty}.
\end{equation}
The injective map \eqref{hom V0 to ZG Spec} therefore corresponds, in the case $\Gamma = \{ 1 \}$, to the restriction of a global $\Lg$-oper in $\Op_{\Lg}(\CP)^{\text{RS}}_{(z_i), \infty}$ to the disc around a regular point $u \in \C \setminus \{ z_1, \ldots, z_N \}$.

\medskip

The first purpose of this article is to conjecture an analog of the isomorphism \eqref{Op CP iso} for describing the spectrum of the cyclotomic Gaudin algebra $\mathscr Z_{(z_i)}^\Gamma(\g)$. To this end we will introduce a notion of \emph{cyclotomic $\g$-oper} on $\CP$, which reduces to the usual notion of $\g$-oper on $\CP$, as given in \cite{MV}, when $\Gamma = \{ 1 \}$ (see also \cite{Frenkel, Frenkel2, BD}).
Fix a diagram automorphism $\nu \in \Aut\g$ whose order divides $T$. Given any automorphism $\up \in \Aut \g$ in the same class as $\nu$ in $\Aut \g/ \Inn \g$ and with the property that $\up^T = \Id$, such as the automorphism $\sigma$ entering the definition of the cyclotomic Gaudin model, we may consider the corresponding space $\Conn_\g^{\hat\up}(\CP)$ of $\Gamma$-equivariant meromorphic $\g$-valued connections on $\CP$ where $\omega \in \Gamma$ acts on $\CP$ by multiplication and on $\g$ by $\up$.
By constrast, consider meromorphic $\g$-valued connections, or $\g$-connections for short, of the form $d + p_{-1} dt + \ms v$ with $\ms v \in \Omega(\b) \coloneqq \b(\M) dt$, where $\M$ is the algebra of meromorphic functions on $\CP$. Requiring the existence of $\Gamma$-equivariant $\g$-connections of this form forces us to work with the representative $\vs = \Ad_{\omega^{-\co\rho}} \circ \nu \in \Aut \g$ of the class of $\nu$ in $\Aut \g/ \Inn \g$.
We therefore define the space of cyclotomic $\g$-opers as equivalence classes of $\Gamma$-equivariant $\g$-connections of the above form with $\omega \in \Gamma$ acting on $\g$ as $\vs$, modulo the gauge action of a $\Gamma$-invariant subgroup $N^{\hat\vs}(\M) \subset N(\M)$ whose definition is given in \S\ref{TwistedAlgebra}. Each such class also admits a unique canonical representative of the form $d + p_{-1} dt + \ms c$, with $\ms c \in \Omega^{\hat\vs}(\a)$,
where $\Omega^{\hat\vs}(\a)$ denotes the space of $\Gamma$-equivariant $\a$-valued differentials on $\CP$.
In \S\ref{GaudinModels} we will use the notion of cyclotomic $\g$-opers to conjecture the analog of the isomorphism \eqref{Op CP iso} for the cyclotomic Gaudin algebra $\mathscr Z_{(z_i)}^\Gamma(\g)$ in Conjecture \ref{conj: cyclo opers}.

\medskip

To motivate the central question concerning cyclotomic $\g$-opers addressed in this article we will use the analogy with classical finite $W$-algebras.
Recall the Slodowy slice $p_{-1} + \a$ which is transverse to the $N$-orbits in $p_{-1} + \b$. Another transverse slice is given by $p_{-1} + \h$ where $\h = \g_0$ is a Cartan subalgebra. Unlike the Slodowy slice, however, it doesn't intersect each $N$-orbit uniquely. In other words, the canonical map $p_{-1} + \h \to \Op_\g^{\rm fin}$ which sends an element $X \in p_{-1} + \h \subset p_{-1} + \b$ to its class $[X]_\g \in \Op_\g^{\rm fin}$ is surjective, but not injective. By comparison with the affine case discussed below, we shall refer to $p_{-1} + \h$ as the space of \emph{finite Miura $\g$-opers}. If $\co\lambda \in \h$ is dominant, \emph{i.e.} $\langle \alpha_i, \co\lambda \rangle \geq 0$ for every simple root $\alpha_i$, $i \in I \coloneqq \{ 1, \ldots, \rk \g \}$ of $\g$, then the subset of all finite Miura $\g$-opers whose class in $\Op_\g^{\rm fin}$ coincides with $[p_{-1} - \co\lambda - \co\rho]_\g$
can be shown using results of Kostant \cite{KostantPrincipal, KostantPolynomial} to be in bijection with the Weyl group $W$ of $\g$.

Following the standard terminology in the non-cyclotomic setting \cite{MV, Frenkel}, we define a \emph{cyclotomic Miura $\g$-oper} on $\CP$ as a cyclotomic $\g$-connection of the form $d + p_{-1} dt + \ms u$ with $\ms u \in \Omega^{\hat\nu}(\h)$. Note that a cyclotomic Miura $\g$-oper $\nabla$ on $\CP$ is \emph{not} a cyclotomic $\g$-oper, but we can associate to it a cyclotomic $\g$-oper by taking its gauge equivalence class $[\nabla]_\Gamma$ under $N^{\hat\vs}(\M)$. Let us fix a cyclotomic Miura $\g$-oper $\nabla$ with trivial monodromy representation and whose underlying cyclotomic $\g$-oper $[\nabla]_\Gamma$ has regular singularities at most at the points $z_i \in \C^\times$, $i = 1, \ldots, N$ (and their $\Gamma$-orbits), the origin and infinity. We will further assume that $\nabla$ is of the form
\begin{equation} \label{cyclo Miura intro}
\nabla = d + p_{-1} dt - \frac{\co\lambda_0}{t} dt + \ms r
\end{equation}
for any $\nu$-invariant integral dominant coweight $\co\lambda_0 \in \h^\nu$ and with $\ms r \in \Omega^{\hat\nu}(\h)$ regular at the origin. We will show in \S\ref{sec: non-integral} how the integrality assumption on $\co\lambda_0$ can be weakened by going to a cover of $\CP$, which allows us to treat also the case of a $\nu$-invariant dominant coweight $\co\lambda_0$ such that $\langle \alpha_i, \co\lambda_0 \rangle \in \mathbb Q$ for all $i \in I$.
Note that in general $\ms r$, and hence the cyclotomic $\g$-connection $\nabla$, may have simple poles at points $x_j \in \C^\times$, $j = 1, \ldots, m$ (and their $\Gamma$-orbits) other than $z_i$, $i = 1, \ldots, N$. Cyclotomic Miura $\g$-opers of the above form whose residue at each $x_j$, $j = 1, \ldots, m$ is a simple coroot and whose residue at each $z_i$, $i = 1, \ldots, N$ is minus an integral dominant coweight, were shown in \cite{VY1} to be described by solutions of the \emph{cyclotomic Bethe ansatz equations} first introduced in \cite{Skryp2}.
Now as in the finite setting described above, a natural question in the present context is then to describe the space of all cyclotomic Miura $\g$-opers whose underlying cyclotomic $\g$-oper coincides with $[\nabla]_\Gamma$.
In the non-cyclotomic case, $\Gamma = \{ 1 \}$, it was proved independently by Mukhin and Varchenko \cite{MVb, MV} and Frenkel \cite{Frenkel} that this space is isomorphic to the flag variety $G/B_-$ where $G$ is the adjoint group of $\g$ and $B_-$ is the Borel subgroup with Lie algebra $\b_- \coloneqq \bigoplus_{i \leq 0} \g_i$.

In applications to the cyclotomic Gaudin model, it follows from the construction of \cite{VY1} that each cyclotomic Miura $\Lg$-oper $\nabla$ (not necessarily monodromy-free) with residues $- \lambda_i \in \h^\ast$ at $z_i$ for $i =1, \ldots, N$, $- \lambda_0 \in \h^{\ast, \nu}$ at the origin and $\lambda_\infty \in \h^{\ast, \nu}$ at infinity, built from a solution of the cyclotomic Bethe ansatz equations, corresponds to a joint eigenvector of $\mathscr Z_{(z_i)}^\Gamma(\g)$ of weight $\lambda_\infty \in \h^\ast$ in the tensor product $\bigotimes_{i=1}^N M_{\lambda_i}$ of Verma modules $M_{\lambda_i}$ with highest weights $\lambda_i$. We will conjecture in \S\ref{GaudinModels} that the corresponding common eigenvalues of the cyclotomic Gaudin algebra $\mathscr Z_{(z_i)}^\Gamma(\g)$ are determined by the underlying cyclotomic $\Lg$-oper $[\nabla]_\Gamma$, see Conjecture \ref{conj: Miura oper}. In the non-cyclotomic setting it was conjectured in \cite{Frenkel2} that for integral dominant highest weights $\lambda_i \in \h^\ast$, $i = 1, \ldots, N$ and $\lambda_\infty \in \h^\ast$, there is a bijection between monodromy-free $\Lg$-opers in $\Op_{\Lg}(\CP)^{\text{RS}}_{(z_i), \infty}$ with residues determined by the weights $\lambda_i$ at $z_i$ and $\lambda_\infty$ at infinity, and eigenvalues of the Gaudin algebra $\mathscr Z_{(z_i)}(\g)$ on the tensor product $\bigotimes_{i=1}^N V_{\lambda_i}$ of finite dimensional irreducible modules $V_{\lambda_i}$ with highest weights $\lambda_i$. In particular, it is believed that when the Bethe ansatz is incomplete some of these eigenvalues correspond to Miura $\Lg$-opers which do not arise from solutions of the Bethe ansatz equations.

\medskip

Before proving the cyclotomic counterpart of the isomorphism between the space of Miura $\g$-opers with a given underlying $\g$-oper and the flag variety $G/B_-$, we begin in \S\ref{Repro} by considering the effect of gauge transformations on $\nabla$ by elements of the form $e^{f E_k} \in N(\M)$ for some $f \in \M$ and $k \in I$. Here $E_i, H_i$ and $F_i$ for $i \in I$ denote the Chevelley-Serre generators of $\g$. It is well know that after such a gauge transformation, the new $\g$-connection $e^{f E_k} \nabla e^{- f E_k}$ is still a Miura $\g$-oper if and only if $f$ is a solution of some Riccati equation. The result of going from the old Miura $\g$-oper $\nabla$ to the new one $e^{f E_k} \nabla e^{- f E_k}$ is known as a \emph{reproduction} in the direction of the simple root $\alpha_k$ \cite{MVb}. In general, however, the new Miura $\g$-oper $e^{f E_k} \nabla e^{- f E_k}$ will no longer be cyclotomic. The problem of defining a cyclotomic version of the reproduction procedure, taking one cyclotomic Miura $\g$-oper to another, was first studied by Varchenko and Young in \cite{CharlesVar} who considered so called populations of solutions to the cyclotomic Bethe ansatz equations associated with an arbitrary Kac-Moody algebra $\g$. It was shown there that under certain conditions on the coweight $\co\lambda_0$ (or rather the \emph{weight} $\lambda_0 \in \h^{\ast, \nu}$ since solutions of the cyclotomic Bethe ansatz equations used in \cite{CharlesVar} correspond to cyclotomic $\Lg$-opers), it is possible to take one solution of the cyclotomic Bethe ansatz equations to another by performing a sequence of reproductions in the direction of other simple roots $\{ \alpha_{i_j} \}_{j=2}^n$ for some $n \in \Z_{\geq 2}$ with $i_j \in I$ in the orbit $\mathcal I \in I/\nu$ of $i_1 \coloneqq k \in I$ under the diagram automorphism $\nu : I \to I$.

When the Lie algebra $\g$ is semisimple, as we are considering, there are only two possible types of orbits $\mathcal I$. We address the issue of existence of cyclotomic reproductions in the present language by studying gauge transformations by elements in the subgroup generated by simple roots along the orbit $\mathcal I$. Specifically, we consider Riccati equations built from the cyclotomic Miura $\g$-oper $\nabla$ given in \eqref{cyclo Miura intro}. These are satisfied by the collection of meromorphic functions $f_j \in \M$ which appear in the individual gauge transformation parameters $e^{f_j E_{i_j}} \in N(\M)$.
Requiring the overall gauge transformation parameter $g = e^{f_n E_{i_n}} \ldots e^{f_1 E_{i_1}}$ to live in the $\Gamma$-invariant subgroup $N^{\hat\vs}(\M)$ imposes certain functional relations among the various functions $f_j$. If the latter are regular at the origin then we find that these functional relations can be satisfied if and only if $\co\lambda_0$ satisfies the following relation
\begin{equation} \label{lambda0 intro}
\ell_{\mathcal I} \langle \alpha_k, \co\lambda_0 + \co\rho \rangle \equiv 0 \;\;\text{mod}\;\; \frac{T}{|\mathcal I|},
\end{equation}
where $\ell_{\mathcal I} = 1$ or $2$ depending on the type of the orbit $\mathcal I$. A key step in our analysis is considering a regularisation of the given Riccati equations at the origin. In this case the new cyclotomic Miura $\g$-oper $g \nabla g^{-1}$ takes the same form as $\nabla$ in \eqref{cyclo Miura intro} but with some new differential $\ms r \in \Omega^{\hat\nu}(\h)$. On the other hand, without imposing any conditions on $\co\lambda_0$, we can always choose the overall gauge transformation parameter to be cyclotomic by letting one of the functions $f_j$ be singular. Moreover, in this case we find that the new cyclotomic Miura $\g$-oper takes the form
\begin{equation*}
g \nabla g^{-1} = d + p_{-1} dt - \frac{s^\nu_{\mathcal I} \cdot \co\lambda_0}{t} dt + \widetilde{\ms r}
\end{equation*}
for some $\widetilde{\ms r} \in \Omega^{\hat\nu}(\h)$, with $s^\nu_{\mathcal I}$ the simple reflection of the $\nu$-invariant Weyl group $W^\nu$ associated with the node $\mathcal I \in I / \nu$ of the folded diagram and where the dot denotes the shifted Weyl action.

\medskip

Let us now consider the effect of a gauge transformation with a more general parameter $g \in N^{\hat\vs}(\M)$ on the cyclotomic Miura $\g$-oper $\nabla$ in \eqref{cyclo Miura intro}, not one corresponding merely to the orbit of a simple root.
We begin in \S\ref{ReproGen} by describing the space of all cyclotomic Miura $\g$-opers $g \nabla g^{-1}$ of the same form as the original cyclotomic Miura $\g$-oper $\nabla$, namely \eqref{cyclo Miura intro} but for some different $\ms r \in \Omega^{\hat\nu}(\h)$. We refer to these cyclotomic Miura $\g$-opers as being \emph{generic} at the origin. Guided by the analysis of \S\ref{Repro} we introduce the following \emph{regularised} $\g$-connection
\begin{equation*}
\nabla_{\rm r} \coloneqq t^{- \co\lambda_0} \nabla t^{\co\lambda_0} = d + \sum_{k \in I} t^{\langle \alpha_k, \co\lambda_0 \rangle} F_k \, dt + \ms r.
\end{equation*}
The regularisation has the effect of modifying the automorphism of $\g$ in the $\Gamma$-equivariance property from $\vs$ to $\vartheta \coloneqq \Ad_{\omega^{-\co\lambda_0 - \co\rho}} \circ \nu \in \Aut \g$. In particular, regularising the new cyclotomic Miura $\g$-oper $g \nabla g^{-1}$ yields the $\g$-connection $(g \nabla g^{-1})_{\rm r} = g_{\rm r} \nabla_{\rm r} g_{\rm r}^{-1}$ where $g_{\rm r} = t^{- \co\lambda_0} \,g\, t^{\co\lambda_0} \in N^{\hat\vartheta}(\M)$. Since we are assuming $g \nabla g^{-1}$ to be of the same form as in \eqref{cyclo Miura intro} it follows that $(g \nabla g^{-1})_{\rm r}$ is regular at the origin and hence so is $g_{\rm r}$. Its initial condition $g_{\rm r}(0)$ then takes value in $N^\vartheta$. Conversely, given any $g_0 \in N^{\vartheta}$ we construct in Theorem \ref{CycReproReg} an element $g \in N^{\hat\vs}(\M)$ such that $g_{\rm r}$ is regular at the origin with $g_{\rm r}(0) = g_0$. The cyclotomic Miura $\g$-oper $g \nabla g^{-1}$ is then of the same form as in \eqref{cyclo Miura intro}. In this language, the relation \eqref{lambda0 intro} on $\co\lambda_0$ for a given orbit $\mathcal I \in I / \nu$ can now be seen as the condition for the existence of a $\vartheta$-invariant element in the subalgebra of $\n$ generated by the $E_i$ with $i \in \mathcal I$.

The residue at the origin of the cyclotomic $\g$-oper $[\nabla]_\Gamma$, as defined in \S\ref{sec: opers RS}, corresponding to the cyclotomic Miura $\g$-oper $\nabla$ given in \eqref{cyclo Miura intro} is the finite $\g^\nu$-oper $[p_{-1} - \co\lambda_0 - \co\rho]_{\g^\nu} \in \Op^{\rm fin}_{\g^\nu}$. The general form of a cyclotomic Miura $\g$-oper $\widetilde \nabla$ whose corresponding cyclotomic $\g$-oper $[\widetilde \nabla]_\Gamma$ has the same residue at the origin is
\begin{equation} \label{cyclo Miura w intro}
\widetilde \nabla = d + p_{-1} dt - \frac{w \cdot \co\lambda_0}{t} dt + \widetilde{\ms r}
\end{equation}
for some $w \in W^\nu$ and $\widetilde{\ms r} \in \Omega^{\hat\nu}(\h)$ regular at the origin. It follows from the above discussion that to reach such cyclotomic Miura $\g$-opers from $\nabla$ given in \eqref{cyclo Miura intro} one should apply a gauge transformation with parameter $g \in N^{\hat\vs}(\M)$ for which $g_{\rm r} = t^{-\co\lambda_0} g t^{\co\lambda_0}$ is singular at the origin. We prove in \S\ref{sec: iso flag} that the space of \emph{all} cyclotomic Miura $\g$-opers of the form $g \nabla g^{-1}$ with $g \in N^{\hat\vs}(\M)$ is isomorphic to the $\vartheta$-invariant subspace $(G/B_-)^\vartheta$ of the flag variety $G/B_-$. We show in \S\ref{sec: cell decomp} that this $\vartheta$-invariant subspace has a cell decomposition
\begin{equation*}
(G/B_-)^\vartheta = \bigsqcup_{w \in W^\nu} N^\vartheta \dot w B_-/ B_-.
\end{equation*}
In particular, the big cell $N^\vartheta B_-/ B_-$ is isomorphic to the space of generic cyclotomic Miura $\g$-opers $g \nabla g^{-1}$, \emph{i.e.} for which $g_{\rm r}$ is regular at the origin. We expect the cyclotomic Miura $\g$-opers $g \nabla g^{-1}$ of the form \eqref{cyclo Miura w intro} to be isomorphic to the cell $N^\vartheta \dot w B_-/ B_-$.

\subsubsection*{Acknowledgements}
We thank C. Young for useful discussions and for collaboration at an early stage of the project. We also thank I. Cherednik, F. Delduc, P. Etingof and A. Varchenko for useful comments, as well as M. Magro for comments on the draft. This work is partially supported by the program PICS 6412 DIGEST of CNRS and by the French Agence Nationale de la Recherche (ANR) under grant ANR-15-CE31-0006 DefIS.

\section{Notations and conventions}
\label{Conventions}

\subsection{Semisimple Lie algebras}
\label{SSLieAlgebras}

Let $\g$ be a finite-dimensional complex semisimple Lie algebra. Let $\h$ be a Cartan subalgebra of $\g$ and denote by $\Phi \subset \h^\ast$ the root system of $(\g, \h)$. Fix a basis of simple roots $\alpha_i$, $i \in I \coloneqq \{ 1, \ldots, \rk \g \}$ and let $\Phi^\pm$ denote the corresponding set of positive and negative roots.
We have the Cartan decomposition
\begin{equation} \label{Cartan decomp}
\g = \n_- \oplus \h \oplus \n, \qquad
\n \coloneqq \bigoplus_{\alpha \in \Phi^+} \C E_\alpha, \quad \n_- \coloneqq \bigoplus_{\alpha \in \Phi^+} \C F_\alpha.
\end{equation}
Introduce also the positive and negative Borel subalgebras $\b \coloneqq \h \oplus \n$ and $\b_- \coloneqq \h \oplus \n_-$. For any positive root $\alpha \in \Phi^+$, we define the root spaces $\g_\alpha \coloneqq \C E_\alpha$ and $\g_{-\alpha} \coloneqq \C F_\alpha$. The Chevalley-Serre generators $\co \alpha_i \in \h$, $E_i \coloneqq E_{\alpha_i}$ and $F_i \coloneqq F_{\alpha_i}$ for $i \in I$ satisfy the relations
\begin{alignat*}{3}
[\co\alpha_i, \co\alpha_j] &= 0, &\qquad\qquad
&[\co\alpha_i, E_j] = a_{ij} E_j, &\qquad\qquad
&[\co\alpha_i, F_j] = - a_{ij} F_j,\\
[E_i, F_j] &= \delta_{ij} \co\alpha_i, &\qquad\qquad
&\ad_{E_i}^{1 - a_{ij}} E_j = 0, &\qquad\qquad
&\ad_{F_i}^{1 - a_{ij}} F_j = 0,
\end{alignat*}
for all $i, j \in I$, where $A \coloneqq (a_{ij})_{i, j \in I} \coloneqq \big( \langle \alpha_j, \co\alpha_i \rangle \big)_{i, j \in I}$ is the Cartan matrix.

We fix a non-degenerate invariant symmetric bilinear form $(\cdot | \cdot)$ on $\g$. Its restriction to the Cartan subalgebra $\h$ is non-degenerate, and hence induces an isomorphism $\h \SimTo \h^\ast$. We use this to define a non-degenerate bilinear pairing on $\h^\ast$ which we also denote $(\cdot | \cdot)$. For every $i \in I$ and $\lambda \in \h^\ast$ we then have $\langle \lambda, \co\alpha_i \rangle = 2 (\lambda | \alpha_i) / (\alpha_i | \alpha_i)$. A coweight $\co\lambda \in \h$ is said to be \emph{integral} if $\langle \alpha_i, \co\lambda \rangle \in \Z$ for every $i \in I$ and \emph{dominant} if $\langle \alpha_i, \co\lambda \rangle \geq 0$ for every $i \in I$.

Let $\up \in \Aut \g$ be an automorphism of $\g$. One can always find a Cartan subalgebra $\h$ of $\g$ such that $\up$ stabilises the associated Cartan decomposition \eqref{Cartan decomp}. We will say that such a Cartan subalgebra is \emph{adapted} to the automorphism $\up$. The restriction of $\up$ to $\h$ then induces a linear map $\up \in \GL(\h)$. We define a linear map $\up \in \GL(\h^*)$ by requiring the $\up$-invariance of the bilinear pairing between $\h$ and $\h^\ast$, \emph{i.e.} $\langle \up \lambda, \up \co\mu \rangle = \langle \lambda, \co\mu \rangle$ for every $\lambda \in \h^\ast$ and $\co\mu \in \h$. This linear map leaves invariant the set of all roots $\Phi \subset \h^\ast$, the subset of positive and negative roots $\Phi^\pm$ and the subset of simple roots $\{ \alpha_i \}_{i\in I}$. It therefore defines an automorphism $\up:\Phi\to\Phi$ of the set of roots and a diagram automorphism $\up:I\to I$. The action of the automorphism $\up \in \Aut \g$ on the Cartan-Weyl basis associated with the Cartan decomposition \eqref{Cartan decomp} is then given by
\begin{equation} \label{AdaptedBasis}
\up(E_\alpha)=\tau_\alpha E_{\up(\alpha)}, \qquad
\up(\co \alpha_i) = \co \alpha_{\up(i)}, \qquad
\up(F_\alpha)=\tau_\alpha^{-1}F_{\up(\alpha)},
\end{equation}
for each $i \in I$ and $\alpha \in \Phi$, where $\tau_\alpha$ are complex numbers. The automorphism is entirely specified by its diagram part $\up : I\to I$ and the collection of numbers $\tau_{\alpha_i}$ corresponding to the simple roots.  

The bilinear form $(\cdot | \cdot)$ on $\g$ is $\up$-invariant, \emph{i.e.} $(\up X | \up Y) = (X | Y)$ for all $X, Y \in \g$. Indeed, it can be expressed as a linear combination of the Killing forms $(X, Y) \mapsto \tr_\g(\ad_X \circ \ad_Y)$ on each simple factor of $\g$ and for every $X \in \g$ we have $\ad_{\up X} = \up \circ \ad_X \circ \up^{-1}$. If the Cartan subalgebra $\h$ is adapted to $\up$ then the restriction of $(\cdot | \cdot)$ to $\h$ and the induced bilinear form on $\h^\ast$ are both $\up$-invariant.

\subsection{Adjoint group}
\label{AdGroup}

Let $G = (\Aut \g)^\circ$ be the adjoint group associated with $\g$, \emph{i.e.} the connected component of the identity in the automorphism group $\Aut \g$. It is a semisimple algebraic group, with Lie algebra $\Lie(G) = \Lie(\Aut \g) = \Der \g = \ad_\g$ isomorphic to $\g$. It is generated by the inner automorphisms of the form $\exp(\ad_X)$, for all $\ad$-nilpotent $X \in \g$.

By construction, $G$ acts on the Lie algebra $\g$ by the adjoint action. We define $H \coloneqq Z_G(\h)$, the centraliser of the Cartan subalgebra $\h$, \emph{i.e.} the subgroup of $G$ that fixes $\h$ pointwise. It is the unique maximal torus of $G$ with Lie algebra $\h$, and is a closed connected abelian subgroup of $G$. In the same way, let $B \coloneqq N_G(\b)$ and $B_- \coloneqq N_G(\b_-)$ be the normalisers of the Borel subalgebras $\b$ and $\b_-$, \emph{i.e.} the subgroups of $G$ whose adjoint action stabilises $\b$ and $\b_-$. They form Borel subgroups of $G$ (maximal closed connected solvable subgroups) and have Lie algebras $\b$ and $\b_-$.

We define the respective derived subgroups $N \coloneqq (B,B)$ and $N_- \coloneqq (B_-,B_-)$. They are closed connected unipotent subgroups of $G$ with Lie algebras $\n$ and $\n_-$. The Borel subgroups $B$ and $B_-$ then uniquely factorise as $B=HN$ and $B_-=HN_-$. We have the root decomposition
\begin{equation*}
N= \prod_{\alpha \in \Phi^+} G_\alpha \; \; \; \; \text{and} \; \; \; \; N_- = \prod_{\alpha \in \Phi^-} G_\alpha,
\end{equation*}
where, for any $\alpha \in \Phi^+$, $G_\alpha \coloneqq \left\{ \exp\left( x \ad_{E_\alpha} \right) \right\}_{x\in\C}$ and $G_{-\alpha} \coloneqq \left\{ \exp\left( x \ad_{F_\alpha} \right) \right\}_{x\in\C}$ are one-dimensional closed connected subgroups of $G$ with Lie algebras $\ad_{\g_\alpha}$ and $\ad_{\g_{-\alpha}}$. The above decompositions hold for any ordering of the roots $\alpha\in\Phi^\pm$.

Let $\upsilon \in \Aut \g$ be any given automorphism of $\g$. Since $G=(\Aut\g)^\circ$ is normal in $\Aut \g$, we can lift $\upsilon$ to an automorphism of $G$, acting by conjugation $\eta \mapsto \upsilon \circ \eta \circ \upsilon^{-1}$, which by abuse of notation we shall also denote $\upsilon \in \Aut G$. The induced Lie algebra automorphism of $\Lie(G)=\ad_\g$ coincides with $\upsilon\in\Aut\g$ via the isomorphism between $\ad_\g$ and $\g$.
If the Cartan subalgebra $\h$ is chosen to be adapted to $\up \in \Aut\g$ (cf. \S\ref{SSLieAlgebras}), then the lift $\upsilon \in \Aut G$ stabilises the corresponding subgroups $H$, $N$ and $N_-$ of $G$. As a consequence, it also stabilises the Borel subgroups $B=HN$ and $B_-=HN_-$.

\subsection{Weyl group and flag varieties}
\label{WeylFlag}

The Weyl group $W \subset \GL(\h^\ast)$ is generated by reflections $\lambda \mapsto s_i \lambda \coloneqq \lambda - \langle \lambda, \co\alpha_i \rangle \alpha_i$ for all $i \in I$. We define an action of $W$ on $\h$ by requiring the invariance of the pairing between $\h$ and $\h^\ast$, namely $\langle w(\lambda), w(\co\mu) \rangle = \langle \lambda, \co\mu \rangle$ for any $\lambda \in \h^\ast$, $\co\mu \in \h$ and $w \in W$. Explicitly, for every $i \in I$ and $\co\mu \in \h$ we have $s_i \co\mu = \co\mu - \langle \alpha_i, \co\mu \rangle \co\alpha_i$.

Let $N_G(\h)$ be the normaliser of $\h$ for the adjoint action of $G$ on the Lie algebra $\g$. The restriction of the coadjoint action $G \to \GL(\g^\ast)$ to $N_G(\h)$ gives rise to a map $N_G(\h) \to \GL(\h^\ast)$ which induces an isomorphism $\pi^{-1} : N_G(\h) / H \SimTo W$. Given any $w \in W$ we fix a representative $\dot{w} \in N_G(\h)$ of the class $\pi(w) \in N_G(\h) / H$ (the other representatives are then of the form $\dot{w}h$ with $h\in H$). For any $\alpha \in \Phi$, the action of $\dot{w}$ by conjugation on the root  subgroup $G_\alpha$ is simply $\dot{w}G_{\alpha}\dot{w}^{-1} = G_{w\alpha}$.

Consider the flag variety $G/B_-$ associated with the group $G$. The Gauss decomposition of $G$ into a disjoint union of cells $N \dot{w} B_-$ over $w \in W$ gives rise to the following cell decomposition of the flag variety
\begin{equation}\label{Schubert}
G/B_- = \bigsqcup_{w \in W} N \dot{w} B_-/B_- =: \bigsqcup_{w \in W} C_w.
\end{equation}
Since the representative $\dot{w}$ of $w\in W$ in $N_G(\h)$ only differs from other choices of representatives through right multiplication by an element of the Cartan subgroup $H \subset B_-$, it is clear that the cell $C_w = N \dot{w} B_-/B_-$ does not depend on the choice of this representative. The big cell $C_\Id = N B_-/B_-$, which is isomorphic to $N$, is dense in $G/B_-$.

Let $\up \in \Aut\g$ be an automorphism of the Lie algebra $\g$. Choose a Cartan subalgebra $\h$ adapted to $\up$. We use the corresponding action of $\up$ on $\h^\ast$ to define an action of $\up$ on $W$ by conjugation $w \mapsto \up(w) \coloneqq \up \circ w \circ \up^{-1}$. By definition of the reflection $s_i \in W$ we have in particular $\up(s_i) = s_{\up(i)}$ for every $i \in I$. The map $\upsilon:W\rightarrow W$ constructed in this way is an automorphism of the group $W$. Note that it only depends on the restriction of $\upsilon$ to $\h$, \emph{i.e.} on the associated diagram automorphism $\up: I \rightarrow I$. We introduce the subgroup
\begin{equation*}
W^\upsilon \coloneqq \{ w \in W \,|\, \upsilon(w) = w \}
\end{equation*}
of $\upsilon$-invariant elements in $W$.

The lift $\upsilon \in \Aut G$ of the automorphism $\up \in \Aut \g$ to the adjoint group $G$ stabilises both the torus $H$ and the normaliser $N_G(\h)$. It therefore descends to an automorphism $\up$ of the subquotient $N_G(\h)/H$. The isomorphism $\pi:W \SimTo N_G(\h)/H$ is equivariant with respect to the action of $\up$ on $W$ and $N_G(\h) / H$, \emph{i.e.} $\pi\bigl(\upsilon(w)\bigr)=\upsilon\bigl(\pi(w)\bigr)$ for any $w\in W$.
Since $\upsilon$ stabilises the Borel subgroup $B_-$, one can also define an automorphism $\up$ of the flag variety $G/B_-$. This will play a central role in \S\ref{ReproFlag} so we postpone its definition and the study of its properties until then.

\section{Finite $\g$-opers and finite Miura $\g$-opers}

We follow the conventions and notations of \S\ref{Conventions}.

\subsection{Principal $\sl_2$ and $\Z$-grading}
\label{Principal}

Consider the regular nilpotent element
\begin{equation} \label{pm1 def}
p_{-1} \coloneqq \sum_{i \in I} F_i.
\end{equation}
By the Jacobson-Morosov theorem it can be embedded into an $\sl_2$-triple.
Let $\co \omega_i \in \h$, $i \in I$ be the fundamental coweights of $\g$ defined by $\alpha_i(\co \omega_j) = \delta_{ij}$ for all $i, j \in I$. The Weyl covector
\begin{equation} \label{corho def}
\co \rho \coloneqq \sum_{i \in I} \co \omega_i
\end{equation}
then satisfies $\alpha_i(\co \rho) = 1$ for all $i \in I$, so that $[2 \co\rho, p_{-1}] = - 2 p_{-1}$. The pair $\{ p_{-1}, 2 \co\rho \}$ extends uniquely to an $\sl_2$-triple $\{ p_{-1}, 2 \co \rho, p_1 \}$, \emph{i.e.} with the relations
\begin{equation*}
[p_1, p_{-1}] = 2 \co \rho, \qquad
[\co \rho, p_{\pm 1}] = \pm p_{\pm 1}.
\end{equation*}
Specifically, the regular nilpotent element $p_1$ is obtained by writing $2 \co \rho = \sum_{i \in I} c_i \co \alpha_i$ in the basis of simple coroots $\co \alpha_i$ of $\g$ and defining $p_1 \coloneqq \sum_{i \in I} c_i  E_i$.
We will refer to $\{ p_{-1}, 2 \co \rho, p_1 \}$ as the \emph{principal $\sl_2$-triple} of $\g$ and to its span $\langle p_{-1}, \co \rho, p_1 \rangle \subset \g$ as the \emph{principal $\sl_2$ subalgebra} of $\g$.

The \emph{height} $\hgt(\alpha) \in \Z$ of a root $\alpha \in \Phi$ is the eigenvalue of the corresponding root vector under the adjoint action of the Cartan element $\co \rho$, \emph{i.e.} $[\co \rho, E_{\alpha}] = \hgt(\alpha) E_{\alpha}$ and $[\co \rho, F_{\alpha}] = \hgt(-\alpha) F_{\alpha}$ for any $\alpha \in \Phi^+$. The eigenspace decomposition of $\ad_{\co \rho}$ defines a $\Z$-grading of the Lie algebra $\g$,
\begin{equation} \label{Z-grading}
\g = \bigoplus_{i \in \Z} \g_i = \bigoplus_{i=-h+1}^{h-1} \g_i,
\end{equation}
with $[ \g_i, \g_j ] \subset \g_{i+j}$ for all $i, j \in \Z$. Here $h$ denotes the Coxeter number of $\g$, defined as the maximum of the Coxeter numbers of the simple factors in the semi-simple decomposition of $\g$ (for simple $\g$ we have $h \coloneqq \hgt(\theta) + 1$ with $\theta$ the maximal root in $\Phi$). The eigenspace $\g_j \coloneqq \{ X \in \g \,|\, [\co \rho, X] = j X \}$ of $\ad_{\co \rho}$ for each $j \in \Z$ is spanned by all root vectors of height $j$, explicitly
\begin{equation*}
\g_0 = \h, \qquad \g_i = \bigoplus_{\alpha \in \Phi^+ \, | \, \hgt(\alpha) = i} \g_\alpha, \qquad \g_{-i} = \bigoplus_{\alpha \in \Phi^+ \, | \, \hgt(\alpha) = i} \g_{-\alpha}
\end{equation*}
for all $i \in \Z_{\geq 1}$. In particular we have $\g_i = 0$ if $|i| \geq h$. Note also that $\b = \oplus_{i \geq 0} \g_i$ and $\n = \oplus_{i > 0} \g_i$.

The centraliser of $p_1$ in $\g$ is a $\Z_{\geq 1}$-graded abelian subalgebra $\a \coloneqq \ker(\ad_{p_1})$ of dimension $\rk \g$.
Let $\Vc{j} \coloneqq \a \cap \g_j$ denote its grade-$j$ subspace for each $j \in \{ 1, \ldots, h-1 \}$. We let $E$ denote the multiset containing each $j \in \{ 1, \ldots, h-1 \}$ with multiplicity $\dim \Vc{j}$. Its elements are called the exponents of $\g$.
From the representation theory of $\sl_2$ we then have, for every $i \in \{ 1, \ldots, h-1 \}$,
\begin{equation} \label{gi decomp Vc}
\g_i = \left\{ \begin{array}{ll} 
[p_{-1}, \g_{i+1}] \qquad & \text{if}\quad i \not\in E,\\
{[p_{-1}, \g_{i+1}]} \oplus \Vc{i} \qquad & \text{if}\quad i \in E.
\end{array} \right.
\end{equation}

\subsection{Finite $\g$-opers}
\label{finite opers}

Let $G$ be the adjoint group of $\g$, cf. \S\ref{AdGroup}. If $g \in G$ and $X \in \g$ then we denote by $g X g^{-1}$ the adjoint action of $g$ on $X$. The affine subspace $p_{-1} + \b \subset \g$ is stabilised by the adjoint action of the unipotent subgroup $N$ of $G$. Consider the quotient
\begin{equation*}
\Op^{\rm fin}_\g \coloneqq (p_{-1} + \b)/ N
\end{equation*}
and denote the class of any $X \in p_{-1} + \b$ as $[X]_\g$. We refer to the elements of
$\Op^{\rm fin}_\g$ as \emph{finite $\g$-opers}.
The so called \emph{Slodowy slice} through $p_{-1}$ is defined as the affine subspace $p_{-1} + \a \subset p_{-1} + \b$.
It is transversal at every point to the adjoint orbit of $N$ in $p_{-1} + \b$. In fact, it is an important result of Kostant that the Slodowy slice intersects each $N$-orbit in $p_{-1} + \b$ exactly once. For our purpose it is convenient to formulate this statement as follows.

\begin{theorem} \label{prop: mini opers}
Every finite $\g$-oper has a unique representative in $p_{-1} + \a$. We shall refer to it as the \emph{canonical representative} of the finite $\g$-oper.
\begin{proof}
This follows from \cite[Theorem 1.2]{KostantWhittaker} which states that the map $N \times (p_{-1} + \a) \rightarrow p_{-1} + \b$, given by the adjoint action $(n, X) \mapsto n X n^{-1}$ is an isomorphism of affine varieties. Alternatively, it could also be proved by a simpler version of the argument used in the proof of Theorem \ref{TwistedDS} below. However, since the argument is so similar to the one given there we do not repeat it here for the sake of brevity.
\end{proof}
\end{theorem}

We define a \emph{finite Miura $\g$-oper} as an element of the affine subspace $p_{-1} + \h \subset \g$. To any finite Miura $\g$-oper is associated a finite $\g$-oper, namely its class in $\Op^{\rm fin}_\g$. It follows from the second half of the proof of \cite[Theorem 1.2]{KostantWhittaker} that the corresponding map
\begin{equation} \label{MOpfin to Opfin}
\begin{split}
p_{-1} + \h &\longrightarrow \Op_\g^{\rm fin}, \\
X &\longmapsto [X]_\g.
\end{split}
\end{equation}
is surjective.
The next theorem gives a necessary and sufficient condition for two finite Miura $\g$-opers to correspond to the same finite $\g$-oper.

\begin{proposition} \label{prop: mini miura opers}
For any $\co\lambda, \co\mu \in \h$, we have $[p_{-1} - \co\lambda]_\g = [p_{-1} - \co\mu]_\g$ if and only if there exists $w \in W$ such that $\co\lambda = w(\co\mu)$.
\begin{proof}
Let $\co\lambda, \co\mu \in \h$. Denote by $\bm u = (u_1, \ldots, u_{\rk \g}) : \g \to \C^{\rk \g}$ the collection of $\rk \g$ fundamental homogenous $G$-invariant polynomials on $\g$.
By \cite[Lemma 9.2]{KostantPrincipal}, there exists $w \in W$ such that $\co\lambda = w(\co\mu)$ if and only if $\bm u(\co\lambda) = \bm u(\co\mu)$. The latter is equivalent to $\bm u(p_{-1} - \co\lambda) = \bm u(p_{-1} - \co\mu)$ by \cite[Proposition 17]{KostantPolynomial}. Now since $p_{-1} - \co\lambda$ and $p_{-1} - \co\mu$ are regular by \cite[Lemma 10]{KostantPolynomial}, this in turn is equivalent to $p_{-1} - \co\lambda$ and $p_{-1} - \co\mu$ being conjugate under the adjoint action of $G$ by \cite[Theorem 3]{KostantPolynomial}. It is clear that if $p_{-1} - \co\lambda$ and $p_{-1} - \co\mu$ are conjugate under the adjoint action of $N$ then they are also conjugate under the adjoint action of $G$. The converse of this statement follows from the second half of the proof of \cite[Theorem 1.2]{KostantWhittaker}. Therefore $p_{-1} - \co\lambda$ and $p_{-1} - \co\mu$ are conjugate under the adjoint action of $G$ if and only if $[p_{-1} - \co\lambda]_\g = [p_{-1} - \co\mu]_\g$, and the theorem follows.
\end{proof}
\end{proposition}

We define the \emph{shifted action} of $W$ on $\h$ by letting $w \in W$ send $\co\lambda \in \h$ to $w \cdot \co\lambda \coloneqq w(\co\lambda + \co\rho) - \co\rho$. Two coweights $\co\lambda, \co\mu \in \h$ are then said to be \emph{$W$-linked} if $\co\mu = w \cdot \co\lambda$ for some $w \in W$. The \emph{$W$-linkage class} of a coweight $\co\lambda \in \h$ is its orbit under this shifted action of $W$, which we denote by $[\co\lambda]_W$.

\begin{corollary}
There is a bijection between the set $\h / (W, \cdot)$ of $W$-linkage classes in $\h$ and the set $\Op_\g^{\rm fin}$ of finite $\g$-opers, given explicitly by
$[\co\lambda]_W \mapsto [p_{-1} - \co\lambda - \co\rho]_\g$.
\begin{proof}
The given map is well defined and injective by Proposition \ref{prop: mini miura opers}. It is surjective since \eqref{MOpfin to Opfin} is.
\end{proof}
\end{corollary}

By abuse of notation, in what follows we always identify the $W$-linkage class $[\co\lambda]_W$ of a coweight $\co\lambda \in \h$ with the corresponding finite $\g$-oper $[p_{-1} - \co\lambda - \co\rho]_\g \in \Op_\g^{\rm fin}$.
Let $\co\lambda \in \h$ and consider the subset
$\MOp_{\g, [\co\lambda]_W}^{\rm fin} \coloneqq \big\{ X \in p_{-1} + \h \,\big|\, [X]_\g = [\co\lambda]_W \big\}$
of the set $p_{-1} + \h$ of all finite Miura $\g$-opers.

\begin{theorem} \label{MOp fin iso W}
If $\co\lambda \in \h$ is dominant then $\MOp_{\g, [\co\lambda]_W}^{\rm fin} \simeq W$.
\begin{proof}
By Proposition \ref{prop: mini miura opers}, every $X \in \MOp_{\g, [\co\lambda]_W}^{\rm fin}$ is of the form $X = p_{-1} - w(\co\lambda + \co\rho)$ for some $w \in W$. We obtain a surjective map $W \twoheadrightarrow \MOp_{\g, [\co\lambda]_W}^{\rm fin}$ given by $w \mapsto p_{-1} - w(\co\lambda + \co\rho)$. To prove injectivity, it is enough to show that if $\co\lambda + \co\rho$ is fixed by $w \in W$ then $w = \Id$.
But since $\co\lambda$ is assumed dominant, we have $\langle \alpha_i, \co\lambda + \co\rho \rangle \geq 1$ for all $i \in I$. In other words, the coweight $\co\lambda + \co\rho$ lies inside the open fundamental Weyl chamber of $\h$. The result now follows from the fact that $W$ acts simply transitively on the set of all Weyl chambers.
\end{proof}
\end{theorem}

\subsection{Folding and finite $\g^\nu$-opers}
\label{folding}

Let $\nu: I \to I$ be a permutation of the simple roots of $\g$ which preserves the Cartan matrix, \emph{i.e.} such that
\begin{equation} \label{diag aut def}
a_{\nu(i) \nu(j)} = a_{ij}
\end{equation}
for all $i, j \in I$. Consider the associated diagram automorphism of $\g$, which we also denote $\nu \in \Aut \g$, defined by its action on the Chevalley-Serre generators as
\begin{equation}\label{Nu}
\nu(E_i) = E_{\nu(i)}, \qquad
\nu(\co \alpha_i) = \co \alpha_{\nu(i)}, \qquad
\nu(F_i) = F_{\nu(i)}.
\end{equation}
We are interested in describing the $\nu$-invariant subalgebra $\g^\nu \coloneqq \{ X \in \g \,|\, \nu(X) = X \}$.

Let $I / \nu$ denote the set of orbits in $I$ under $\nu$.
For each orbit $\mathcal I \in I / \nu$ we define (cf. \cite{FSS})
\begin{equation*}
\ell_{\mathcal I} \coloneqq 3 - \sum_{i \in \mathcal I} a_{ij},
\end{equation*}
where $j$ is any point in the orbit $\mathcal I$. The right hand side depends only on the orbit $\mathcal I$ itself, and not on $j$, by virtue of the assumption \eqref{diag aut def} on $\nu$. Moreover, by the properties $a_{ii} = 2$ and $a_{ij} \in \Z_{\leq 0}$ for all $i \neq j$ of the Cartan matrix it follows that $\ell_{\mathcal I} \in \Z_{\geq 1}$. Since $\g$ is semisimple we have $\ell_{\mathcal I} = 1$ or $2$ for every $\mathcal I \in I / \nu$ (see \emph{e.g.} \cite[Lemma 4.3]{Naito}). If $\ell_{\mathcal I} = 1$ then the restriction of the Dynkin diagram to the orbit $\mathcal I$ is isomorphic to the Dynkin diagram of $A_1^{\times |\mathcal I |}$. On the other hand, if $\ell_{\mathcal I} = 2$ then $|\mathcal I |$ is even and the Dynkin diagram restricted to the orbit $\mathcal I$ is isomorphic to the Dynkin diagram of $A_2^{\times |\mathcal I | / 2}$.
Finally, note that in either case the relation $\ell_{\mathcal I} \sum_{i \in \mathcal I} a_{ij} = 2$ holds for any $j \in \mathcal I$.

The following result is well known (see for instance \cite[\S 9.5]{CarterBook}).

\begin{proposition} \label{g nu semisimple}
The $\nu$-invariant subalgebra $\g^\nu$ is semisimple. Its Chevalley-Serre generators are
\begin{equation*}
\co\alpha^\nu_{\mathcal I} \coloneqq \ell_{\mathcal I} \sum_{i \in \mathcal I} \co\alpha_i, \qquad
E^\nu_{\mathcal I} \coloneqq \ell_{\mathcal I} \sum_{i \in \mathcal I} E_i, \qquad
F^\nu_{\mathcal I} \coloneqq \sum_{i \in \mathcal I} F_i,
\end{equation*}
for $\mathcal I \in I/\nu$ with corresponding system of simple roots $\alpha^\nu_{\mathcal I} \coloneqq \frac{1}{|\mathcal I |} \sum_{i \in \mathcal I} \alpha_i \subset \h^{\ast,\nu} \coloneqq (\h^\ast)^\nu$. The Cartan matrix $A^\nu \coloneqq (a_{\mathcal I \mathcal J}^\nu)_{\mathcal I, \mathcal J \in I/\nu}$
is given by
$a_{\mathcal I \mathcal J}^\nu \coloneqq \langle \alpha^\nu_{\mathcal J}, \co\alpha^\nu_{\mathcal I} \rangle = \ell_{\mathcal I} \sum_{i \in \mathcal I} a_{ij}$ for any $j \in \mathcal J$,
and the Weyl group is the $\nu$-invariant subgroup
$W^\nu$, cf. \S\ref{WeylFlag}.
\end{proposition}

We shall need the following explicit description of the Weyl group $W^\nu$ of $\g^\nu$ (see \emph{e.g.} \cite{CharlesVar}).

\begin{lemma} \label{lem: Weyl nu}
The simple reflections $\{ s^\nu_{\mathcal I} \}_{\mathcal I \in I/\nu}$ generating the $\nu$-invariant subgroup $W^\nu \subset W$ read 
\begin{equation*}
s^\nu_{\mathcal I} = \left\{ \begin{array}{ll}
\prod_{i \in \mathcal I} s_i & \;\;\text{for} \quad \ell_{\mathcal I} = 1,\\
\prod_{i \in \mathcal I / 2} s_i s_{\bar\imath} s_i & \;\;\text{for} \quad \ell_{\mathcal I} = 2,
\end{array}
\right.
\end{equation*}
where the product in the case $\ell_{\mathcal I} = 2$ is over half of the orbit, namely $\mathcal I / 2 \coloneqq \{ k, \nu(k), \ldots, \nu^{|\mathcal I|/2 - 1}(k) \}$ for any $k \in \mathcal I$, and we define $\bar\imath \coloneqq \nu^{|\mathcal I|/2}(i)$ for every $i \in \mathcal I$.
\begin{proof}
Let $\mathcal I \in I/\nu$. We determine $s^\nu_{\mathcal I}$ by computing its action on an arbitrary $\nu$-invariant coweight $\co\mu \in \h^\nu$. By definition we have
\begin{equation} \label{si orbit}
s^\nu_{\mathcal I} \co\mu = \co\mu - \langle \alpha^\nu_{\mathcal I}, \co\mu \rangle \co\alpha^\nu_{\mathcal I}
= \co\mu - \frac{\ell_{\mathcal I}}{|\mathcal I|} \sum_{i, j \in \mathcal I} \langle \alpha_j, \co\mu \rangle \co\alpha_i
= \co\mu - \ell_{\mathcal I} \sum_{i \in \mathcal I} \langle \alpha_i, \co\mu \rangle \co\alpha_i,
\end{equation}
where in the last step we used the relation $\langle \alpha_j, \co\mu \rangle = \langle \alpha_i, \co\mu \rangle$ which follows from the $\nu$-invariance of $\co\mu$ and the fact that $i, j$ lie on the same orbit $\mathcal I$.

If $\ell_{\mathcal I} = 1$ then the right hand side of \eqref{si orbit} can be rewritten as $\big( \prod_{i \in \mathcal I} s_i \big) (\co\mu)$ since $\langle \alpha_j, \co\alpha_i \rangle = 0$ for all distinct $i, j \in \mathcal I$. On the other hand, if $\ell_{\mathcal I} = 2$ then a direct computation using the fact that $\langle \alpha_i, \co\alpha_{\bar\imath} \rangle = \langle \alpha_{\bar\imath}, \co\alpha_i \rangle = -1$ yields $s_i s_{\bar\imath} s_i (\co\mu) = \co\mu - 2 \langle \alpha_i, \co\mu \rangle \co\alpha_i - 2 \langle \alpha_{\bar\imath}, \co\mu \rangle \co\alpha_{\bar\imath}$. The result now follows.
\end{proof}
\end{lemma}

The fundamental coweights of $\g^\nu$, defined by $\alpha^\nu_{\mathcal I}(\co\omega^\nu_{\mathcal J}) = \delta_{\mathcal I \mathcal J}$ for all $\mathcal I, \mathcal J \in I / \nu$, are
\begin{equation*}
\co\omega^\nu_{\mathcal I} \coloneqq \sum_{i \in \mathcal I} \co\omega_i.
\end{equation*}

\begin{proposition}\label{ExpoNu}
The automorphism $\nu$ preserves the $\Z$-grading \eqref{Z-grading} and fixes the principal $\sl_2$. In particular, $\{ p_{-1}, 2 \co\rho, p_1 \}$ is the principal $\sl_2$-triple of $\g^\nu$.
Moreover, each $\Vc{i}$ for $i \in E$, is $\nu$-stable.
\begin{proof}
Since $\hgt(\nu(\alpha)) = \hgt(\alpha)$ for any $\alpha \in \Phi$, it follows that $\nu$ stabilises the eigenspaces $\g_i$, $i \in \Z$.
Moreover, $\nu(\co\omega_i) = \co\omega_{\nu(i)}$ for all $i \in I$ since $\langle \alpha_j, \nu(\co\omega_i) \rangle = \langle \alpha_{\nu^{-1}(j)}, \co\omega_i \rangle = \delta_{\nu^{-1}(j) i} = \delta_{j \nu(i)} = \langle \alpha_j, \co\omega_{\nu(i)} \rangle$ for any $j \in I$. Therefore $\co\rho$ is $\nu$-invariant.
By definition of $p_{-1}$ and $\nu$, we clearly have $\nu(p_{-1}) = p_{-1}$. It follows that
\begin{equation*}
[p_1,p_{-1}] = 2 \co \rho = 2 \nu(\co \rho) = \nu\big( [p_1,p_{-1}] \big) = [ \nu(p_1), \nu(p_{-1}) ] = [ \nu(p_1), p_{-1} ],
\end{equation*}
so $\ad_{p_{-1}}(p_1 - \nu(p_1)) = 0$. As $\nu$ stabilises $\g_1$, $p_1$ and $\nu(p_1)$ are both in $\g_1$. Yet $\ad_{p_{-1}}: \g_{i+1} \rightarrow \g_i$ is injective for $i\geq 0$, and hence $p_1 = \nu(p_1)$.

Recall the expressions \eqref{pm1 def} and \eqref{corho def} for $p_{-1}$ and $\co\rho$, respectively. Noting that the sum over $i \in I$ in these expressions can be rewritten as a double sum over orbits $\mathcal I \in I / \nu$ and elements $i \in \mathcal I$ in each orbit, we can write these as
\begin{equation*}
p_{-1} = \sum_{\mathcal I \in I / \nu} F^\nu_{\mathcal I}, \qquad
\co\rho = \sum_{\mathcal I \in I / \nu} \co\omega^\nu_{\mathcal I}.
\end{equation*}
It follows that $\{ p_{-1}, 2 \co\rho, p_1 \}$ is the principal $\sl_2$-triple of $\g^\nu$.

Let $X \in \a$. Then $[\nu(X), p_1] = [\nu(X), \nu(p_1)] = \nu\big( [X, p_1] \big) = 0$ so that $\nu(X) \in \a$. Thus $\nu$ stabilises both $\a$ and $\g_i$, and hence also $\Vc{i} = \a \cap \g_i$.
\end{proof}
\end{proposition}

As a result of Proposition \ref{g nu semisimple} and Proposition \ref{ExpoNu}, the constructions of \S\ref{Principal} and \S\ref{finite opers} apply directly to the $\nu$-invariant subalgebra $\g^\nu$. In particular, $\a^\nu \coloneqq \a \cap \g^\nu = \ker(\ad_{p_1} : \g^\nu \to \g^\nu)$ is a $\Z_{\geq 1}$-graded abelian subalgebra in $\g^\nu$ of dimension $\rk \g^\nu$, and the corresponding multiset $E^\nu$ of exponents of $\g^\nu$ forms a sub-multiset of $E$.

Since $\nu$ fixes $p_{-1}$ and stabilises $\b$, we can consider the affine subspace $p_{-1}+\b^\nu \subset \g^\nu$ of $\nu$-invariant elements in $p_{-1}+\b$. It is stabilised by the adjoint action of $N^\nu$, the $\nu$-invariant subgroup of $N$. We may therefore define the corresponding space of finite $\g^\nu$-opers
\begin{equation*}
\Op^{\rm fin}_{\g^\nu} = (p_{-1}+\b^\nu)/N^\nu.
\end{equation*}
If $X \in p_{-1}+\b^\nu$ then we denote by $[X]_{\g^\nu}$ the associated finite $\g^\nu$-oper. In particular, we recover the notion of finite $\g$-opers from \S\ref{finite opers} when $\nu=\Id$. For completeness let us end by stating the analogs of Theorems \ref{prop: mini opers}, Proposition \ref{prop: mini miura opers} and Theorem \ref{MOp fin iso W} for finite $\g^\nu$-opers.

\begin{theorem} \label{prop: mini miura opers cyc}
Every finite $\g^\nu$-oper has a unique representative in $p_{-1} + \a^\nu$.

For any $\nu$-invariant coweights $\co\lambda, \co\mu \in \h^\nu$, we have $[p_{-1} - \co\lambda]_{\g^\nu} = [p_{-1} - \co\mu]_{\g^\nu}$ if and only if there exists $w \in W^\nu$ such that $\co\lambda = w(\co\mu)$.

If $\co\lambda \in \h^\nu$ is dominant then $\MOp_{\g^\nu, [\co\lambda]_{W^\nu}}^{\rm fin} \simeq W^\nu$.
\end{theorem}

\section{Cyclotomic $\g$-opers and canonical representatives}
\label{DS}

We pick and fix, once and for all, a diagram automorphism $\nu : I \to I$ as in \S\ref{folding}. Let $T \in \Z_{\geq 1}$ be a multiple of the order of $\nu$ and choose a primitive $T^{\rm th}$-root of unity $\omega$. It generates a copy of the cyclic group of order $T$ under multiplication, which we denote
\begin{equation*}
\Gamma \coloneqq \langle\omega\rangle \subset \C^\times.
\end{equation*}

In this section we will define an action of $\Gamma$ on various spaces defined over the Riemann sphere, including the spaces of $\g$-valued meromorphic functions, differentials and connections. Elements of these spaces which are invariant under this action will be called \emph{cyclotomic}.

\subsection{Cyclotomic $\g$-valued functions and differentials}
\label{cycl func diff}

Let $\CP \coloneqq \C \cup \{ \infty \}$ be the Riemann sphere and fix a global coordinate $t$ on $\C \subset \CP$.
There is a natural action $\mu : \Gamma \to \Aut \CP$, $\alpha \mapsto \mu_\alpha$ of the cyclic group $\Gamma$ on $\CP$ which fixes $\infty$ and with $\alpha \in \Gamma$ acting by multiplication $t \mapsto \alpha t$ on $\C$.

In this subsection and the next, we let $\up \in \Aut \g$ be any automorphism whose order divides $T$, \emph{i.e.} such that $\up^T = \Id$, and with diagram part $\nu$. In other words, let $\up$ be defined on the Chevalley-Serre generators as
\begin{equation} \label{upsilon order T}
\up(E_i) = \tau_{\alpha_i} E_{\nu(i)}, \qquad
\up(\co \alpha_i) = \co \alpha_{\nu(i)}, \qquad
\up(F_i) = \tau_{\alpha_i}^{-1}F_{\nu(i)},
\end{equation}
where $\tau_{\alpha_i}$, $i \in I$ are arbitrary $T^{\rm th}$-roots of unity.
We obtain an action of the cyclic group $\Gamma$ on $\g$ by letting $\omega$ act as $\up$. That is, we have a homomorphism $\Gamma \to \Aut \g$, $\omega \mapsto \up$.

Let $\M$ denote the algebra of meromorphic functions on $\CP$. We introduce an action of $\Gamma$ on $\M$ by letting $\alpha \in \Gamma$ act through the pullback $\mu_{\alpha^{-1}}^\ast : \M \to \M$, namely sending $f \in \M$ to the function $\alpha.f \coloneqq f \circ \mu_{\alpha^{-1}} \in \M$. Let $\Omega \coloneqq \M dt$ be the space of meromorphic differentials on $\CP$. It too comes naturally equipped with an action of $\Gamma$, letting $\alpha \in \Gamma$ act also by the pullback $\mu_{\alpha^{-1}}^\ast : \Omega \to \Omega$.

Consider the Lie algebra $\g(\M) \coloneqq \g \otimes_\C \M$, \emph{i.e.} the set of $\g$-valued meromorphic functions equipped with the pointwise Lie bracket.
Define an action of $\Gamma$ on $\g(\M)$ by combining the action of $\Gamma$ on $\g$ with the above action on $\M$. That is, for any $X \otimes f \in \g(\M)$ we set
\begin{equation} \label{hat tau action}
\hat{\up}(X \otimes f) \coloneqq \up(X) \otimes \mu_{\omega^{-1}}^\ast f.
\end{equation}
We denote the subalgebra of $\Gamma$-invariants as
\begin{equation*}
\g^{\hat\up}(\M) \coloneqq
\{ h \in \g(\M) \,|\, \hat\up (h) = h \},
\end{equation*}
which we shall call the space of cyclotomic $\g$-valued meromorphic functions.
Given any $\up$-stable subspace $\p \subset \g$ we define the $\hat\up$-stable subspace $\p(\M) \coloneqq \p \otimes_\C \M$ of $\g(\M)$. We denote the corresponding subspace of $\Gamma$-invariants by
\begin{equation*}
\p^{\hat\up}(\M) \coloneqq
\{ h \in \p(\M) \,|\, \hat\up(h) = h \}.
\end{equation*}

Given any subspace $\p \subset \g$, we denote by $\Omega(\p) \coloneqq \p \otimes_\C \Omega$ the space of $\p$-valued meromorphic differentials. If $\p$ is $\up$-stable, we equip $\Omega(\p)$ with an action of $\Gamma$ defined by letting $\omega$ act through the pullback by $\mu_{\omega^{-1}}$ on the second tensor factor, namely
\begin{equation*}
\hat\up(X \otimes \varpi) \coloneqq \up(X) \otimes \mu_{\omega^{-1}}^\ast \varpi,
\end{equation*}
for any $X \in \p$ and $\varpi \in \Omega$.
We define the subspace
\begin{equation*}
\Omega^{\hat\up}(\p) \coloneqq \{ \ms A \in \Omega(\p) \,|\, \hat\up(\ms A) = \ms A \}
\end{equation*}
of cyclotomic (\emph{i.e.} $\Gamma$-invariant) $\p$-valued meromorphic differentials.
We shall often abbreviate the notation $X \otimes f$ for an element in $\p(\M)$ to $X f$ for simplicity, and similarly an element $Y \otimes \varpi$ of $\Omega(\p)$ will be denoted simply as $Y \varpi$.

\subsection{Cyclotomic $\g$-connections}
\label{TwistedAlgebra}

Let $d : \M \to \Omega$ denote the de Rham differential and consider the affine space
\begin{equation*}
\Conn_\g(\CP) \coloneqq \{ d + \ms A \,|\, \ms A \in \Omega(\g) \}
\end{equation*}
of $\g$-valued meromorphic connections on $\CP$, or \emph{$\g$-connections} for short. More generally, given any subspace $\p \subset \g$ we define the affine subspace $\Conn_\p(\CP)$ of $\p$-valued mermorphic connections of the form $d + \ms A$ with $\ms A \in \Omega(\p)$. Similarly, we define the affine subspace
\begin{equation*}
\Conn^{\hat\up}_\g(\CP) \coloneqq \{ d + \ms A \,|\, \ms A \in \Omega^{\hat\up}(\g) \} \subset \Conn_\g(\CP)
\end{equation*}
of cyclotomic $\g$-valued meromorphic connections on $\CP$, and if $\p \subset \g$ is a $\up$-stable subspace we define the affine subspace $\Conn^{\hat\up}_\p(\CP) \subset \Conn^{\hat\up}_\g(\CP)$ of connections of the form $d + \ms A$ with $\ms A \in \Omega^{\hat\up}(\p)$. We shall refer to elements of $\Conn^{\hat\up}_\p(\CP)$ as \emph{cyclotomic $\p$-connections}.

Consider the group $G(\M)$ of $\M$-points of $G$, which can be thought of as the group of
$G$-valued meromorphic functions on $\CP$. It can be formally defined as the set of all $\C$-algebra homomorphisms $\sO(G) \to \M$, where $\sO(G)$ is the algebra of regular functions on $G$.
The group structure on $G$ makes $\sO(G)$ into a commutative Hopf algebra with coproduct
$\Delta : \sO(G) \to \sO(G \times G) \simeq \sO(G) \otimes \sO(G)$
given by $(\Delta \varphi)(x, y) = \varphi(xy)$. The counit and antipode maps are given by $\epsilon : \sO(G) \to \C$, $\epsilon(\varphi) = \varphi(e)$ and $s : \sO(G) \to \sO(G)$, $s(\varphi)(x) = \varphi(x^{-1})$ respectively.

The Lie algebra of $G(\M)$ is naturally isomorphic to $\g(\M)$. Lifting $\upsilon \in \Aut \g$ to an automorphism $\up \in \Aut G$ of $G$, cf. \S\ref{AdGroup}, this in turn induces an automorphism of $\sO(G)$ via the pullback $\upsilon^\ast$. In particular, when the order of $\up$ divides $T$ we can endow $G(\M)$ with an action of $\Gamma$. Specifically, we let $\omega$ act on $G(\M)$ by sending any $g : \sO(G) \to \M$ to the composition
\begin{equation*}
\sO(G) \overset{\up^\ast}\longrightarrow \sO(G) \overset{g}\longrightarrow \M \xrightarrow{\mu^\ast_{\omega^{-1}}} \M,
\end{equation*}
which we write $\hat\up(g) \in G(\M)$.
The map $\hat\upsilon : G(\M) \to G(\M)$ so defined is an automorphism of $G(\M)$ since the coproduct of $\sO(G)$ and the multiplication in $\M$ are both $\Gamma$-equivariant. In other words, the following diagram
\begin{equation*}
\begin{tikzpicture}
\matrix (m) [matrix of math nodes, row sep=1em, column sep=3em,text height=1.5ex, text depth=0.25ex]    
{
            & \sO(G) \otimes \sO(G) &                                     &                        & \M \otimes \M & \\
\sO(G) &                                     & \sO(G) \otimes \sO(G) & \M \otimes \M &                         & \M \\
            & \sO(G)                         &                                     &                        & \M                    &\\
};
\path[->] (m-2-1) edge node[above left]{\small $\Delta$} (m-1-2);
\path[->] (m-2-1) edge node[below left]{\small $\up^\ast$} (m-3-2);
\path[->] (m-1-2) edge node[above right]{\small $\up^\ast \otimes \up^\ast$} (m-2-3);
\path[->] (m-3-2) edge node[below right]{\small $\Delta$} (m-2-3);
\path[->] (m-2-3) edge node[below]{$g_1 \otimes g_2$} (m-2-4);
\path[->] (m-2-4) edge node[above left=-1mm]{\small $\mu_{\omega^{-1}}^\ast \otimes \mu_{\omega^{-1}}^\ast$} (m-1-5);
\path[->] (m-2-4) edge (m-3-5);
\path[->] (m-1-5) edge (m-2-6);
\path[->] (m-3-5) edge node[below right]{\small $\mu_{\omega^{-1}}^\ast$} (m-2-6);
\end{tikzpicture}
\end{equation*}
commutes for any $g_1, g_2 \in G(\M)$, with the sequence of maps from left to right along the top of the diagram representing the element $\hat\upsilon (g_1) \hat\upsilon (g_2) \in G(\M)$ and those along the bottom corresponding to $\hat\upsilon(g_1 g_2) \in G(\M)$.
The subgroup $G^{\hat\up}(\M)$ of $\Gamma$-invariant elements in $G(\M)$ then consisting of all $\Gamma$-equivariant $\C$-algebra homomorphisms $g : \sO(G) \to \M$, \emph{i.e.} for which the diagram
\begin{equation*}
\begin{tikzpicture}
\matrix (m) [matrix of math nodes, row sep=2em, column sep=3em,text height=1.5ex, text depth=0.25ex]    
{
\sO(G) & \M\\
\sO(G) & \M\\
};
\path[->] (m-1-1) edge node[above]{$g$} (m-1-2);
\path[->] (m-2-1) edge node[below]{$g$} (m-2-2);
\path[->] (m-1-1) edge node[left]{$\upsilon^\ast$} (m-2-1);
\path[->] (m-1-2) edge node[right]{$\mu^\ast_\omega$} (m-2-2);
\end{tikzpicture}
\end{equation*}

\vspace{-3mm}
\noindent commutes. Given any $t \in \CP$, by abuse of notation we write $g(t)$ for the composition $\text{ev}_t \circ g$ with the evaluation at $t$ map $\text{ev}_t : \M \to \CP$. The property of $\Gamma$-equivariance of $g \in G^{\hat\up}(\M)$ may then be expressed as $\up\big( g(t) \big) = g(\omega^{-1} t)$.
The Lie algebra of $G^{\hat\up}(\M)$ is naturally isomorphic to $\g^{\hat \up}(\M)$.

Now for any $g \in G(\M)$ we define the $\C$-linear map $dg : \sO(G) \to \Omega$ by $\varphi \mapsto d\big( g(\varphi) \big)$. In particular we then have $dg g^{-1} \in \Omega(\g)$, where the product of $dg$ with $g^{-1}$ is defined as for $G(\M)$ in terms of the coproduct on $\sO(G)$.
The affine space $\Conn_\g(\CP)$ is equipped with an action of $G(\M)$ by gauge transformations
\begin{equation} \label{gauge tr def}
d + \ms A \longmapsto g (d + \ms A) g^{-1} \coloneqq d - dg g^{-1} + \Ad_g \ms A,
\end{equation}
where $\Ad : G(\M) \to \GL\big( \Omega(\g) \big)$, $g \mapsto \Ad_g$ denotes the action of the adjoint group $G(\M)$ on $\Omega(\g)$.
If $g \in G^{\hat\up}(\M)$ then $dg$ is seen to be $\Gamma$-equivariant so that $dg g^{-1} \in \Omega^{\hat\up}(\g)$. It follows that \eqref{gauge tr def} restricts to an action of $G^{\hat\up}(\M)$ on the affine subspace $\Conn^{\hat\up}_\g(\CP)$ of cyclotomic $\g$-connections.

We say that $d + \ms A \in \Conn_\g(\CP)$ has a \emph{pole} (resp. is \emph{regular}) at some point $x \in \CP$ if the $\g$-valued differential $\ms A \in \Omega(\g)$ has a pole (resp. is regular) there. We define the residue of $d + \ms A$ at $x \in \CP$ as
\begin{equation*}
\res_x (d + \ms A) \coloneqq \res_x \ms A.
\end{equation*}
Let $\{ x_i \}_{i=1}^n \subset \CP$ be the set of all poles of the $\g$-connection $d + \ms A$. Since a meromorphic connection on $\CP$ is always flat, $d + \ms A$ gives rise to a group homomorphism
\begin{equation*}
M : \pi_1\big( \CP \setminus \{ x_i \}_{i=1}^n \big) \longrightarrow G,
\end{equation*}
called the \emph{monodromy representation}.
The monodromy of $d + \ms A$ at one of its poles $x_i$ is defined as $M([\gamma_i]) \in G$ where $[\gamma_i]$ is the homotopy class in $\CP \setminus \{ x_i \}_{i=1}^n$ of a small loop $\gamma_i$ encircling the point $x_i$. The $\g$-connection $d + \ms A$ is said to have \emph{trivial monodromy at} $x_i$ if $M([\gamma_i]) = \Id$. Moreover, $d + \ms A$ has \emph{trivial monodromy} or is \emph{monodromy-free} if the homomorphism $M$ is trivial.

\begin{proposition}\label{PropMonodromy}
The $\g$-connection $d + \ms A \in \Conn_\g(\CP)$ has trivial monodromy at $x \in \CP$ if and only if there exists $g\in G(\M)$ such that the gauge transformed $\g$-connection $g (d + \ms A) g^{-1}$ is regular at $x$. Moreover, it is monodromy-free if and only if there exists a solution $g \in G(\M)$ to $dg g^{-1} = - \ms A$.\end{proposition}

Given a $\g$-connection $\nabla = d + \ms A \in \Conn_\g(\CP)$ without monodromy, we shall say that $Y \in G(\M)$ is a solution of the equation $\nabla Y = 0$ if $dY Y^{-1} = - \ms A$. We will need the following later.

\begin{lemma}\label{BorelSol}
Let $\nabla \in \Conn_\g(\CP)$ be monodromy-free and be regular at $x \in \C$. Let $Y \in G(\M)$ be a solution of $\nabla Y = 0$ such that $Y(x) \in B_-$. Then $\nabla \in \Conn_{\b_-}(\CP)$ if and only if $Y \in B_-(\M)$. \qed
\end{lemma}

\subsection{Cyclotomic $\g$-opers} \label{sec: cyclo op}

Recall the principal $\sl_2$-triple of $\g$ defined in \S\ref{Principal} and the associated notion of finite $\g$-oper introduced in \S\ref{finite opers}.

Consider the following subset of $\g$-valued connections
\begin{equation*}
\op_\g(\CP) \coloneqq \big\{ d + p_{-1} dt + \ms v \,\big|\, \ms v \in \Omega(\b) \big\} \subset \Conn_\g(\CP),
\end{equation*}
where we recall that $p_{-1} dt$ is used as a shorthand notation for $p_{-1} \otimes dt$. This set is stable under the gauge action of the unipotent subgroup $N(\M) \subset G(\M)$. The space of \emph{$\g$-opers} is defined as the quotient space
\begin{equation*}
\Op_\g(\CP) \coloneqq \op_\g(\CP) \big/ N(\M).
\end{equation*}
We shall denote by $[\nabla]$ the class of $\nabla \in \op_\g(\CP)$ in $\Op_\g(\CP)$.

In order to define a cyclotomic version of the space $\Op_\g(\CP)$, it is natural to first seek an analog of the subspace $\op_\g(\CP)$ within the space of cyclotomic $\g$-connections $\Conn_\g^{\hat\up}(\CP)$. However, it is clear that for a generic automorphism $\up \in \Aut \g$ of the form \eqref{upsilon order T}, the $\g$-valued differential $p_{-1} dt$ does not live in the $\hat\up$-invariant subspace $\Omega^{\hat\up}(\n_-)$. In fact, the requirement that $p_{-1} dt \in \Omega^{\hat\up}(\n_-)$ fixes the automorphism $\up$ uniquely to be equal to
\begin{equation} \label{Tau}
\vs \coloneqq \Ad_{\omega^{-\co\rho}} \circ \, \nu \in \Aut \g.
\end{equation}
Explicitly, the action of this automorphism on the Chevalley-Serre generators is given by
\begin{equation*}
\vs(E_i) = \omega^{-1} E_{\nu(i)}, \qquad
\vs(\co \alpha_i) = \co \alpha_{\nu(i)}, \qquad
\vs(F_i) = \omega F_{\nu(i)}.
\end{equation*}
In the remainder of this section we shall specialise the constructions of \S\ref{cycl func diff} and \S\ref{TwistedAlgebra} for an arbitrary automorphism $\up$ of the form \eqref{upsilon order T} to the specific automorphism $\vs$ defined by \eqref{Tau}. We will need the following properties of the latter.

\begin{proposition}\label{ExpoSigma}
The automorphism $\vs$ preserves the $\Z$-grading \eqref{Z-grading} and its action on the principal $\sl_2$-triple is given by $\vs(p_{\pm 1}) = \omega^{\mp 1} p_{\pm 1}$ and $\vs(\co \rho) = \co \rho$. Moreover, each $\Vc{i}$, $i \in E$ is $\vs$-stable.
\begin{proof}
This is a direct consequence of Proposition \ref{ExpoNu} and the definition \eqref{Tau} of $\vs$.
\end{proof}
\end{proposition}

It follows from Proposition \ref{ExpoSigma} that $\hat\vs(p_{-1} dt) = p_{-1} dt$ and hence $p_{-1} dt \in \Omega^{\hat\vs}(\n_-)$, as required.

We may now consider the following subset of $\Gamma$-invariant $\g$-valued connections
\begin{equation}\label{CycConn}
\op_\g^\Gamma(\CP) \coloneqq \big\{ d + p_{-1} dt + \ms v \,\big|\, \ms v \in \Omega^{\hat\vs}(\b) \big\} \subset \Conn^{\hat\vs}_\g(\CP).
\end{equation}
It is stable under the gauge action of the unipotent subgroup $N^{\hat\vs}(\M) \subset G^{\hat\vs}(\M)$ on $\Conn^{\hat\vs}_\g(\CP)$. We define the space of \emph{cyclotomic $\g$-opers} as the quotient space
\begin{equation*}
\Op^\Gamma_\g(\CP) \coloneqq \op_\g^\Gamma(\CP) \big/ N^{\hat\vs}(\M).
\end{equation*}
The class of $\nabla$ in $\Op^\Gamma_\g(\CP)$ is denoted by $[\nabla]_\Gamma$. We note that the space $\Op_\g^\Gamma(\CP)$ only depends on the diagram automorphism $\nu : I \to I$ and the choice of $T^{\rm th}$-root of unity $\omega$ through the automorphism $\vs$ in \eqref{Tau}. Any cyclotomic $\g$-connection $\nabla \in \op^\Gamma_\g(\CP)$ is also an element of $\op_\g(\CP)$. As such, we can also consider its orbit in $\op_\g(\CP)$ under the action of $N(\M)$, namely its class $[\nabla] \in \Op_\g(\CP)$. Moreover, it is clear that if $\nabla,\nabla' \in\op^\Gamma_\g(\CP)$ are such that $[\nabla]_\Gamma=[\nabla']_\Gamma$ then $[\nabla]=[\nabla']$. This gives rise to a map
\begin{equation} \label{Op cyclo to Op}
\begin{split}
 \Op^\Gamma_\g(\CP) & \longrightarrow \Op_\g(\CP), \\
 {[\nabla]_\Gamma}    & \longmapsto     [\nabla].
\end{split}
\end{equation}
It follows from Theorem \ref{TwistedDS} below that this map is injective.

\subsection{Canonical representatives}

Recall from Theorem \ref{prop: mini opers} that each finite $\g$-oper admits a unique representative in the Slodowy slice $p_{-1} + \a$. Similarly, it is well known that a transverse slice in the space $\Op_\g(\CP)$ is given by the so called Drinfel'd-Sokolov gauge \cite{Drin}. We will prove a cyclotomic version of this result in Theorem \ref{TwistedDS} below.

The following lemma can be regarded as a (cyclotomic) affine analog of equation \eqref{gi decomp Vc}.

\begin{lemma}\label{LemmeAlgo}
Let $\ms w \in \Omega^{\hat\vs}(\g_i)$.

If $i \not \in E$ then there exists a unique $m \in \g_{i+1}^{\hat\vs}(\M)$ such that $\ms w = [m, p_{-1} dt]$.

If $i\in E$ then there exists unique $m \in \g_{i+1}^{\hat\vs}(\M)$ and $\ms c \in \Omega^{\hat\vs}(\Vc{i})$ such that $\ms w = [m, p_{-1} dt] + \ms c$.
\begin{proof}
We will show this for $i \in E$, the case $i \not\in E$ being similar. Suppose $\ms w \in \Omega^{\hat\vs}(\g_i)$. By definition we can write $\ms w$ as a finite sum of terms of the form $X_j \otimes \varpi_j$ with $X_j \in \g_i$ and $\varpi_j \in \Omega$. Applying the decomposition \eqref{gi decomp Vc} to each $X_j$ we may write
\begin{equation} \label{decomp LemmeAlgo}
\ms w = [m, p_{-1}dt] + \ms c
\end{equation}
for some unique $\ms c \in \Omega(\Vc{i})$ and $[m, p_{-1}dt] \in [\g_{i+1}(\M), p_{-1} dt]$. In turn, $m \in \g_{i+1}(\M)$ is then also unique by the injectivity of $\ad_{p_{-1}} : \g_{i+1} \to \g_i$ for $i \geq 0$. Applying $\hat\vs$ to both sides of \eqref{decomp LemmeAlgo} we obtain $\ms w = \big[ \hat\vs(m), p_{-1}dt \big] + \hat\vs(\ms c)$ using the fact that $\ms w$ and $p_{-1} dt$ are both $\hat\vs$-stable. It now follows by the uniqueness of $\ms c$ and $m$ in \eqref{decomp LemmeAlgo} together with the fact that $\a_i$ and $\g_{i+1}$ are $\vs$-stable by Proposition \ref{ExpoSigma}, that $\ms c \in \Omega^{\hat\vs}(\Vc{i})$ and $m \in \g_{i+1}^{\hat\vs}(\M)$, as required.
\end{proof}
\end{lemma}

\begin{theorem}\label{TwistedDS}
The action of $N^{\hat\vs}(\M)$ on $\op^\Gamma_\g(\CP)$ (resp. of $N(\M)$ on $\op_\g(\CP)$) is free.
Moreover, every cyclotomic $\g$-oper $[\nabla]_\Gamma \in \Op^\Gamma_\g(\CP)$ (resp. $\g$-oper $[\nabla] \in \Op_\g(\CP)$) has a unique representative of the form $d + p_{-1} dt + \ms c$ with $\ms c \in \Omega^{\hat\vs}(\a)$ (resp. $\ms c \in \Omega(\a)$). We shall call it the \emph{canonical representative} of $[\nabla]_\Gamma$ (resp. $[\nabla]$).
\begin{proof}
We consider only the case of a cyclotomic $\g$-oper $[\nabla]_\Gamma \in \Op^\Gamma_\g(\CP)$ since the proof in the case of a $\g$-oper $[\nabla] \in \Op_\g(\CP)$ is standard \cite{Drin} and completely analogous.

Let $\nabla = d + p_{-1} dt + \ms v \in \op^\Gamma_\g(\CP)$. We want to find $m \in \n^{\hat\vs}(\M)$ such that
\begin{equation} \label{eq sol DS}
d + p_{-1} dt + \ms v = e^m ( d + p_{-1} dt + \ms c ) e^{- m}
\end{equation}
for some $\ms c \in \Omega^{\hat\vs}(\a)$. Let us decompose $m$, $\ms v$ and $\ms c$ as follows
\begin{equation}
m = \sum_{i>0} m_i, \qquad \ms v = \sum_{i\geq 0} \ms v_i, \qquad \ms c = \sum_{i \in E} \ms c_i 
\end{equation}
with $m_i \in \g_i^{\hat\vs}(\M)$, $\ms v_i \in \Omega^{\hat\vs}(\g_i)$ and $\ms c_i \in \Omega^{\hat\vs}(\Vc{i})$.
We can then write the three different terms on the right hand side of \eqref{eq sol DS} as
\begin{align*}
d - d(e^m) e^{-m} &= d - \sum_{i>0} d m_i - \frac{1}{2} \sum_{i,j>0} [m_j, d m_i] - \dfrac{1}{6} \sum_{i,j,k} \left[m_k,[m_j, d m_i]\right] + \ldots \\
e^{\ad_m} p_{-1} dt &= p_{-1} dt + \sum_{i>0} [m_i, p_{-1} dt] + \frac{1}{2} \sum_{i,j>0} \left[m_j,[m_i, p_{-1} dt] \right] + \ldots \\
e^{\ad_m} \ms c &= \sum_{i\in E} \ms c_i + \sum_{i\in E} \sum_{j>0} [m_j,\ms c_i] + \frac{1}{2} \sum_{i\in E}\sum_{j,k>0} \left[m_k,[m_j,\ms c_i] \right] + \ldots
\end{align*}
where the ellipses represent only finitely many terms since $m \in \n^{\hat\vs}(\M)$ is nilpotent. Here we note that $\Ad_{e^m} = e^{\ad_m}$.

We determine $m_i$ and $\ms c_i$ inductively. For this we note that if $\ms w \in \Omega^{\hat\vs}(\g_j)$ for some $j \in \Z$ then
\begin{equation} \label{ad m on w}
\big[ m_{j_1}, [m_{j_2}, \ldots [m_{j_n}, \ms w] \ldots ] \big] \in \Omega^{\hat\vs}(\g_{j+j_1+ \ldots +j_n}),
\end{equation}
for any $j_1, \ldots, j_n \in \Z$.
The terms of degree $i \in \Z_{\geq 0}$ in the equation \eqref{eq sol DS} then read
\begin{equation} \label{AlgoDrinApp}
\ms c_i + [m_{i+1}, p_{-1} dt] = \ms v_i + F_i \big( \{ m_{j+1},d m_{j+1},\ms c_j \}_{j<i} \big)
\end{equation}
where $F_i$ is a sum of commutators of the form \eqref{ad m on w} with $\ms w = \ms c_k$, $d m_k$ or $p_{-1} dt$, and hence lives in $\Omega^{\hat\vs}(\g_i)$.
When $i=0$ the relation \eqref{AlgoDrinApp} simply says that $[m_1, p_{-1} dt] = \ms v_0$. And since $\ms v_0 \in \Omega^{\hat\vs}(\g_0)$, Lemma \ref{LemmeAlgo} yields $m_1 \in \g_1^{\hat\vs}(\M)$ uniquely.
Now suppose that we have constructed the $\ms c_j \in \Omega^{\hat\vs}(\Vc{j})$ and $m_{j+1} \in \g_{j+1}^{\hat\vs}(\M)$ for $j<i$. By applying Lemma \ref{LemmeAlgo} we can then construct, in a unique way, $m_{i+1} \in \g_{i+1}^{\hat\vs}(\M)$ and $\ms c_i \in \Omega^{\hat\vs}(\Vc{i})$ for all $i \geq 1$ (with $\ms c_i = 0$ if $i \not\in E$). In particular, the freeness of the action of $N^{\hat\vs}(\M)$ on $\op^\Gamma_\g(\CP)$ follows from the uniqueness of $m_{i+1} \in \g_{i+1}^{\hat\vs}(\M)$.
\end{proof}
\end{theorem}

The proof of Theorem \ref{TwistedDS} provides an algorithm for constructing all the $\ms c_i$ and $m_i$ from the $\ms v_i$. In particular, a formula for $\ms c_1$ can be obtained by considering the first few steps of the algorithm. For this we need the explicit form of \eqref{AlgoDrinApp} for $i=0,1$, which read
\begin{subequations} \label{AlgoDrinFirst}
\begin{align}
\label{AlgoDrinFirst a} [m_1, p_{-1} dt] &= \ms v_0,\\
\label{AlgoDrinFirst b} \ms c_1 + [m_2, p_{-1} dt] &= \ms v_1 - \ha \big[ m_1, [m_1, p_{-1} dt] \big] + d m_1.
\end{align}
\end{subequations}

\begin{proposition}\label{u1}
Write $\ms v_i = v_i dt$ with $v_i \in \g_i(\M)$ for $i \in \Z_{\geq 0}$. Then $\ms c_1 = u_1 p_1 dt$ with
\begin{equation*}
u_1 = \frac{1}{2 (\co\rho | \co\rho)}
\big( \ha (v_0|v_0) + (\co \rho|v_0') + (p_{-1} |v_1) \big).
\end{equation*}
\begin{proof}
Let $\ms w_1 = w_1 dt$ with $w_1 \in \g_1(\M)$ denote the left hand side of equation \eqref{AlgoDrinFirst b}.
Then
\begin{equation*}
(p_{-1} | w_1) = u_1 (p_{-1}|p_1) + \big( p_{-1} \big| [m_2, p_{-1}] \big) = u_1 (p_{-1}|p_1) - \big( [p_{-1},p_{-1}] \big| m_2 \big) = 2 u_1(\co \rho|\co \rho),
\end{equation*}
where in the last step we used $(p_{-1}|p_1) = (p_{-1}|[\co \rho,p_1]) = ([p_1,p_{-1}]|\co \rho) = 2(\co \rho|\co \rho)$. It now follows that
\begin{equation*}
u_1 = \frac{(p_{-1} | w_1 )}{2 (\co \rho|\co \rho)}.
\end{equation*}
Similarly, equation \eqref{AlgoDrinFirst a} implies
$(p_{-1} | m_1) = ([p_{-1},\co \rho] | m_1) = (\co \rho | [m_1, p_{-1}]) = (\co \rho | v_0)$.
So taking the derivative with respect to $t$ yields $(p_{-1} | m'_1 ) = (\co \rho | v'_0)$.
Next, we also have
\begin{equation*}
\big( p_{-1} \big| \big[m_1, [m_1,p_{-1}] \big] \big) = - \big( [m_1,p_{-1}] | [m_1,p_{-1}] \big) = - (v_0|v_0).
\end{equation*}
Equation \eqref{AlgoDrinFirst b} then implies $(p_{-1} | w_1) = \ha (v_0|v_0) + (\co \rho|v'_0) + (p_{-1} |v_1)$, as required.
\end{proof}
\end{proposition}

\subsection{Regular points and regular singularities} \label{sec: opers RS}

We shall say that a (cyclotomic) $\g$-connection $\nabla = d + p_{-1} dt + \ms v \in \op^\Gamma_\g(\CP)$ is \emph{regular} at $x \in \C$ if $\ms v \in \Omega^{\hat\vs}(\b)$ has no pole at $x$.

\begin{proposition} \label{OrbitReg}
Let $\nabla \in \op^\Gamma_\g(\CP)$ and $x \in \C$. The following are equivalent:
\begin{enumerate}
\item[$(i)$] the canonical representative of $[\nabla]_\Gamma$ is regular at $x$
\item[$(ii)$] there exists a representative of $[\nabla]_\Gamma$ which is regular at $x$
\item[$(iii)$] there exists a representative of $[\nabla]$ which is regular at $x$
\end{enumerate}
In this case we will say that the cyclotomic $\g$-oper $[\nabla]_\Gamma$ itself is \emph{regular} at the point $x$.
\begin{proof}
It is obvious that $(i)$ implies $(ii)$ and $(iii)$. To see that $(ii)$ implies $(i)$, we have to look back at the construction of the canonical representative (cf. Theorem \ref{TwistedDS}).

Let $d + p_{-1} dt + \ms v \in \op^\Gamma_\g(\CP)$ be regular at $x$. In the inductive step of the proof of Theorem \ref{TwistedDS}, if the right hand side of \eqref{AlgoDrinApp} is regular at $x$ then so are $\ms c_i$ and $m_{i+1}$ appearing on the left hand side. We therefore conclude that the canonical representative of $[d + p_{-1} dt + \ms v]_\Gamma$ is regular at $x$.
Now if $\nabla$ and $d + p_{-1} dt + \ms v$ lie in the same $N^{\hat\vs}(\M)$-orbit then they share the same canonical representative, which proves that $(ii)$ implies $(i)$.

Finally, we prove that $(iii)$ implies $(i)$. Let $d + p_{-1} dt + \ms c$ and $d + p_{-1} dt + \tilde{\ms c}$ denote the canonical representatives of $[\nabla]_\Gamma$ and $[\nabla]$ respectively, which exist and are unique by Theorem \ref{TwistedDS}. As $N^{\hat\vs}(\M)$ is a subgroup of $N(\M)$, $d + p_{-1} dt + \ms c$ is also a representative of $[\nabla]$ with $\ms c \in \Omega^{\hat\vs}(\a) \subset \Omega(\a)$. Yet, $d + p_{-1} dt + \tilde{\ms c}$ is the unique such representative of $[\nabla]$, hence $\tilde{\ms c} = \ms c$.
Now, if there is a $\g$-connection in the $N(\M)$-orbit of $\nabla$ which is regular at $x \in \C$, then the same argument as in the above proof that $(ii)$ implies $(i)$ applies so that $\tilde{\ms c}$ is regular at $x$, which concludes the proof since $\tilde{\ms c} = \ms c$.
\end{proof} 
\end{proposition}

At any point on $\CP$ where a cyclotomic $\g$-oper $[\nabla]_\Gamma \in \Op^\Gamma_\g(\CP)$ is not regular, the mildest possible singularity it can have is a \emph{regular singularity}. In the remainder of this section we define the notion of regular singularity at a non-zero finite point $x \in \C^\times$ and then at the origin and infinity. We shall say that $[\nabla]_\Gamma$ has \emph{at most a regular singularity} at $x \in \CP$ if it is either regular at $x$ or has a regular singularity there.

Suppose for the moment that $x \in \C^\times$ and consider the subset of $\op_\g(\CP)$ given by
\begin{equation*}
\op_\g(\CP)^{\text{RS}}_x \coloneqq \bigg\{ d + p_{-1} dt + \sum_{i \geq 0} (t-x)^{-i-1} \,\ms v_i \in \op_\g(\CP) \,\bigg|\, \ms v_i \in \Omega(\g_i) \; \text{regular at} \; x \bigg\}.
\end{equation*}
Let $\mathcal O_x \subset \M$ be the algebra of meromorphic functions on $\CP$ which are regular at $x$. One checks using the definition of the $\Z$-grading in \S\ref{Principal} that elements of $\op_\g(\CP)^{\text{RS}}_x$ can equivalently be characterised as $\g$-connections $\nabla \in \op_\g(\CP)$ whose gauge transformation by $(t - x)^{\co\rho}$ take the form
\begin{equation} \label{gauge tr RS}
(t-x)^{\co\rho} \, \nabla \, (t-x)^{-\co\rho} = d + \frac{1}{t-x} ( p_{-1} - \co\rho + v ) dt
\end{equation}
with $v = \sum_{i \geq 0} v_i \in \b(\mathcal O_x)$. It is clear that the space $\op_\g(\CP)^{\text{RS}}_x$ is stabilised by the gauge action of the group $(t-x)^{-\co\rho} N(\mathcal O_x) (t-x)^{\co\rho}$.
We define the space $\Op_\g(\CP)^{\text{RS}}_x$ of \emph{$\g$-opers with regular singularity at $x$} as the corresponding quotient
\begin{equation*}
\Op_\g(\CP)^{\text{RS}}_x \coloneqq \op_\g(\CP)^{\text{RS}}_x \Big/ (t-x)^{-\co\rho} N(\mathcal O_x) (t-x)^{\co\rho}.
\end{equation*}

Since $(t-x)^{-\co\rho} N(\mathcal O_x) (t-x)^{\co\rho} \subset N(\M)$, we obtain a well defined canonical map
\begin{equation} \label{Op RS to Op}
\Op_\g(\CP)^{\textup{RS}}_x \longrightarrow \Op_\g(\CP)
\end{equation}
for each $x \in \C^\times$, which sends the class in $\Op_\g(\CP)^{\textup{RS}}_x$ of a $\g$-connection $\nabla \in \op_\g(\CP)^{\text{RS}}_x \subset \op_\g(\CP)$ to its class $[\nabla]$ in $\Op_\g(\CP)$. By virtue of the next lemma we will usually identify a class in $\Op_\g(\CP)^{\textup{RS}}_x$ with its image under \eqref{Op RS to Op}.

\begin{lemma} \label{lem: Op RS to Op}
The map \eqref{Op RS to Op} is injective. Moreover, its image consists of $\g$-opers whose canonical representative $d + p_{-1} dt + \sum_{i \in E} \ms c_i$ is such that for each $i \in E$, $\ms c_i \in \Omega(\a_i)$ has a pole of order at most $i + 1$ at $x$.

We will say that $x \in \C^\times$ is a \emph{regular singularity} of a cyclotomic $\g$-oper $[\nabla]_\Gamma \in \Op^\Gamma_\g(\CP)$ if it is not regular there, cf. Proposition \ref{OrbitReg}, and if the corresponding $\g$-oper $[\nabla]$ lies in the image of \eqref{Op RS to Op}.
\begin{proof}
In order to prove the lemma it suffices to show that for any $\nabla \in \op_\g(\CP)^{\text{RS}}_x$, the canonical representative of $[\nabla] \in \Op_\g(\CP)$ also lives in $\op_\g(\CP)^{\text{RS}}_x$ and that the gauge transformation that brings $\nabla$ to its canonical representative belongs to $(t-x)^{-\co\rho} N(\mathcal O_x) (t-x)^{\co\rho}$.

Let $\nabla = d + p_{-1} dt + \sum_{i \geq 0} (t-x)^{-i-1} \,\ms v_i \in \op_\g(\CP)^{\text{RS}}_x$. The canonical representative of the $\g$-oper $[\nabla]$ is obtained by the algorithm described in the proof of Theorem \ref{TwistedDS}. Recall that one proceeds by induction on the degree by solving equation \eqref{AlgoDrinApp}, where $\ms v_i$ there is replaced by $(t-x)^{-i-1} \,\ms v_i$ in the present case. As the latter has a pole of order at most $i+1$ at $x$, one can check by induction that, in this construction, $m_i$ as a pole at most of order $i$ and $\ms c_i$ as a pole at most of order $i+1$. By definition this means that $e^{\sum_{i > 0} m_i} \in (t-x)^{-\co\rho} N(\mathcal O_x) (t-x)^{\co\rho}$ and $d + p_{-1} dt + \sum_{i \geq 0} \ms c_i \in \op_\g(\CP)^{\text{RS}}_x$, as required.
\end{proof}
\end{lemma}

Note that the gauge action by $(t - x)^{\co\rho}$ on classes in $\Op_\g(\CP)^{\text{RS}}_x$ induces an obvious bijection
\begin{equation*}
\Op_\g(\CP)^{\text{RS}}_x \overset{\sim}\longrightarrow \bigg\{ d + \frac{1}{t - x} \big( p_{-1} - \co\rho + v \big) dt \,\bigg|\, v \in \b(\mathcal O_x) \bigg\} \bigg/ N(\mathcal O_x).
\end{equation*}
The representative $d + \frac{1}{t - x} \big( p_{-1} - \co\rho + v \big) dt$ of a class in the latter quotient has residue $p_{-1} - \co\rho + v(x)$. Under the gauge action by an element $g \in N(\mathcal O_x)$ this residue transforms under the adjoint action of $g(x) \in N$. This yields a notion of \emph{residue} at $x \in \C^\times$ for elements of $\Op_\g(\CP)^{\text{RS}}_x$
defined by
\begin{equation} \label{residue op}
\begin{split}
\res_x : \Op_\g(\CP)^{\text{RS}}_x &\longrightarrow \Op_\g^{\rm fin},\\
[\nabla] &\longmapsto [p_{-1} - \co\rho + v(x)]_\g,
\end{split}
\end{equation}
where $v \in \b(\mathcal O_x)$ is defined in terms of $\nabla$ through \eqref{gauge tr RS}.

We now turn to the subset $\{ 0, \infty \} \subset \CP$ of fixed points under the action $\mu : \Gamma \to \Aut \CP$ introduced in \S\ref{cycl func diff}. For each $z \in \{ 0, \infty \}$ we consider the subset of cyclotomic $\g$-connections defined by
\begin{equation*}
\op^\Gamma_\g(\CP)^{\text{RS}}_z \coloneqq \bigg\{ d + p_{-1} dt + \sum_{i \geq 0} t^{-i-1} \,\ms v_i \in \op^\Gamma_\g(\CP) \,\bigg|\, \ms v_i \in \Omega(\g_i) \; \text{regular at} \; z \bigg\}.
\end{equation*}
By definition of the $\Z$-grading, an element of $\op^\Gamma_\g(\CP)^{\text{RS}}_z$ can alternatively be seen as a cyclotomic $\g$-connection $\nabla \in \op^\Gamma_\g(\CP)$ whose gauge transformation by $t^{\co\rho}$ reads
\begin{equation} \label{gauge tr RS 0}
t^{\co\rho} \, \nabla \, t^{-\co\rho} = d + \frac{1}{t} ( p_{-1} - \co\rho + v ) dt
\end{equation}
with $v = \sum_{i \geq 0} v_i \in \b^{\hat\nu}(\mathcal O_z)$. The space $\op^\Gamma_\g(\CP)^{\text{RS}}_z$ is stabilised by the gauge action of the subgroup $t^{-\co\rho} N^{\hat\nu}(\mathcal O_z) t^{\co\rho} \subset N^{\hat\vs}(\M)$. The corresponding quotient
\begin{equation*}
\Op^\Gamma_\g(\CP)^{\text{RS}}_z \coloneqq \op^\Gamma_\g(\CP)^{\text{RS}}_z \Big/ t^{-\co\rho} N^{\hat\nu}(\mathcal O_z) t^{\co\rho},
\end{equation*}
defines the space of \emph{cyclotomic $\g$-opers with regular singularity at $z$}.

For each $z \in \{ 0, \infty \}$ we have a canonical map
\begin{equation} \label{Op RS to Op 0}
\Op^\Gamma_\g(\CP)^{\textup{RS}}_z \longrightarrow \Op^\Gamma_\g(\CP)
\end{equation}
which sends the class in $\Op^\Gamma_\g(\CP)^{\textup{RS}}_z$ of a cyclotomic $\g$-connection $\nabla \in \op^\Gamma_\g(\CP)^{\text{RS}}_z \subset \op^\Gamma_\g(\CP)$ to its class $[\nabla]_\Gamma$ in $\Op^\Gamma_\g(\CP)$. This is well defined since $t^{-\co\rho} N^{\hat\nu}(\mathcal O_z) t^{\co\rho} \subset N^{\hat\vs}(\M)$. The following lemma is proved in exactly the same way as Lemma \ref{lem: Op RS to Op}. Using this lemma we will usually identify a class in $\Op^\Gamma_\g(\CP)^{\textup{RS}}_z$ with its image under \eqref{Op RS to Op 0}.

\begin{lemma} \label{lem: Op RS to Op 0}
The map \eqref{Op RS to Op 0} is injective. Moreover, its image consists of cyclotomic $\g$-opers whose canonical representative $d + p_{-1} dt + \sum_{i \in E} \ms c_i$ is such that for each $i \in E$, $\ms c_i \in \Omega^{\hat\vs}(\a_i)$ has a pole of order at most $i + 1$ at the origin (if $z = 0$) or a zero of order at least $i-1$ at infinity (if $z = \infty$).

We will say that $z \in \{ 0, \infty \}$ is a \emph{regular singularity} of a cyclotomic $\g$-oper $[\nabla]_\Gamma \in \Op^\Gamma_\g(\CP)$ if it is not regular there, cf. Proposition \ref{OrbitReg}, and it lies in the image of \eqref{Op RS to Op 0}.
\end{lemma}

The gauge action by $t^{\co\rho}$ on classes in $\Op^\Gamma_\g(\CP)^{\text{RS}}_z$ induces a bijection
\begin{equation*}
\Op^\Gamma_\g(\CP)^{\text{RS}}_z \overset{\sim}\longrightarrow \bigg\{ d + \frac{1}{t} ( p_{-1} - \co\rho + v ) dt \,\bigg|\, v \in \b^{\hat\nu}(\mathcal O_z) \bigg\} \bigg/ N^{\hat\nu}(\mathcal O_z).
\end{equation*}
When $z = 0$ (resp. $z = \infty$), if we gauge transform the $\hat\nu$-invariant $\g$-connection $d + \frac{1}{t} ( p_{-1} - \co\rho + v ) dt$ by an element $g \in N^{\hat\nu}(\mathcal O_0)$ (resp. $g \in N^{\hat\nu}(\mathcal O_\infty)$) then its residue $p_{-1} - \co\rho + v(0) \in p_{-1} + \b^\nu$ (resp. $- p_{-1} + \co\rho - v(\infty) \in - p_{-1} + \b^\nu$) transforms by the adjoint action of $g(0) \in N^\nu$ (resp. $g(\infty) \in N^\nu$).
It is therefore natural to define the notion of \emph{residue} at $z$ on an element in $\Op^\Gamma_\g(\CP)^{\text{RS}}_z$ to be (minus) a finite $\g^\nu$-oper. Specifically, we define
\begin{subequations} \label{residue op 0 inf}
\begin{equation} \label{residue op 0}
\begin{split}
\res_0 : \Op^\Gamma_\g(\CP)^{\text{RS}}_0 &\longrightarrow \Op_{\g^\nu}^{\rm fin},\\
[\nabla]_\Gamma &\longmapsto [p_{-1} - \co\rho + v(0)]_{\g^\nu},
\end{split}
\end{equation}
for the origin, where $v \in \b^{\hat\nu}(\mathcal O_0)$ is defined by \eqref{gauge tr RS 0}. At infinity we define instead
\begin{equation} \label{residue op inf}
\begin{split}
\res_\infty : \Op^\Gamma_\g(\CP)^{\text{RS}}_\infty &\longrightarrow - \Op_{\g^\nu}^{\rm fin},\\
[\nabla]_\Gamma &\longmapsto - [p_{-1} - \co\rho + v(\infty)]_{\g^\nu},
\end{split}
\end{equation}
\end{subequations}
where $- \Op^{\rm fin}_{\g^\nu}$ denotes the set $(-p_{-1} + \b^\nu) / N^\nu$ and $v \in \b^{\hat\nu}(\mathcal O_\infty)$ is defined also by \eqref{gauge tr RS 0}.

Note that if a cyclotomic $\g$-oper $[\nabla]_\Gamma \in \Op^\Gamma_\g(\CP)$ is regular at $x \in \C^\times$ then the corresponding $\g$-oper $[\nabla] \in \Op_\g(\CP)$, defined using the map \eqref{Op cyclo to Op}, belongs to $\Op_\g(\CP)^{\text{RS}}_x$ (recall that we identify this space with its image under the injection \eqref{Op RS to Op}). Likewise, if $[\nabla]_\Gamma \in \Op^\Gamma_\g(\CP)$ is regular at the origin then it belongs to $\Op^\Gamma_\g(\CP)^{\text{RS}}_0$.

\begin{lemma} \label{lem: reg op}
If $[\nabla]_\Gamma \in \Op_\g^\Gamma(\CP)$ is regular at $x \in \C^\times$ (resp. at the origin) then $\res_x [\nabla]_\Gamma = [0]_W$ (resp. $\res_0 [\nabla]_\Gamma = [0]_{W^\nu}$).
\begin{proof}
Let $d + p_{-1} dt + \sum_{i \in E} \ms c_i$ denote the canonical representative of $[\nabla]_\Gamma$, so that $\ms c_i \in \Omega(\a_i)$ are all regular at $x$, cf. Proposition \ref{OrbitReg}$(i)$. Then
\begin{equation*}
(t - x)^{\co\rho} \bigg( d + p_{-1} dt + \sum_{i \in E} \ms c_i \bigg) (t - x)^{-\co\rho} = d + \frac{1}{t - x} \bigg( (p_{-1} - \co\rho) dt + \sum_{i \in E} (t - x)^{i+1} \ms c_i \bigg),
\end{equation*}
from which the result follows.
\end{proof}
\end{lemma}

\subsection{Cyclotomic Miura $\g$-opers}
\label{CycMiuraOpers}

Recall the notion of a finite Miura $\g$-oper from \S\ref{finite opers}.

We will call \emph{Miura $\g$-oper} any $\g$-connection of the form $\nabla = d + p_{-1} dt + \ms u$ where $\ms u \in \Omega(\h)$. Let $\MOp_\g(\CP)$ denote the set of all Miura $\g$-opers. Given a Miura $\g$-oper $\nabla \in \MOp_\g(\CP)$ we will refer to $[\nabla]$ as the underlying $\g$-oper.

Similarly, we define a \emph{cyclotomic Miura $\g$-oper} as a $\g$-connection of the form $\nabla = d + p_{-1} dt + \ms u$ with $\ms u \in \Omega^{\hat\vs}(\h) = \Omega^{\hat\nu}(\h)$, where the last equality follows from the definition \eqref{Tau} of $\vs \in \Aut \g$ and the fact that $\h = \g_0$ in the $\Z$-grading of \S\ref{Principal} defined by $\ad_{\co\rho}$. Denote by $\MOp^\Gamma_\g(\CP)$ the set of all cyclotomic Miura $\g$-opers. If $\nabla \in \MOp^\Gamma_\g(\CP)$ then we call $[\nabla]_\Gamma$ its underlying cyclotomic $\g$-oper. In \S\ref{ReproFlag} we shall be interested in describing the preimage of a given cyclotomic $\g$-oper under
\begin{equation} \label{MOp to Op}
\begin{split}
\MOp_\g^\Gamma(\CP) &\longrightarrow \Op_\g^\Gamma(\CP),\\
\nabla &\longmapsto [\nabla]_\Gamma,
\end{split}
\end{equation}
which is a direct analog of \eqref{MOpfin to Opfin} in the finite case.

There is an obvious bijection between cyclotomic $\h$-connections and cyclotomic Miura $\g$-opers
\begin{equation} \label{Conn to MOp}
\begin{split}
\Conn_\h^{\hat\nu}(\CP) &\overset{\sim}\longrightarrow \MOp_\g^\Gamma(\CP),\\
\nabla &\longmapsto \nabla + p_{-1} dt.
\end{split}
\end{equation}
Given a cyclotomic Miura $\g$-oper $\nabla \in \MOp_\g^\Gamma(\CP)$ we denote the associated cyclotomic $\h$-connection by $\overline{\nabla} \coloneqq \nabla - p_{-1} dt \in \Conn_\h^{\hat\nu}(\CP)$.

In the remainder of this section we turn to the study of singularities of cyclotomic Miura $\g$-opers.
For any $x \in \CP$, we let $\Conn_\h^{\hat\nu}(\CP)^{\text{RS}}_x$ denote the subspace of cyclotomic $\h$-connections with a simple pole at $x$. Consider the composition of the above maps \eqref{MOp to Op} and \eqref{Conn to MOp}, namely
\begin{alignat*}{2}
\Conn_\h^{\hat\nu}(\CP) &\overset{\sim}\longrightarrow \MOp_\g^\Gamma(\CP) & &\longrightarrow \Op_\g^\Gamma(\CP),\\
\nabla &\longmapsto \nabla + p_{-1} dt & & \longmapsto [\nabla + p_{-1} dt]_\Gamma.
\end{alignat*}
We will call this the \emph{cyclotomic Miura transform} by analogy with the non-cyclotomic case \cite{Frenkel}.
Its restriction to $\Conn_\h^{\hat\nu}(\CP)^{\text{RS}}_x$ induces a map $\Conn_\h^{\hat\nu}(\CP)^{\text{RS}}_x \to \Op_\g^\Gamma(\CP)^{\text{RS}}_x$ for each $x \in \CP$.
Recalling the notation $[\co\lambda]_W$ for the $W$-linkage class of a coweight $\co\lambda \in \h$ introduced in \S\ref{finite opers}, we obtain the following commutative diagram
\begin{equation*}
\begin{tikzpicture}
\matrix (m) [matrix of math nodes, row sep=3em, column sep=3em,text height=1.5ex, text depth=0.25ex]    
{
\Conn_\h^{\hat\nu}(\CP)^{\text{RS}}_x & \Op_\g^\Gamma(\CP)^{\text{RS}}_x\\
\h & \Op^{\rm fin}_\g\\
};
\path[->] (m-1-1) edge (m-1-2);
\path[->] (m-2-1) edge node[below]{$[\cdot]_W$} (m-2-2);
\path[->] (m-1-1) edge node[left]{$-\res_x$} (m-2-1);
\path[->] (m-1-2) edge node[right]{$\res_x$} (m-2-2);
\end{tikzpicture}
\end{equation*}
for each $x \in \C^\times$, and for the origin and infinity we have
\begin{equation*}
\begin{tikzpicture}
\matrix (m) [matrix of math nodes, row sep=3em, column sep=3em,text height=1.5ex, text depth=0.25ex]    
{
\Conn_\h^{\hat\nu}(\CP)^{\text{RS}}_0 & \Op_\g^\Gamma(\CP)^{\text{RS}}_0\\
\h^\nu & \Op^{\rm fin}_{\g^\nu}\\
};
\path[->] (m-1-1) edge (m-1-2);
\path[->] (m-2-1) edge node[below]{$[\cdot]_{W^\nu}$} (m-2-2);
\path[->] (m-1-1) edge node[left]{$-\res_0$} (m-2-1);
\path[->] (m-1-2) edge node[right]{$\res_0$} (m-2-2);
\end{tikzpicture}
\qquad
\begin{tikzpicture}
\matrix (m) [matrix of math nodes, row sep=3em, column sep=3em,text height=1.5ex, text depth=0.25ex]    
{
\Conn_\h^{\hat\nu}(\CP)^{\text{RS}}_\infty & \Op_\g^\Gamma(\CP)^{\text{RS}}_\infty\\
\h^\nu & - \Op^{\rm fin}_{\g^\nu}\\
};
\path[->] (m-1-1) edge (m-1-2);
\path[->] (m-2-1) edge node[below]{$-[\cdot]_{W^\nu}$} (m-2-2);
\path[->] (m-1-1) edge node[left]{$\res_\infty$} (m-2-1);
\path[->] (m-1-2) edge node[right]{$\res_\infty$} (m-2-2);
\end{tikzpicture}
\end{equation*}

\begin{proposition} \label{MiuraRegSing}
Let $\nabla \in \MOp_\g^\Gamma(\CP)$.
\begin{itemize}
  \item[$(i)$] The underlying cyclotomic $\g$-oper $[\nabla]_\Gamma$ has at most a regular singularity at $x \in \C^\times$ if and only if
\begin{equation*}
\nabla = d + p_{-1} dt - \frac{w \cdot \co\lambda}{t - x} \, dt + \ms r
\end{equation*}
for some $w \in W$, $\ms r \in \Omega(\h)$ regular at $x$ and $\co\lambda \in \h$ such that $\res_x [\nabla]_\Gamma = [\co\lambda]_W$ and $\co\lambda + \co\rho$ is dominant.

  \item[$(ii)$] The underlying cyclotomic $\g$-oper $[\nabla]_\Gamma$ has at most a regular singularity at $0$ (resp. $\infty$) if and only if
\begin{equation*}
\nabla = d + p_{-1} dt - \frac{w \cdot \co\lambda}{t} \, dt + \ms r
\end{equation*}
for some $w \in W^\nu$, $\ms r \in \Omega^{\hat\nu}(\h)$ regular at $0$ (resp. $\infty$) and $\co\lambda \in \h^\nu$ such that $\res_0 [\nabla]_\Gamma = [\co\lambda]_{W^\nu}$ (resp. $\res_\infty [\nabla]_\Gamma = - [\co\lambda]_{W^\nu}$) and $\co\lambda + \co\rho$ is dominant.
\end{itemize} 
\begin{proof}
Let $x \in \C^\times$.
Note that $\nabla \in \MOp_\g^\Gamma(\CP)$ lies in $\op_\g(\CP)^{\text{RS}}_x$ if and only if it has the form
\begin{equation} \label{reg Miura mu}
\nabla = d + p_{-1} dt - \frac{\co\mu}{t - x} \, dt + \ms r
\end{equation}
for some $\co\mu \in \h$ and $\ms r \in \Omega(\h)$ regular at $x$. It follows from Lemma \ref{lem: Op RS to Op} that $[\nabla]_\Gamma$ has at most a regular singularity at $x$ if and only if $\nabla$ is of the form \eqref{reg Miura mu}.

By definition of the residue at $x \in \C^\times$ in \eqref{residue op} we have $\res_x [\nabla]_\Gamma = [p_{-1}-\co\rho-\co\mu]_\g=[\co\mu]_W$, from which it follows by Proposition \ref{prop: mini miura opers} that $\res_x [\nabla]_\Gamma = [\co\lambda]_W$ if and only if $\co\mu = w \cdot \co\lambda$ for some $w \in W$, where $\co\lambda \in \h$ can be chosen such that $\co\lambda + \co\rho$ is dominant.

Similarly, for $z \in \{ 0, \infty \}$ we have that $\nabla \in \MOp_\g^\Gamma(\CP)$ lies in $\op^\Gamma_\g(\CP)^{\text{RS}}_z$ if and only if it takes the form
\begin{equation} \label{reg Miura mu 0}
\nabla = d + p_{-1} dt - \frac{\co\mu}{t} \, dt + \ms r
\end{equation}
for some $\co\mu \in \h^\nu$ and $\ms r \in \Omega^{\hat\nu}(\h)$ regular at $z$. It follows from Lemma \ref{lem: Op RS to Op 0} that $[\nabla]_\Gamma$ has at most a regular singularity at $z$ if and only if $\nabla$ is of the form \eqref{reg Miura mu 0}.

In this case, from the definition of the residue at $0$ (resp. $\infty$) in \eqref{residue op 0 inf} we have $\res_0 [\nabla]_\Gamma = [\co\mu]_{W^\nu}$ (resp. $\res_\infty [\nabla]_\Gamma = - [\co\mu]_{W^\nu}$). It then follows from Theorem \ref{prop: mini miura opers cyc} that $\res_0 [\nabla]_\Gamma = [\co\lambda]_{W^\nu}$ (resp. $\res_\infty [\nabla]_\Gamma = - [\co\lambda]_{W^\nu}$) if and only if $\co\mu = w \cdot \co\lambda$ for some $w \in W^\nu$, where $\co\lambda \in \h^\nu$ can be chosen such that $\co\lambda + \co\rho$ is dominant.
\end{proof}
\end{proposition}

As a special case of Proposition \ref{MiuraRegSing} suppose $\nabla \in \MOp_\g^\Gamma(\CP)$ is a cyclotomic Miura $\g$-oper whose underlying cyclotomic $\g$-oper $[\nabla]_\Gamma$ is regular at $x \in \C^\times$.
By Lemma \ref{lem: reg op} we have $\res_x [\nabla]_\Gamma = [0]_W$, and so applying Proposition \ref{MiuraRegSing} we deduce that $\nabla$ is of the form
\begin{equation} \label{Miura w dot 0}
\nabla = d + p_{-1} dt - \frac{w \cdot 0}{t - x} \, dt + \ms r,
\end{equation}
for some $w \in W$ and $\ms r \in \Omega(\h)$ regular at $x$.
In particular, it is possible for a cyclotomic Miura $\g$-oper $\nabla$ to have a simple pole at a point $x$ even if its underlying cyclotomic $\g$-oper $[\nabla]_\Gamma$ is regular there. In the next proposition we give a necessary condition for this to happen.

\begin{proposition} \label{RegBethe}
Let $\nabla \in \MOp^\Gamma_\g(\CP)$ be of the form \eqref{Miura w dot 0}. If $[\nabla]_\Gamma$ is regular at $x$ then
\begin{equation}\label{Bethe2}
( w \cdot 0 | \ms r(x) ) = 0.
\end{equation}
\begin{proof}
If $[\nabla]_\Gamma$ is regular at $x$ then, in particular, the component $\ms c_1$ of its canonical representative is regular at $x$. Writing $\nabla = d + (p_{-1} + u) dt$ and $\ms r = r\, dt$, by Proposition \ref{u1} the latter is proportional to
\begin{equation*}
\ha ( u | u ) + ( \co\rho | u' )
= \frac{(w \cdot 0 | w \cdot 0 + 2 \co\rho)}{2 (t-x)^2} - \frac{(w \cdot 0 | r(t))}{t-x} + \ha \left( r(t)| r(t)\right) + \left(\co\rho \big| r'(t)\right).
\end{equation*}
The double pole term vanishes identically using the identity $( w \cdot \co\mu | w \cdot \co\mu + 2 \co\rho ) = ( \co\mu | \co\mu + 2 \co\rho )$ which is valid for any $\co\mu \in \h$, and the last two terms are both regular at $x$. It follows that $\ms c_1$ is regular at $x$ only if \eqref{Bethe2} holds.
\end{proof} 
\end{proposition}

The above condition \eqref{Bethe2} on the cyclotomic Miura $\g$-oper $\nabla$ of the form \eqref{Miura w dot 0} is, in general, not sufficient for $[\nabla]_\Gamma$ to be regular at a point $x \in \C^\times$. However, if the Weyl group element $w$ is of length one, \emph{i.e.} if $w$ is the reflection $s_i$ with respect to a simple root $\alpha_i$ (in which case $w\cdot 0=-\co\alpha_i$), then \eqref{Bethe2} becomes a necessary and sufficient condition for the underlying cyclotomic $\g$-oper $[\nabla]_\Gamma$ to be regular at $x$.

\begin{proposition}\label{RegBetheSimple}
Let $x \in \C^\times$. Suppose $\nabla \in \MOp^\Gamma_\g(\CP)$ has the form
\begin{equation*}
\nabla = d + p_{-1} dt + \frac{\co\alpha_i}{t - x} \, dt + \ms r
\end{equation*}
for some simple root $\alpha_i$ and $\ms r \in \Omega(\h)$ regular at $x$. Then $[\nabla]_\Gamma$ is regular at $x$ if and only if
\begin{equation}\label{BetheSimple}
\langle \alpha_i, \ms r(x) \rangle = 0.
\end{equation}
\begin{proof}
We have $\co\alpha_i = -s_i \cdot 0$ so that $\nabla$ is of the form \eqref{Miura w dot 0}. Therefore, if $[\nabla]_\Gamma$ is regular at $x$ then applying Proposition \ref{RegBethe} we deduce $(\co\alpha_i| \ms r(x) )=0$, which is equivalent to \eqref{BetheSimple}.

Conversely, suppose \eqref{BetheSimple} holds and let us show $[\nabla]_\Gamma$ is regular at $x$. According to Proposition \ref{OrbitReg}, it is enough to find a (possibly non-cyclotomic) representative of the $\g$-oper $[\nabla]$ which is regular at $x$. Applying a gauge transformation by $g = \exp \big( - \frac{1}{t - x} E_{\alpha_i} \big)$ to the given Miura $\g$-oper $\nabla$, cf. the proof of Lemma \ref{TransfSimple} below, we obtain
\begin{equation*}
g \nabla g^{-1} = d + p_{-1}dt + \ms r(t) + \frac{\langle \alpha_i, r(t) \rangle}{t-x} E_{\alpha_i} dt.
\end{equation*}
This is regular at $x$ by virtue of equation \eqref{BetheSimple}, hence the proposition.
\end{proof}
\end{proposition}

Let $N \in \Z_{\geq 1}$ and fix a finite subset of $N+2$ points $\bm z \coloneqq \{ 0, z_1, \ldots, z_N, \infty \} \subset \CP$. We assume these points have disjoint $\Gamma$-orbits, \emph{i.e.} $\omega^r z_i \neq z_j$ for all $i \neq j$ and $r = 0, \ldots, T-1$.
Let $\Op^\Gamma_\g(\CP)^{\text{RS}}_{\bm z}$ denote the set of all cyclotomic $\g$-opers with at most regular singularities at the points in $\Gamma \bm z$ and which are regular elsewhere.

Let $\bm{\co\lambda} \coloneqq \{ \co\lambda_0, \co\lambda_1, \ldots, \co\lambda_N, \co\lambda_\infty \} \subset \h$ be a collection of $N+2$ integral dominant coweights, which we think of as being attached to the corresponding points in $\bm z$. We assume that the coweights at the origin and infinity are both $\nu$-invariant, namely $\co\lambda_0, \co\lambda_\infty \in \h^\nu$.
We denote by $\Op^\Gamma_\g(\CP)^{\text{RS}}_{\bm z; \bm{\co\lambda}}$ the subset of cyclotomic $\g$-opers $[\nabla]_\Gamma \in \Op^\Gamma_\g(\CP)^{\text{RS}}_{\bm z}$ whose residues at the points of $\bm z$ are given by
\begin{gather*}
\res_0 [\nabla]_\Gamma = [\co\lambda_0]_{W^\nu}, \qquad
\res_{z_i} [\nabla]_\Gamma = [\co\lambda_i]_W, \qquad
\res_\infty [\nabla]_\Gamma = - [\co\lambda_\infty]_{W^\nu}
\end{gather*}
for $i = 1, \ldots, N$.
Let $\Op^\Gamma_\g(\CP)_{\bm z; \bm{\co\lambda}}$ denote the further subset consisting of those cyclotomic $\g$-opers in $\Op^\Gamma_\g(\CP)^{\text{RS}}_{\bm z; \bm{\co\lambda}}$ which are also monodromy-free.

\begin{theorem}\label{GeneralForm}
Let $\bm z \subset \CP$ and $\bm{\co\lambda} \subset \h$ be as above and let $\nabla \in \MOp^\Gamma_\g(\CP)$. The underlying cyclotomic $\g$-oper $[\nabla]_\Gamma$ lives in $\Op^\Gamma_\g(\CP)^{\textup{RS}}_{\bm z; \bm{\co\lambda}}$ if and only if the following conditions hold:
\begin{itemize}
\item[$(i)$] $\nabla$ is of the form
\begin{equation} \label{general cycl Miura}
\nabla = d + p_{-1} dt - \frac{w_0 \cdot \co\lambda_0}{t}dt - \sum_{r=0}^{T-1} \left( \sum_{i=1}^N \frac{\nu^r (w_i \cdot \co \lambda_i)}{t - \omega^r z_i} + \sum_{j=1}^m \frac{\nu^r y_j \cdot 0}{t - \omega^r x_j} \right)dt
\end{equation}
for some $m \in \Z_{\geq 0}$,  $w_0 \in W^\nu$, $w_i \in W$ for each $i = 1, \ldots, N$ and $x_j\in\CP \setminus \Gamma \bm z$, $y_j \in W$ for each $j = 1, \ldots, m$,
\item[$(ii)$] there exists $w_\infty \in W^\nu$ such that
\begin{equation} \label{residue condition}
w_0 \cdot \co\lambda_0 + \sum_{r=0}^{T-1} \left( \sum_{i=1}^N \nu^r (w_i \cdot \co \lambda_i) + \sum_{j=1}^m \nu^r y_j \cdot 0 \right) = w_\infty \cdot \co\lambda_\infty,
\end{equation}
\item[$(iii)$] $[\nabla]_\Gamma$ is regular at $x_j$ for each $j=1,\ldots,m$.
\end{itemize}
\begin{proof}
Suppose that the cyclotomic $\g$-oper $[\nabla]_\Gamma$ belongs to $\Op^\Gamma_\g(\CP)^{\textup{RS}}_{\bm z; \bm{\co\lambda}}$. Then $[\nabla]_\Gamma$ has regular singularities at the points of $\Gamma\bm z$ and is regular everywhere else. However, as we have seen above, the cyclotomic Miura $\g$-oper $\nabla$ itself can have other simple poles at a set of points $\bm{\widetilde{x}} \subset \CP\setminus\Gamma\bm z$. By virtue of cyclotomy, if $\nabla$ has a pole at a point in $\CP$, it also has a pole at each element in the orbit of this point under the action of $\Gamma$. Thus, this set $\widetilde{\bm x}$ can be seen as the image $\Gamma \bm x$ of some minimal set $\bm x \subset \CP\setminus \Gamma\bm z$. Let $m\in\Z_{\geq 0}$ be the size of $\bm x$ and write $\bm x = \{ x_1, \ldots, x_m \}$.

According to Proposition \ref{MiuraRegSing}, there exist $w_0\in W^\nu$ and $w_i \in W$ for each $i=1,\ldots,N$ such that $\res_0 \overline{\nabla} = - w_0\cdot\co\lambda_0$ and $\res_{z_i} \overline{\nabla} = -w_i\cdot\co\lambda_i$. Furthermore, from the discussion around \eqref{Miura w dot 0} there exist $y_j \in W$ for each $j=1,\ldots,m$ such that $\res_{x_j} \overline{\nabla} = -y_j \cdot 0$. Moreover, since $\nabla$ is cyclotomic, one has $\res_{\omega^r z_i} \overline{\nabla} = -\nu^r\left( w_i\cdot\co\lambda_i \right)$ and $\res_{\omega^r x_j} \overline{\nabla} = -\nu^r\left( y_j\cdot 0 \right)=-\nu^r y_j \cdot 0$ for any $r=0,\ldots,T-1$, $i=1,\ldots,N$ and $j=1,\ldots,m$. Thus $\nabla$ is of the form \eqref{general cycl Miura}.

The residue at infinity of the corresponding $\h$-connection is then
\begin{equation*}
\res_\infty \overline{\nabla} = w_0 \cdot \co\lambda_0 + \sum_{r=0}^{T-1} \left( \sum_{i=1}^N \nu^r (w_i \cdot \co \lambda_i) + \sum_{j=1}^m \nu^r y_j \cdot 0 \right).
\end{equation*}
This is related to the residue of the cyclotomic $\g$-oper $[\nabla]_\Gamma$ by $\res_\infty [\nabla]_\Gamma = - \left[ \res_\infty \overline{\nabla} \right]_{W^\nu}$. Yet, as $[\nabla]_\Gamma$ belongs to $\Op^\Gamma_\g(\CP)^{\textup{RS}}_{\bm z; \bm{\co\lambda}}$, we have $\res_\infty [\nabla]_\Gamma = - [\co\lambda_\infty]_{W^\nu}$, hence the existence of $w_\infty \in W^\nu$ such that equation \eqref{residue condition} holds.

Conversely, suppose that $\nabla$ is of the form \eqref{general cycl Miura} and that we have the condition \eqref{residue condition}. It is then clear that, for $i = 1, \ldots, N$,
\begin{gather*}
\res_0 [\nabla]_\Gamma = [\co\lambda_0]_{W^\nu}, \qquad
\res_{z_i} [\nabla]_\Gamma = [\co\lambda_i]_W, \qquad
\res_\infty [\nabla]_\Gamma = -[\co\lambda_\infty]_{W^\nu}.
\end{gather*}
Moreover, suppose that $[\nabla]_\Gamma$ is regular at the points $x_j$, for $j=1,\ldots,m$. Then, by virtue of the cyclotomy, $[\nabla]_\Gamma$ is regular at all the points $\omega^r x_j$, for $r=0,\ldots,T-1$ and $j=1,\ldots,m$. Hence $[\nabla]_\Gamma$ belongs to $\Op^\Gamma_\g(\CP)^{\textup{RS}}_{\bm z; \bm{\co\lambda}}$.
\end{proof}
\end{theorem}

\begin{rem} \label{apres proposition 0}
Let $\nabla \in \MOp^\Gamma_\g(\CP)$ be as in Theorem \ref{GeneralForm}. According to Proposition \ref{RegBethe}, a necessary condition for $[\nabla]_\Gamma$ to be regular at $x_k$, $k \in \{1,\ldots,m\}$ is the \textit{generalised cyclotomic Bethe ansatz equation}
\begin{equation*}
\frac{\bigl( y_k \cdot 0 | w_0 \cdot \co\lambda_0 \bigr)}{x_k} + \sum_{r=0}^{T-1} \sum_{i=1}^N \frac{\bigl(y_k\cdot 0 | \nu^r (w_i \cdot \co \lambda_i) \bigr)}{x_k - \omega^r z_i} + \underset{(r,j) \neq (0, k)}{\sum_{r=0}^{T-1} \sum_{j=1}^m} \frac{\bigl( y_k\cdot 0 | \nu^r y_j \cdot 0 \bigr)}{x_k - \omega^r x_j} = 0.
\end{equation*}
Moreover, Proposition \ref{RegBetheSimple} states that this condition is sufficient if $y_k \in W$ is a simple reflection.
\end{rem}

\begin{rem} \label{apres proposition}
Let $\nabla \in \MOp^\Gamma_\g(\CP)$ be as in Theorem \ref{GeneralForm}. If $w_i = \Id$ for some $i=1,\ldots,N$ (resp. $w_0=\Id$), then $[\nabla]_\Gamma$ has trivial monodromy at $z_i$ (resp. at the origin). Indeed, in this case a gauge transformation of $\nabla$ by $(t-z_i)^{\co\lambda_i}$ (resp. $t^{\co\lambda_0}$) yields a $\g$-connection regular at $z_i$ (resp. at the origin) and hence with trivial monodromy there by Proposition \ref{PropMonodromy}.
\end{rem}

To end this section we present two explicit examples of cyclotomic Miura $\g$-opers, for the simple Lie algebras $\g=\sl_3$ and $\g=\sl_4$, that we shall use to illustrate the various results from \S\ref{Repro} and \S\ref{ReproFlag}.

\begin{exmp}\label{Sl3}
Consider the simple Lie algebra $\g=\sl_3$ in the fundamental representation. The two fundamental coweights of $\sl_3$ are represented by the following diagonal matrices
\begin{equation*}
\co\omega_1=\frac{1}{3}\begin{pmatrix}
2 & 0 & 0 \\ 0 & -1 & 0 \\ 0 & 0 & -1
\end{pmatrix}, \;\;\;\;
\co\omega_2=\frac{1}{3}\begin{pmatrix}
1 & 0 & 0 \\ 0 & 1 & 0 \\ 0 & 0 & -2
\end{pmatrix}.
\end{equation*}
Moreover, the principal nilpotent element $p_{-1}$ is given by
\begin{equation*}
p_{-1} = \begin{pmatrix}
0 & 0 & 0 \\ 1 & 0 & 0 \\ 0 & 1 & 0
\end{pmatrix}.
\end{equation*}

The unique non-trivial diagram automorphism $\nu : I \to I$ of $\sl_3$ exchanges the labels 1 and 2. Let us consider the simplest possible cyclotomic Miura $\sl_3$-oper $\nabla$ with a pole only at the origin and whose residue there is given by a $\nu$-invariant integral dominant coweight. The general form of such a coweight is $\co\lambda_0=\eta(\co\omega_1+\co\omega_2)$ with $\eta\in\Z_{\geq 0}$. Therefore the cyclotomic Miura $\sl_3$-oper in question takes the form
\begin{equation}\label{ExpleSl3}
\nabla = d + \left(p_{-1}-\frac{\co\lambda_0}{t}\right)dt = d +  \begin{pmatrix}
-\frac{\eta}{t} & 0 & 0 \\ 1 & 0 & 0 \\ 0 & 1 & \frac{\eta}{t}
\end{pmatrix}dt.
\end{equation}
Note that the residue at infinity of the associated $\h$-connection $\overline{\nabla}$ is simply $\co\lambda_0=\eta(\co\omega_1+\co\omega_2)$. The canonical representative of the associated cyclotomic $\sl_3$-oper $[\nabla]_\Gamma$ reads
\begin{equation*}
g\nabla g^{-1} = d + \begin{pmatrix}
0 & \frac{\eta(\eta+2)}{2t} & 0 \\ 1 & 0 & \frac{\eta(\eta+2)}{2t} \\ 0 & 1 & 0
\end{pmatrix}dt, \;\;\;\;\; \text{with } \;\;\;\; g(t) = \begin{pmatrix}
1 & \frac{\eta}{t} & \frac{\eta^2}{2t^2} \\ 0 & 1 & \frac{\eta}{t} \\ 0 & 0 & 1
\end{pmatrix}.
\end{equation*}

Since we are working in the fundamental representation, the above group element $g$ belongs to the three dimensional representation of the special linear group $\mathsf{SL}_3$. However, in this article we consider only the adjoint group $G = (\Aut\g)^\circ$ of the Lie algebra $\g$, cf. \S\ref{AdGroup}. In the present example we have $\g=\sl_3$, whose adjoint group is the projective linear group $\mathsf{PGL}_3$ (the quotient of the linear group $\mathsf{GL}_3$ by its center $\lbrace \lambda \, \Id, \; \lambda \in\C^\times \rbrace$). The adjoint group $\mathsf{PGL}_3$ can also be seen as the quotient of $\mathsf{SL}_3$ by its center $Z=\lbrace \lambda \, \Id, \; \lambda \in\C^\times \, | \, \lambda^3=1 \rbrace$. Here and in the examples \ref{ReproSl3} and \ref{ReproGenSl3} below, all group elements written as $3 \times 3$ matrices are to be understood as the class of these matrices in the three dimensional representation of $\mathsf{SL}_3$ modulo multiplication by elements of $Z$.
\end{exmp}

\begin{exmp}\label{Sl4}
Consider the simple Lie algebra $\g=\sl_4$. Denote by $\alpha_i$ and $\co\omega_i$ for $i=1,2,3$, the simple roots and the associated coweights.

The unique non-trivial diagram automorphism of $\sl_4$ is $\nu: \, 1 \mapsto 3, \, 2 \mapsto 2, \, 3 \mapsto 1$, of order 2. We choose the primitive $T^{\rm th}$-root of unity $\omega = e^{\frac{2 \pi i}{T}}$, with $T=2S$ and $S \in \Z_{\geq 1}$. The general $\nu$-invariant integral dominant coweight is of the form $\co\lambda_0 = \eta(\co\omega_1+\co\omega_3) + \kappa \co\omega_2$ with $\eta,\kappa \in \Z_{\geq 0}$. We shall consider the cyclotomic Miura $\sl_4$-oper specified by the coweight $\co\lambda_0$ at the origin and the coweight $\co\omega_1$ at a point $z \in \C^\times$, namely
\begin{equation*}
\nabla = d + \left( p_{-1} - \frac{\co\lambda_0}{t} - \sum_{r=0}^{T-1} \frac{\nu^r \co\omega_1}{t-\omega^r z} \right)dt.
\end{equation*}
This can be re-expressed as
\begin{equation}\label{ExpleSl4}
\nabla = d + p_{-1}dt - \left( \frac{\eta}{t} + \frac{St^{S-1}}{t^S-z^S} \right)\co\omega_1 dt - \frac{\kappa}{t}\co\omega_2 dt - \left( \frac{\eta}{t} + \frac{St^{S-1}}{t^S+z^S} \right)\co\omega_3 dt.
\end{equation}
The residues at the origin and infinity of the associated $\h$-connection are
\begin{equation*}
\res_0 \overline{\nabla} = -\eta(\co\omega_1+\co\omega_3) - \kappa \co\omega_2, \qquad
\res_\infty \overline{\nabla} = \left( \eta + S \right)(\co\omega_1+\co\omega_3) + \kappa\co\omega_2. \qedhere
\end{equation*}
\end{exmp}

\section{Cyclotomic Miura $\g$-opers and reproductions} 
\label{Repro}

Given a Miura $\g$-oper $\nabla \in \MOp_\g(\CP)$ and any simple root $\alpha_k$, $k \in I$, the following lemma provides a way of constructing a new Miura $\g$-oper $g \nabla g^{-1}$, where the gauge transformation parameter is of the form $g=e^{fE_k}$ with $f\in\M$. This is known as the \emph{reproduction} or \emph{generation procedure} in the direction $\alpha_k$, see \emph{e.g.} \cite{MV}.

\begin{lemma}\label{TransfSimple}
Let $\nabla = d + (p_{-1} + u) dt \in \MOp_\g(\CP)$ and $g = e^{f E_k} \in N(\M)$ for some $f \in \M$. We have $g \nabla g^{-1} \in \MOp_\g(\CP)$ if and only if the function $f$ satisfies the Riccati equation
\begin{equation*}
f'(t) + f(t)^2 + f(t) \langle \alpha_k, u(t) \rangle = 0.
\end{equation*}
In this case $g \nabla g^{-1} = \nabla + f \co\alpha_k \, dt$.
Moreover, for any $x \in \C$ we have
\begin{equation*}
\res_x \overline{g \nabla g^{-1}} = \left\{
\begin{array}{ll}
-s_k \cdot (-\res_x \overline{\nabla}) & \; \text{if } f dt \in \Omega \text{ has a pole at } x,\\ 
\res_x \overline{\nabla} & \; \text{otherwise}.
\end{array}
\right.
\end{equation*}
Similarly, at infinity, we have
\begin{equation*}
\res_\infty \overline{g \nabla g^{-1}} = \left\{
\begin{array}{ll}
s_k \cdot (\res_\infty \overline{\nabla}) & \; \text{if } f dt \in \Omega \text{ has a pole at } \infty,\\ 
\res_\infty \overline{\nabla} & \; \text{otherwise}.
\end{array}
\right.
\end{equation*}
Moreover, if $\langle \alpha_k, \res_\infty \overline{\nabla} + \co\rho \rangle \geq 0$ and $f$ is not identically zero, then $f dt$ has a pole at $\infty$. 
\begin{proof}
By a direct calculation we find that
\begin{gather*}
- dg g^{-1} = - f'(t) E_k dt, \qquad
g p_{-1} g^{-1} = p_{-1} + f(t) \co\alpha_k - f(t)^2 E_k,\\
g u(t) g^{-1} = u(t) - f(t) \langle \alpha_k, u(t) \rangle E_k,
\end{gather*}
from which we obtain
$g \nabla g^{-1} = d + \left[ p_{-1} + u(t) + f(t) \co\alpha_k - \big( f'(t) + f(t)^2 + f(t) \langle \alpha_k, u(t) \rangle \big) E_k \right] dt$.
The latter is a Miura $\g$-oper if and only if the coefficient of $E_k$ vanishes, which is equivalent to the Riccati equation. We then have $g \nabla g^{-1} = \nabla + f \co\alpha_k \, dt$ and hence also $\overline{g \nabla g^{-1}} = \overline{\nabla} + f \co\alpha_k \, dt$.

To prove the statement about the residue, suppose first that $x \in \C$ and define $\co\mu \coloneqq - \res_x \overline{\nabla}$ so that $u(t) = - \frac{\co\mu}{t - x} + r(t)$ for some $r \in \h(\M)$ regular at $x$. By considering the Laurent series expansion at $x$ of the solution $f \in \M$ we find that it is consistent with the Riccati equation provided $f$ has at most a simple pole there, namely
\begin{equation*}
f(t) = \frac{a}{t - x} + \mathcal O\big( (t-x)^0 \big)
\end{equation*}
with either $a = 0$ or $a = \langle \alpha_k, \co\mu + \co\rho \rangle$. In the first case, $f dt$ is regular at $x$ so that $\res_x f dt = 0$ and hence $\res_x \overline{g \nabla g^{-1}} = \res_x \overline{\nabla}$. In the second case, $f dt$ has a simple pole at $x$ with $\res_x f dt = \langle \alpha_k, \co\mu + \co\rho \rangle$ so that $\res_x \overline{g \nabla g^{-1}} = - \co\mu + \langle \alpha_k, \co\mu + \co\rho \rangle \co\alpha_k = - s_k \cdot \co\mu = -s_k \cdot (-\res_x \overline{\nabla})$.

Consider now the point at infinity and define $\co\mu \coloneqq \res_\infty \overline{\nabla}$ so that $u(t) = - \frac{\co\mu}{t} + \mathcal O \big( \frac{1}{t^2} \big)$.
It then follows that the asymptotic behaviour of the solution $f$ to the Riccati equation for large $t$ reads
\begin{equation*}
f(t) = \frac{a}{t} + \mathcal O\left( \frac{1}{t^2} \right),
\end{equation*}
with either $a = 0$ or $a = \langle \alpha_k, \co\mu + \co\rho \rangle$. As above, in the first case, $f dt$ is regular at infinity and we find $\res_\infty \overline{g \nabla g^{-1}} = \res_\infty \overline{\nabla}$. In the second case, $f dt$ has a simple pole at infinity and we deduce that $\res_\infty \overline{g \nabla g^{-1}} = s_k \cdot (\res_\infty \overline{\nabla})$, as required.

Finally, let us suppose that $\langle \alpha_k, \res_\infty{\overline{\nabla}} + \co\rho \rangle$ is non-negative, \textit{i.e.} that $\langle \alpha_k, \co\mu \rangle \geq -1$. In order to prove the last statement of the lemma, we shall show that if $f$ does not have a pole at infinity (\textit{i.e.} if $a=0$), then $f$ is identically zero. Indeed, if $f$ were different from zero, we could write
\begin{equation*}
f(t) = \frac{b}{t^p} + \mathcal O\left( \frac{1}{t^{p+1}} \right),
\end{equation*}
for some non-zero $b\in\C^\times$ and $p \in\Z_{\geq 2}$. The cancellation of the terms of order $t^{-p-1}$ in the Ricatti equation yields the condition $\bigl(p+\langle\alpha_k,\co\mu \rangle\bigr)b = 0$. However, $p+\langle\alpha_k,\co\mu\rangle \geq 2-1=1$ so that $b=0$, hence a contradiction.
\end{proof}
\end{lemma}

Suppose now that the Miura $\g$-oper $\nabla$ we start with is cyclotomic. More precisely, let us consider $\nabla \in \MOp^\Gamma_\g(\CP)$ such that the underlying cyclotomic $\g$-oper $[\nabla]_\Gamma$ lives in $\Op^\Gamma_\g(\CP)^{\textup{RS}}_{\bm z; \bm{\co\lambda}}$, as defined in \S\ref{CycMiuraOpers}. Then $\nabla$ is described by Theorem \ref{GeneralForm} and in particular is of the form \eqref{general cycl Miura}. Moreover, using the notations of Theorem \ref{GeneralForm}, we suppose $w_0=\Id$ so that $-\res_0 \overline{\nabla}=\co\lambda_0$ is a $\nu$-invariant integral dominant coweight. We can thus write $\nabla = d + (p_{-1} + u) dt$ with
\begin{equation} \label{cycl Miura lambda0}
u(t) \coloneqq - \frac{\co\lambda_0}{t} - \sum_{r=0}^{T-1} \left( \sum_{i=1}^N \frac{\nu^r (w_i \cdot \co \lambda_i)}{t - \omega^r z_i} + \sum_{j=1}^m \frac{\nu^r y_j \cdot 0}{t - \omega^r x_j} \right).
\end{equation}

The difficulty with applying the generation procedure in the cyclotomic setting is that the new Miura $\g$-oper $g \nabla g^{-1}$ will in general no longer be cyclotomic, since the parameter $g \in N(\M)$ of the gauge transformation we applied need not be $\hat\vs$-invariant. This problem was first studied in \cite{CharlesVar}, when $w_i = \Id$ for each $i = 1, \ldots, N$ and $y_j$ are simple reflections for each $j = 1, \ldots, m$, using the language of populations of solutions to the Bethe ansatz equations (the collection of algebraic equations ensuring the regularity of \eqref{cycl Miura lambda0} at the set of points $\{ x_j \}_{j=1}^m$, cf. Proposition \ref{RegBetheSimple}). As observed there, it is possible to move from one solution of the \emph{cyclotomic} Bethe equations \cite{VY1} to another by applying instead a sequence of reproductions in the directions of the various roots in the $\nu$-orbit of $\alpha_k$. Moreover, the existence of such a \emph{cyclotomic generation procedure} was shown to require imposing certain conditions on the coweight $\co\lambda_0$, cf. Theorems \ref{CommutingOrbit} and \ref{NonCommutingOrbit} below. Note that for the case discussed in \S\ref{sec: A2 orbit} below $\co\lambda_0$ was only assumed \emph{half}-integral dominant in \cite{CharlesVar}. Although we assume throughout this section that $\co\lambda_0$ is integral dominant, we shall show later in \S\ref{sec: non-integral} that the integrality condition can in fact be relaxed.

In this section we will describe cyclotomic reproduction within the framework of cyclotomic Miura $\g$-opers by studying the properties of Riccati-type equations. Specifically, we will prove cyclotomic generalisations of Lemma \ref{TransfSimple}, corresponding to reproduction in the direction of a simple root $\alpha^\nu_{\mathcal I}$ of the folded Lie algebra $\g^\nu$, in the case when the orbit $\mathcal I \in I/\nu$ is of type $A_1^{\times |\mathcal I|}$ or $A_2^{\times |\mathcal I| /2}$.

Recall the setting of \S\ref{folding}. Let $\mathcal I \in I/\nu$ be an orbit in $I$ under the diagram automorphism $\nu : I \to I$. Consider the corresponding orbit of simple roots
\begin{equation*}
\{ \alpha_i \}_{i \in \mathcal I} \subset \Phi^+.
\end{equation*}
Unless the orbit is of size one, \emph{i.e.} $|\mathcal I| = 1$, applying a single gauge transformation by an element of the form $e^{f_i E_i}$ with $i \in \mathcal I$ and $f_i \in \M$ will certainly break the $\Gamma$-invariance of the cyclotomic Miura $\g$-oper $\nabla$. In order to restore it, we will have to apply a series of successive gauge transformations by $e^{f_i E_i}$, with some $f_i \in \M$, corresponding to the other points $i$ in the orbit $\mathcal I$.
The specific combination of gauge transformations required will depend on whether the orbit $\mathcal I$ is of type $A_1^{\times |\mathcal I|}$ or $A_2^{\times |\mathcal I|/2}$, cf. \S\ref{folding}. We will treat these two cases separately.

\subsection{Orbit of type $A_1^{\times |\mathcal I|}$} \label{sec: A1 orbit}

Let $\mathcal I$ be an orbit of type $A_1^{\times |\mathcal I|}$, \emph{i.e.} such that $\ell_{\mathcal I} = 1$ in the notation of \S\ref{folding}. As mentioned above the case $|\mathcal I| = 1$ is easy to treat, so we shall assume $|\mathcal I| \geq 2$. 

Let $i, j \in \mathcal I$ be distinct. If we perform a gauge transformation of $\nabla = d + (p_{-1} + u)dt \in \MOp^\Gamma_\g(\CP)$ by $e^{f_i E_i}$ with $f_i \in \M$ a solution to the Riccati equation
\begin{equation*}
f'(t) + f(t)^2 + f(t) \langle \alpha_i, u(t) \rangle = 0,
\end{equation*}
then by Lemma \ref{TransfSimple} we obtain the new Miura $\g$-oper
$e^{f_i E_i} \nabla e^{- f_i E_i} = \nabla + f_i \co\alpha_i \, dt$.
If we perform a second gauge transformation of the latter by $e^{f_j E_j}$, then by Lemma \ref{TransfSimple} the resulting $\g$-connection will again be a Miura $\g$-oper provided $f_j \in \M$ is chosen to satisfy the Riccati equation
\begin{equation*}
f'(t) + f(t)^2 + f(t) \langle \alpha_j, u(t) + f_i \co\alpha_i \rangle = 0.
\end{equation*}
However, since we are assuming $\mathcal I$ to be of type $A_1^{\times |\mathcal I|}$, we have $\langle \alpha_j, \co\alpha_i \rangle = 0$ and hence the Riccati equation satisfied by $f_j$ only depends on the original Miura $\g$-oper $\nabla$ and not on $f_i$.

We are therefore seeking a gauge transformation by $g \coloneqq \prod_{i \in \mathcal I} e^{f_i E_i}$ where the functions $f_i \in \M$ satisfy the Riccati equations
\begin{equation}\label{RiccatiCommuting}
(R^i): \qquad
f_i'(t) + f_i(t)^2 + f_i(t) \langle \alpha_i, u(t) \rangle = 0,
\end{equation}
for $i \in \mathcal I$. Note that the order of the factors in the expression for $g$ doesn't matter by the assumption on the nature of the orbit $\mathcal I$.
By Lemma \ref{TransfSimple}, the resulting Miura $\g$-oper is given by
\begin{equation} \label{new Miura comm}
g \nabla g^{-1} = \nabla + \sum_{i \in \mathcal I} f_i(t) \, \co\alpha_i \, dt.
\end{equation}
The sum over the orbit $\mathcal I$ is $\hat\nu$-invariant, and hence $g \nabla g^{-1}$ is cyclotomic, if and only if
\begin{equation} \label{cycl cond comm}
\omega^{-1}f_i(\omega^{-1}t) = f_{\nu(i)}(t), \qquad \text{for all} \;\; i \in \mathcal I.
\end{equation}

Let us fix a reference point $k \in \mathcal I$ on the orbit so that
\begin{equation*}
\mathcal I = \big\{ k, \nu(k), \ldots,\nu^{|\mathcal I|-1}(k) \big\}.
\end{equation*}
Noting that $\omega^{-1} \nu \big( u(\omega^{-1} t) \big) = u(t)$, we have
$\omega^{-1} \langle \alpha_i, u(\omega^{-1} t)\rangle = \langle \alpha_{\nu(i)}, u(t)\rangle$,
from which it follows that $\omega^{-1} f_i(\omega^{-1} t)$ is a solution of the Riccati equation $(R^{\nu(i)})$.
So starting from \emph{any} solution $f_k \in \M$ of $(R^k)$ we may define solutions $f_i \in \M$ of $(R^i)$ recursively by
\begin{equation}\label{RecurrenceOrbit}
f_i(t) \coloneqq \omega^{-1}f_{\nu^{-1}(i)}(\omega^{-1}t),
\end{equation}
for $i \in \mathcal I \setminus \{ k \}$. By construction these satisfy the relations in \eqref{cycl cond comm} for all $i \in \mathcal I \setminus \{ \nu^{|\mathcal I|-1}(k) \}$, so we just need to ensure that it also holds for $i = \nu^{|\mathcal I|-1}(k)$, namely
\begin{equation*}
\omega^{-1}f_{\nu^{|\mathcal I| - 1}(k)}(\omega^{-1}t) = f_k(t).
\end{equation*}
Using the recurrence relation \eqref{RecurrenceOrbit} to rewrite the left hand side of this relation we can express the condition for \eqref{new Miura comm} to be cyclotomic as the following functional relation on $f_k$ alone
\begin{equation} \label{f1 relation}
\omega^{-|\mathcal I|} f_k(\omega^{-|\mathcal I|} t) = f_k(t).
\end{equation}
Now since both sides are solutions of the same differential equation $(R^k)$, to check they are equal it suffices to compare their values at a single point. However, to avoid having to explicitly solve $(R^k)$ it would be preferable to evaluate them at the origin since $t=0$ is a fixed point of $t \mapsto \omega^{-|\mathcal I|} t$.

If the original cyclotomic Miura $\g$-oper $\nabla$ happens to be regular at the origin, \emph{i.e.} $\co\lambda_0 = 0$, so that $u(t)$ is also, then the Riccati equation $(R^k)$ is regular at $0$ and therefore admits solutions regular at the origin for any initial value $f_k(0) \in \C$. In particular, two solutions of \eqref{f1 relation} are equal if and only if they agree at the origin. It therefore follows that the gauge transformed Miura $\g$-oper $g \nabla g^{-1}$ is cyclotomic if and only if $\omega^{|\mathcal I|} = 1$, that is to say $|\mathcal I| = T$ since $\omega$ is a primitive $T^{\rm th}$-root of unity.

In order to treat the general case of an arbitrary $\nu$-invariant integral dominant coweight $\co\lambda_0 \in \h^\nu$ we will make use of the following.

\begin{lemma} \label{Riccati singular}
Let $f$ be a non-zero meromorphic solution of the Riccati equation
\begin{equation} \label{reg Riccati lemma}
f'(t) + f(t)^2 + f(t) \left( - \frac{\eta}{t} + r(t) \right) = 0,
\end{equation}
where $r \in \M$ is regular at the origin. If $\eta \in \Z_{\geq 0}$ then either
\begin{itemize}
  \item[$(i)$] $f$ is regular at the origin, in which case it reads
$f(t) = t^\eta h(t)$
for some $h \in \M$ regular at $0$, with $h(0) \in \C^\times$ arbitrary, and satisfying
\begin{equation} \label{reg Riccati def}
h'(t) + t^\eta h(t)^2 + h(t) r(t) = 0.
\end{equation}
We call \eqref{reg Riccati def} the \emph{regularised Riccati equation} associated with \eqref{reg Riccati lemma}.
  \item[$(ii)$] $f$ is singular at the origin, in which case it takes the form $f(t) = \frac{\eta + 1}{t} + k(t)$ for some unique $k \in \M$ regular at $0$.
\qed
\end{itemize}
\begin{proof}
Consider the Laurent expansion of $f \in \M$ at the origin. As in the proof of Lemma \ref{TransfSimple}, for this to be consistent with equation \eqref{reg Riccati lemma} we find that $f$ can have at most a first order pole at the origin, namely
\begin{equation*}
f(t) = \frac{a}{t} + \sum_{n \geq 0} a_n t^n
\end{equation*}
with either $a = 0$ or $a = \eta + 1$.

Let us now discuss how the Ricatti equation \eqref{reg Riccati lemma} constrains the coefficients $a_n$. Taylor expanding $r$ as $r(t)=\sum_{n\geq 0} r_nt^n$, we can write the Laurent expansion of the different terms appearing in this equation as
\begin{align*}
f'(t)  &= -\frac{a}{t^2} + \sum_{n\geq 0} (n+1)a_{n+1}t^n \\
f(t)^2 &= \frac{a^2}{t^2} + \frac{2aa_0}{t} + \sum_{n\geq 0} \left( 2aa_{n+1}+\sum_{k=0}^n a_ka_{n-k} \right) t^n \\
f(t)\left(-\frac{\eta}{t}+r(t) \right) & = -\frac{\eta a}{t^2} + \frac{ar_0-\eta a_0}{t} + \sum_{n\geq 0} \left( a r_{n+1} - \eta a_{n+1} + \sum_{k=0}^n a_kr_{n-k} \right) t^n.
\end{align*}
The sum of these three terms must vanish by virtue of the Ricatti equation \eqref{reg Riccati lemma}. The cancellation of the double pole gives $a=0$ or $a=\eta+1$, as stated above. Those of the simple pole and of the Taylor expansion are equivalent to
\begin{equation*}
\left\lbrace \begin{array}{ll}
(2a-\eta)a_0 = -a r_0, \\
(2a + n - \eta)a_{n} = -a r_{n} - \displaystyle\sum_{k=0}^{n-1} a_k \left(a_{n-1-k}+r_{n-1-k} \right) & \; \; \; \text{ for } n \geq 1.
\end{array} \right.
\end{equation*}

If $a=0$ then $f$ is regular at $0$, corresponding to part $(i)$. The first equation then gives $a_0=0$. For $0 < n < \eta$, the second equation also yields $a_n = 0$ by induction, as $n-\eta \neq 0$. For $n=\eta$, the equation is verified for arbitrary $a_\eta \in \C$. Finally, all $a_n$ for $n > \eta$ are fixed in terms of $a_\eta$ and the coefficients in the Taylor expansion of $r(t)$ at the origin. In particular $f$ is of the form $f(t) = t^\eta h(t)$ for some $h \in \M$ with $h(0) = a_\eta \in \C$ arbitrary. Moreover, by direct substitution into \eqref{reg Riccati lemma} we find that $h$ satisfies \eqref{reg Riccati def}.

On the other hand, if $a = \eta + 1$ then $f$ is singular at $0$, which corresponds to part $(ii)$. In this case the first equation reads $a_0=-\frac{\eta+1}{\eta+2}r_0$. By induction on $n >0$, the second one fixes uniquely the coefficients $a_n$ in terms of the coefficients $r_k$, as for any $n>0$ the coefficient $2a+n-\eta=n+\eta+2$ is non-zero.
\end{proof}
\end{lemma}

\subsubsection{Singular reproduction procedure}
Since $\co\lambda_0$ is assumed integral dominant we may apply Lemma \ref{Riccati singular} to the solution $f_k$ of $(R^k)$, with $\eta = \langle \alpha_k, \co\lambda_0 \rangle$. If $f_k$ has a pole at the origin then so does the function $t \mapsto \omega^{-|\mathcal I|} f_k(\omega^{-|\mathcal I|} t)$. However, since both are solutions to the same Riccati equation of the form \eqref{reg Riccati lemma}, by the uniqueness of such a solution in Lemma \ref{Riccati singular}$(ii)$ it follows that they are equal, \emph{i.e.} the condition \eqref{f1 relation} holds.

Moreover, in this case, all the functions $f_i$ for $i\in \mathcal{I}$ have a pole at the origin. Thus, using the explicit form $g = \prod_{i \in \mathcal I} e^{f_i E_i}$ of the gauge transformation parameter and applying Lemma \ref{TransfSimple} for each $i\in\mathcal{I}$, followed by Lemma \ref{lem: Weyl nu}, we find that
\begin{equation*}
-\res_0 \overline{g \nabla g^{-1}} = \left( \prod_{i \in \mathcal I} s_i \right) \cdot (-\res_0 \overline{\nabla}) = s^\nu_{\mathcal I} \cdot (-\res_0 \overline{\nabla}).\vspace{5pt}
\end{equation*}

\subsubsection{Regular reproduction procedure}
Suppose now that $f_k$ is regular at the origin, corresponding to case $(i)$ of Lemma \ref{Riccati singular}. It remains to check when the condition \eqref{f1 relation} holds. Introducing the function
\begin{equation*}
h_k(t) \coloneqq t^{- \eta} f_k(t),
\end{equation*}
regular at the origin, the condition \eqref{f1 relation} for $g \nabla g^{-1}$ to be cyclotomic may then be rewritten as
\begin{equation} \label{h1 relation}
\omega^{- |\mathcal I| (\eta+1)} h_k(\omega^{-|\mathcal I|}t) = h_k(t).
\end{equation}
Now since both sides of \eqref{h1 relation} satisfy the same regularised Riccati equation of the form \eqref{reg Riccati def}, they are equal if and only if they agree at the origin. Note that we may as well assume $h_k(0) \neq 0$ since if $h_k(0) = 0$ then by the homogeneity of the regularised Riccati equation $h_k$ would be identically zero and hence $g$ the identity. It follows that $g \nabla g^{-1}$ is cyclotomic if and only if $\omega^{|\mathcal I|(\eta + 1)} = 1$. Recalling that $\omega$ is a primitive $T^{\rm th}$-root of unity, this condition is equivalent to $|\mathcal I| (\eta + 1) \equiv 0 \; \textup{mod} \; T$. Or since the size $|\mathcal I|$ of the orbit $\mathcal I$ necessarily divides the order of $\nu$ which by definition divides $T$, we can also rewrite this condition as
\begin{equation*}
\langle \alpha_k, \co\lambda_0 + \co\rho \rangle \equiv 0 \;\;\textup{mod}\;\; \frac{T}{|\mathcal I|}.
\end{equation*}
Note that, in this case, all the functions $f_i$ are regular at the origin, so that $\res_0 \overline{g\nabla g^{-1}}=\res_0 \overline{\nabla}$.

\subsubsection{Residue at infinity}
Finally, let us discuss the residue at infinity of the new connection $g\nabla g^{-1}$. We suppose that $\langle \alpha_k, \res_\infty \overline{\nabla} +\co\rho \rangle$ (and thus all the $\langle \alpha_i, \res_\infty \overline{\nabla} +\co\rho \rangle$ for $i\in\mathcal{I}$) is non-negative. One can see the gauge transformation by $g$ as successive gauge transformations by $e^{f_iE_i}$, for each $i\in \mathcal{I}$. As $\langle \alpha_i, \co\alpha_j \rangle =0$ for any distinct $i,j\in\mathcal{I}$, the condition on the residue at infinity of the connection still holds after each step. Thus, applying the last part of Lemma \ref{TransfSimple} for each $i\in\mathcal{I}$, followed by Lemma \ref{lem: Weyl nu}, we find that
\begin{equation*}
\res_\infty \overline{g \nabla g^{-1}} = \left( \prod_{i \in \mathcal I} s_i \right) \cdot (\res_\infty \overline{\nabla}) = s^\nu_{\mathcal I} \cdot (\res_\infty \overline{\nabla}).
\end{equation*}
Note that this discussion holds whether $f_k$ is regular or singular at $0$. We have proved the following.

\begin{theorem}\label{CommutingOrbit}
Let $\mathcal I \in I/\nu$ be an orbit of type $A_1^{\times |\mathcal I|}$.
Fix a $k \in \mathcal I$ and let $f_k \in \M$ be any non-zero solution of the Riccati equation
\begin{equation*}
f_k'(t) + f_k(t)^2 + f_k(t) \langle \alpha_k, u(t) \rangle = 0,
\end{equation*}
with $u(t)$ as in \eqref{cycl Miura lambda0}.
Define $f_i \in \M$ for all $i \in \mathcal I \setminus \{ k \}$ recursively by $f_i(t) \coloneqq \omega^{-1}f_{\nu^{-1}(i)}(\omega^{-1}t)$ and let $g = \prod_{i \in \mathcal I} e^{f_i E_i} \in N(\M)$.

If $f_k$ is regular at the origin, then $g \in N^{\hat\vs}(\M)$, and hence $g \nabla g^{-1} \in \MOp^\Gamma_\g(\CP)$, if and only if
\begin{equation*}
\langle \alpha_k, \co\lambda_0 + \co\rho \rangle \equiv 0 \;\;\textup{mod}\;\; \frac{T}{|\mathcal I|}.
\end{equation*}
Moreover, in this case
we have $\res_0 \overline{g \nabla g^{-1}} = \res_0 \overline{\nabla}$.

If $f_k$ has a pole at the origin, then we have $g \nabla g^{-1} \in \MOp^\Gamma_\g(\CP)$ without any condition on $\co\lambda_0$ and $-\res_0 \overline{g \nabla g^{-1}} = s^\nu_{\mathcal I} \cdot (-\res_0 \overline{\nabla})$.

Finally, if $\langle \alpha_k, \res_\infty \overline{\nabla} + \co\rho \rangle$ is non-negative, then $\res_\infty \overline{g \nabla g^{-1}} = s^\nu_{\mathcal I} \cdot (\res_\infty \overline{\nabla})$ (whether the function $f_k$ is regular or singular at the origin).
\qed
\end{theorem}
 
\begin{rem}
It is interesting to note that, just as in the case of an ordinary reproduction, we have here a one-parameter family of transformations (the parameter being the initial value of $h_k$ at the origin). In other words, even though we performed $|\mathcal I|$ different successive gauge transformations, we don't have $|\mathcal I|$ free parameters since the $\Gamma$-equivariance imposes relations between them.
\end{rem}

\begin{exmp}
Let us illustrate Theorem \ref{CommutingOrbit} with an example. Consider the Lie algebra $\g=\sl_4$ and the cyclotomic Miura $\sl_4$-oper $\nabla$ described in \eqref{ExpleSl4}. We use the notations and conventions introduced in the example \ref{Sl4}. To perform a reproduction in the direction of the orbit $\lbrace\alpha_1,\alpha_3\rbrace$, we use a gauge transformation by the unipotent element $g = \exp\bigl(f_3 E_3\bigr)\exp\bigl(f_1 E_1\bigr)$. The functions $f_1$ and $f_3$ must satisfy the Ricatti equations
\begin{equation*}
(R^i) : \; \; \; f_i'(t) + f_i(t)^2 + q_i(t) f_i(t) = 0,
\end{equation*}
with
\begin{equation*}
q_1(t) \coloneqq - \frac{\eta}{t} - \frac{St^{S-1}}{t^S-z^S}, \qquad q_3(t) \coloneqq - \frac{\eta}{t} - \frac{St^{S-1}}{t^S+z^S}.
\end{equation*}
Introducing the functions
\begin{equation*}
Q_1(t) \coloneqq e^{-\int q_1(t) dt} = t^\eta (t^S-z^S), \qquad  R_1(t) \coloneqq \int Q_1(t) dt = \frac{1}{\eta+S+1}t^{\eta+S+1}-\frac{z^S}{\eta+1}t^{\eta+1},
\end{equation*}
one can check that the most general solution of $(R^1)$ is
\begin{equation*}
f_1(t) = \frac{Q_1(t)}{R_1(t)+\widetilde{A}} = \frac{(\eta+1)(S+\eta+1)t^\eta(t^S-z^S)}{(\eta+1)t^{S+\eta+1}-(S+\eta+1)z^St^{\eta+1}+A},
\end{equation*}
where $A \coloneqq (S+\eta+1)(\eta+1)\widetilde{A}$ is an arbitrary integration constant. Following Theorem \ref{CommutingOrbit}, we introduce
\begin{equation*}
f_3(t) \coloneqq \omega^{-1} f_1(\omega^{-1}t) = \frac{(\eta+1)(S+\eta+1)t^\eta(t^S+z^S)}{(\eta+1)t^{S+\eta+1}+(S+\eta+1)z^St^{\eta+1} - \omega^{\eta+1} A},
\end{equation*}
where in the second equality we used the fact that $\omega^{-S} = -1$.
One can check that $f_3$ is a solution of the Ricatti equation $(R^3)$. The resulting unipotent element $g$ is cyclotomic if and only if
\begin{equation*}
f_1(t) = \omega^{-1} f_3(\omega^{-1}t) = \frac{(\eta+1)(S+\eta+1)t^\eta(t^S-z^S)}{(\eta+1)t^{S+\eta+1}-(S+\eta+1)z^St^{\eta+1} + \omega^{2(\eta+1)} A}.
\end{equation*}
If $A = 0$ then the above condition always holds. On the other hand, if $A \neq 0$ then this condition is equivalent to $\omega^{2(\eta+1)}=1$, \textit{i.e.} to
\begin{equation*}
\langle \alpha_1, \co\lambda_0 + \co\rho \rangle \equiv \eta+1 \equiv 0 \text{ mod } \frac{T}{2},
\end{equation*}
in agreement with Theorem \ref{CommutingOrbit}.

Let us also determine the residues at the origin and infinity of the connection $\overline{g\nabla g^{-1}}$ obtained after reproduction. We have seen in example \ref{Sl4} that $\res_\infty \overline{\nabla} = \left( \eta + S \right)(\co\omega_1+\co\omega_3) + \kappa\co\omega_2$. After the reproduction procedure, we get
\begin{equation*}
\res_\infty \overline{g\nabla g^{-1}} = \left( \eta + S \right)(\co\omega_1+\co\omega_3) + \kappa\co\omega_2 - \left( \eta + S + 1 \right) (\co\alpha_1+\co\alpha_3) = s_1s_3 \cdot\left( \res_\infty \overline{\nabla} \right).
\end{equation*}
When $f_1$ and $f_3$ are both regular at the origin (\textit{i.e.} when $A\neq 0$), the residue at the origin does not change after the reproduction procedure: $\res_0 \overline{g\nabla g^{-1}} = \res_0 \overline{\nabla} = - \eta (\co\omega_1+\co\omega_3) - \kappa \co\omega_2$. However, if $f_1$ and $f_3$ are singular at the origin (\textit{i.e.} $A=0$), we find
\begin{equation*}
-\res_0 \overline{g\nabla g^{-1}} =  \eta (\co\omega_1+\co\omega_3) + \kappa \co\omega_2 - (\eta+1)(\co\alpha_1+\co\alpha_3) = s_1s_3 \cdot\left( -\res_0 \overline{\nabla} \right). \qedhere
\end{equation*}
\end{exmp}

\subsection{Orbit of type $A_2^{\times |\mathcal I|/2}$} \label{sec: A2 orbit}

Suppose now that $\mathcal I \in I/\nu$ is an orbit of type $A_2^{\times |\mathcal I|/2}$ with $|\mathcal I|$ even, corresponding to the case $\ell_{\mathcal I} = 2$, cf. \S\ref{folding}.
In the setting of simple Lie algebras of finite type, such orbits occur only in type $A_{2n}$, in which case the orbit is of size $|\mathcal I|=2$.

Let us fix a reference point $k \in \mathcal I$ and denote one half of the orbit by
\begin{equation*}
\mathcal I/2 = \{ k, \nu(k), \ldots, \nu^{|\mathcal I|/2 - 1}(k) \},
\end{equation*}
cf. Lemma \ref{lem: Weyl nu}. Recall the notation $\bar\imath = \nu^{|\mathcal I|/2}(i)$. The full orbit $\mathcal I$ then consists of $|\mathcal I|/2$ distinct pairs $\{ i, \bar\imath \}$ for each $i \in \mathcal I/2$ such that $\langle \alpha_i, \co\alpha_{\bar\imath} \rangle = \langle \alpha_{\bar\imath}, \co\alpha_i \rangle = -1$ and $\langle \alpha_i, \co\alpha_j \rangle = 0$ for any $j \in \mathcal I$ with $j \neq i, \bar\imath$.
For each $i \in \mathcal I/2$, the subalgebra generated by $E_i$ and $E_{\bar\imath}$ has dimension 3 and is spanned by $E_i$, $E_{\bar\imath}$ and $[E_i,E_{\bar\imath}]$.

In this case, performing two successive reproductions in the directions of $\alpha_i$ and $\alpha_{\bar\imath}$ is not enough to restore the $\Gamma$-equivariance of the Miura $\g$-oper. In order to obtain a cyclotomic Miura $\g$-oper we will need to perform a third reproduction in the direction of $\alpha_i$, so that the total gauge transformation parameter is of the form $g = \prod_{i \in \mathcal I/2} e^{k_{i, 3} E_i} e^{k_{i, 2} E_{\bar\imath}} e^{k_{i, 1} E_i}$ for some $k_{i,1}, k_{i,2}, k_{i,3} \in \M$.
Equivalently, we can write the latter as 
 \begin{equation}\label{TransfoNonCommOrbit}
g \coloneqq \prod_{i \in \mathcal I /2} e^{f_{i, 1} (E_i + E_{\bar\imath}) + f_{i,2} (E_i - E_{\bar\imath}) + f_{i,3} [E_i, E_{\bar\imath}]}
\end{equation}
for some $f_{i,1}, f_{i,2}, f_{i,3} \in \M$, where we have used the combinations $E_i \pm E_{\bar\imath}$ rather than $E_i$ and $E_{\bar\imath}$ for later convenience. Note that the order in the above product over $\mathcal I/2$ does not matter since the generators $E_i, E_{\bar\imath}$ commute with $E_j, E_{\bar\jmath}$ for any distinct $i, j \in \mathcal I/2$.
The calculation can be performed using either of the above expressions for $g$. If we consider the product of exponentials of simple root generators then we have to study three Riccati equations but, since now $(\alpha_i | \alpha_{\bar\imath}) \neq 0$, the argument used in \S\ref{sec: A1 orbit} no longer applies. We will use the second form \eqref{TransfoNonCommOrbit} for $g$.

Let $\nabla = d + (p_{-1} + u) dt \in \MOp^\Gamma_\g(\CP)$. A lengthy but straightforward computation gives the expression of the gauge transformed $\g$-connection
\begin{align*}
g \nabla g^{-1} = \nabla + \sum_{i \in \mathcal I/2} \Big( &(f_{i,1}+f_{i,2}) \co\alpha_i + (f_{i,1}-f_{i,2}) \co\alpha_{\bar\imath}\\[-3mm]
& - \ha \big(f_{i,1}^2 + 3f_{i,2}^2 + (q_i + q_{\bar\imath})f_{i,1} + (q_i - q_{\bar\imath})f_{i,2} + 2f'_{i,1} \big) (E_i + E_{\bar\imath})  \\
 & - \ha \big(4 f_{i,1} f_{i,2} + (q_i + q_{\bar\imath})f_{i,2} + (q_i - q_{\bar\imath})f_{i,1} - 2f_{i,3} + 2f'_{i,2} \big) (E_i-E_{\bar\imath})  \\
 & - \big( 2f_{i,1} f_{i,3} - f_{i,1}f'_{i,2} + f_{i,2}f'_{i,1}+f_{i,2}(f_{i,2}^2-f_{i,1}^2) + f'_{i,3}\\
 &\qquad\qquad\qquad\qquad + (q_i + q_{\bar\imath})f_{i,3} - \ha (q_i - q_{\bar\imath})(f_{i,1}^2-f_{i,2}^2) \big) [E_i, E_{\bar\imath}] \Big) dt,
\end{align*}
where $q_i(t) \coloneqq \langle \alpha_i, u(t) \rangle$ and $q_{\bar\imath}(t) \coloneqq \langle \alpha_{\bar\imath}, u(t) \rangle$.

Since we want $g \nabla g^{-1}$ to be a Miura $\g$-oper we should impose that the coefficients of $E_i \pm E_{\bar\imath}$ and $[E_i, E_{\bar\imath}]$ all vanish. This gives rise to a set of three coupled Riccati-type differential equations, which after some rearranging take the form
\begin{equation} \label{RiccatiNonCommuting}
(R^i): \qquad \left\{
\begin{split}
2f'_{i,1} + f_{i,1}^2 + 3f_{i,2}^2 + (q_i + q_{\bar\imath})f_{i,1} + (q_i - q_{\bar\imath})f_{i,2} &= 0, \\
2f'_{i,2} + 4f_{i,1}f_{i,2} - 2f_{i,3} + (q_i + q_{\bar\imath})f_{i,2} + (q_i - q_{\bar\imath})f_{i,1} &= 0,  \\
2f'_{i,3} + 2f_{i,1}f_{i,3} + f_{i,2}(f_{i,1}^2-f_{i,2}^2) + 2(q_i + q_{\bar\imath})f_{i,3} &= 0.
\end{split}
\right.
\end{equation}

Noting that $\omega^{-1} \nu\big( u(\omega^{-1} t) \big) = u(t)$ since $u(t) dt \in \Omega^{\hat\nu}(\h)$, which also implies that for every $i \in \mathcal I / 2$ we have $\omega^{-1} q_{\nu^{-1}(i)}(\omega^{-1} t) = q_i(t)$ and $\omega^{-1} q_{\nu^{-1}(\bar\imath)}(\omega^{-1} t) = q_{\bar\imath}(t)$, we deduce that the three functions $\omega^{-1} f_{i, 1}(\omega^{-1} t)$, $\omega^{-1} f_{i, 2}(\omega^{-1} t)$ and $\omega^{-2} f_{i, 3}(\omega^{-1} t)$ satisfy the coupled system of equations $(R^{\nu(i)})$. Therefore, starting with \emph{any} solution $f_{k,1}, f_{k,2}, f_{k,3} \in \M$ for the reference point $k \in \mathcal I$ we can define a solution $f_{i, 1}, f_{i,2}, f_{i,3} \in \M$ of $(R^i)$ recursively by
\begin{equation}\label{Recursion}
f_{i,1}(t) \coloneqq \omega^{-1} f_{\nu^{-1}(i),1}(\omega^{-1} t), \qquad
f_{i,2}(t) \coloneqq \omega^{-1} f_{\nu^{-1}(i),2}(\omega^{-1} t), \qquad
f_{i,3}(t) \coloneqq \omega^{-2} f_{\nu^{-1}(i),3}(\omega^{-1} t).
\end{equation}
for every $i \in \mathcal I/2 \setminus \{ k \}$. Now since $\nu^{|\mathcal I|/2}(k) = \bar k$ we have
\begin{equation*}
\vs\big( E_{\nu^{|\mathcal I|/2 -1}(k)} \pm E_{\nu^{|\mathcal I|/2 -1}(\bar k)} \big) = \pm \omega^{-1}(E_k \pm E_{\bar k}), \qquad
\vs\big( \big[ E_{\nu^{|\mathcal I|/2 -1}(k)}, E_{\nu^{|\mathcal I|/2 -1}(\bar k)} \big] \big) = - \omega^{-2} [E_k, E_{\bar k}].
\end{equation*}
So by construction, we find that $g \in N(\M)$ defined in \eqref{TransfoNonCommOrbit} using the functions $f_{i,1}, f_{i,2}, f_{i,3}$ for all $i \in \mathcal I/2$ is invariant under $\hat\vs$ if and only if
\begin{gather*}
f_{k,1}(t) = \omega^{-1} f_{\nu^{|\mathcal I|/2 -1}(k),1}(\omega^{-1} t), \qquad
f_{k,2}(t) = - \omega^{-1} f_{\nu^{|\mathcal I|/2 -1}(k),2}(\omega^{-1} t), \\
f_{k,3}(t) = - \omega^{-2} f_{\nu^{|\mathcal I|/2 -1}(k),3}(\omega^{-1} t).
\end{gather*}
Equivalently, using the above recurrence relations to rewrite the right hand sides of these equations we obtain conditions on the functions $f_{k,1}$, $f_{k,2}$ and $f_{k,3}$ alone, namely
\begin{equation} \label{f123 cycl relations}
\begin{split}
f_{k,1}(t) &= \omega^{-|\mathcal I|/2} f_{k,1}(\omega^{-|\mathcal I|/2} t), \quad
f_{k,2}(t) = - \omega^{-|\mathcal I|/2} f_{k,2}(\omega^{-|\mathcal I|/2} t), \\
&\qquad\qquad\qquad f_{k,3}(t) = - \omega^{-|\mathcal I|} f_{k,3}(\omega^{-|\mathcal I|/2} t).
\end{split}
\end{equation}

We now consider the issue of the regularity of the equations \eqref{RiccatiNonCommuting} at the origin. The functions $q_i$ and $q_{\bar\imath}$ have the form $q_i(t) = - \frac{\eta}{t} + r_i(t)$ and $q_{\bar\imath}(t) = - \frac{\eta}{t} + r_{\bar\imath}(t)$, with $\eta = \langle \alpha_i, \co\lambda_0 \rangle = \langle \alpha_{\bar\imath}, \co\lambda_0 \rangle$ and $r_i, r_{\bar\imath} \in \h(\M)$ both regular at the origin. Note that $\eta$ is independent of $i$ and $\bar\imath$ by the $\nu$-invariance of the coweight $\co\lambda_0$ and is a non-negative integer, as $\co\lambda_0$ is a dominant integral coweight. We will need the following result, which is the analog, for an orbit of type $A_2^{\times|\mathcal{I}|/2}$, of Lemma \ref{Riccati singular}.

\begin{lemma}\label{LemmaRegNonCom}
Let $f_1$, $f_2$ and $f_3$ be meromorphic solutions of the system of differential equations
\begin{subequations} \label{RiccatiNonCom2}
\begin{align}
2f'_{1} + f_{1}^2 + 3f_2^2 + (q + \overline{q})f_{1} + (q - \overline{q})f_{2} &= 0, \\
2f'_{2} + 4f_{1}f_{2} - 2f_{3} + (q + \overline{q})f_{2} + (q - \overline{q})f_{1} &= 0,  \\
2f'_{3} + 2f_{1}f_{3} + f_{2}(f_1^2-f_{2}^2) + 2(q + \overline{q})f_3 &= 0.
\end{align}
\end{subequations}
Suppose $q(t)=-\frac{\eta}{t}+r(t)$ and $\overline{q}(t)=-\frac{\eta}{t}+\overline{r}(t)$, with $r,\overline{r}\in\M$ regular at $0$. If $\eta\in\Z_{\geq0}$ then either
\begin{enumerate}
\item[$(i)$] the functions $f_i$ are regular at $0$, in which case they read $f_1(t)=t^\eta h_1(t)$, $f_2(t)=t^\eta h_2(t)$ and $f_3(t)=t^{2\eta}h_3(t)$ for some $h_i \in \M$ regular at $0$, with $h_i(0) \in\C$ arbitrary. 

\item[$(ii)$] at least one of the functions $f_i$ is singular at $0$, in which case they are of the form
\begin{equation} \label{SolSing}
f_1(t) = \frac{a}{t}+k_1(t), \;\;\; f_2(t)=\frac{b}{t}+k_2(t), \;\;\; f_3(t) = \frac{c}{t^2} + \frac{d}{t} + k_3(t),
\end{equation}
where $k_i \in \M$ are regular at 0 and $a,b,c,d \in \C$. Moreover, the coefficients $a,b,c$ fall into one of two possible classes
\begin{itemize}
\item[$(a)$] $(a,b,c)=\bigl(2(\eta+1),0,0\bigr)$, in which case $d\in\C$ and $k_i\in\M$ are unique,
\item[$(b)$] $(a,b,c)=\left(\frac{\eta+1}{2},\pm \frac{\eta+1}{2},0 \right)$ or $(a,b,c)=\left( \frac{3(\eta+1)}{2},\pm\frac{\eta+1}{2},\pm(\eta+1)^2 \right)$.
\end{itemize}
\end{enumerate}
\begin{proof}
Consider the Laurent expansions of the functions $f_i$ around the origin. The consistency of these expansions with the equations \eqref{RiccatiNonCom2} requires $f_1$ and $f_2$ to have at most a simple pole at $0$ and $f_3$ at most a double pole at $0$. Thus the $f_i$'s can be written in the form of equation \eqref{SolSing}. Moreover, one finds that the cancellation of the double poles in the first two equations and of the triple pole in the third one leads to the set of conditions
\begin{equation}\label{Condabc}
a^2-2(\eta+1)a+3b^2=0, \;\;\;\; 2ab-(\eta+1)b-c=0, \;\;\;\; 2ac+b(a^2-b^2)-4(\eta+1)c=0.
\end{equation}
There are 6 solutions to this system of equations: $(a,b,c)=(0,0,0)$, $(a,b,c)=\bigl(2(\eta+1),0,0\bigr)$, $(a,b,c)=\left(\frac{\eta+1}{2},\pm \frac{\eta+1}{2},0 \right)$ and $(a,b,c)=\left( \frac{3(\eta+1)}{2},\pm\frac{\eta+1}{2},\pm(\eta+1)^2 \right)$.

In the case where $(a,b,c)=(0,0,0)$ one finds that $d$ is also zero. This corresponds to case $(i)$ of the lemma, where the functions $f_i$ are all regular at the origin. Consider then the Taylor expansions $f_1(t)=\sum_{n\geq 0} a_nt^n$, $f_2(t)=\sum_{n\geq 0} b_nt^n$ and $f_3(t)=\sum_{n\geq 0} c_nt^n$. Substituting these into \eqref{RiccatiNonCom2} and working order by order in powers of $t$ one finds that $a_0=b_0=c_0=0$ and that the coefficients $a_n$, $b_n$ and $c_n$ for $n\in\Z_{\geq 1}$ must satisfy the recurrence relations
\begin{equation*}
\left\{
\begin{split}
2(n+1-\eta)a_{n+1} &= - \sum_{k=0}^n \Bigl( a_k(a_{n-k}+r_{n-k}+\overline{r}_{n-k}) + b_k(3b_{n-k}+r_{n-k}-\overline{r}_{n-k}) \Bigr), \\
2(n+1-\eta)b_{n+1} &= 2c_n - \sum_{k=0}^n \Bigl( a_k(4b_{n-k} + r_{n-k} - \overline{r}_{n-k}) + b_k(r_{n-k}+\overline{r}_{n-k}) \Bigr), \\
2(n+1-2\eta)c_{n+1} &= - \sum_{k=0}^n c_k(a_{n-k}+2r_{n-k}+2\overline{r}_{n-k}) + \sum_{k_1+k_2+k_3=n}\bigl( b_{k_1}b_{k_2}b_{k_3} - b_{k_1}a_{k_2}a_{k_3} \bigr),
\end{split}
\right.
\end{equation*}
where $r_n$ and $\overline{r}_n$ are coefficients in the Taylor expansion of $r$ and $\overline{r}$.

We have $a_0=b_0=c_0=0$. For $n+1 < \eta$, the coefficients $n+1-\eta$ and $n+1-2\eta$ in the left hand sides of these equations are non-zero, so we conclude by induction that $a_n=b_n=c_n=0$ for $n < \eta - 1$. When $n=\eta-1$, the first two equations are verified for arbitrary values of $a_{n+1} = a_\eta \in \C$ and $b_{n+1} = b_\eta \in \C$. Similarly, the third equation yields $c_n=0$ for all $n<2\eta - 1$ and $c_{2\eta} \in \C$ is arbitrary. The coefficients $a_n$ and $b_n$ for $n>\eta$ and $c_n$ for $n>2\eta$ are then determined uniquely in terms of the $r_n$'s, the $\overline{r}_n$'s and the arbitrary coefficients $a_\eta$, $b_\eta$ and $c_{2\eta}$. This conclude the discussion of case $(i)$ of the lemma.

Next, we turn to the 5 other solutions of equation \eqref{Condabc}, corresponding to the case where at least one of the $f_i$'s is singular, \emph{i.e.} the case $(ii)$ of the lemma. We shall focus on case $(a)$, with $(a,b,c)=\bigl(2(\eta+1),0,0\bigr)$. Substituting the Taylor expansions of the regular functions $k_i$ into \eqref{RiccatiNonCom2}, one finds unique expressions for the coefficients $d$, $a_0$, $b_0$ and $c_0$. After a tedious but straightforward computation, one finds recurrence relations similar to the ones above. In the case at hand, however, the coefficients in front of $a_{n+1}$, $b_{n+1}$ and $c_{n+1}$ are non-zero for any $n\in\Z_{\geq 0}$. Thus, all the coefficients $a_n$, $b_n$ and $c_n$ are determined uniquely by these recursion relations and the initial values $a_0$, $b_0$ and $c_0$ previously obtained. Therefore, in this case, the functions $k_i$ are unique.
\end{proof}

\end{lemma}

Let us now come back to the discussion of cyclotomic reproductions along the orbit $\mathcal{I}$. We began with solutions $f_{k,1}$, $f_{k,2}$ and $f_{k,3}$ of equations $(R^k)$. As explained above, we constructed solutions $f_{i,1}$, $f_{i,2}$ and $f_{i,3}$ of equations $(R^i)$ for any $i\in\mathcal{I}/2$ recursively using the formula \eqref{Recursion}. The condition for the corresponding reproduction to be cyclotomic is then given by equation \eqref{f123 cycl relations}.

\subsubsection{Regular reproduction procedure}

Suppose that the functions $f_{k,1}$, $f_{k,2}$ and $f_{k,3}$ (and thus also the functions $f_{i,1}$, $f_{i,2}$ and $f_{i,3}$ for any $i\in\mathcal{I}/2$) are regular at the origin, so that they are described by case $(i)$ of Lemma \ref{LemmaRegNonCom}. Then $h_{i,1}(t) \coloneqq t^{-\eta} f_{i,1}(t)$, $h_{i,2}(t) \coloneqq t^{-\eta} f_{i,2}(t)$ and $h_{i,3}(t) \coloneqq t^{-2\eta} f_{i,3}(t)$ are regular at the origin. They obey the following regularised equations
\begin{subequations}\label{RegNonComm}
\begin{align}
\label{RRA} 2h'_{i,1} + t^\eta(h_{i,1}^2 + 3h_{i,2}^2) + (r_i + r_{\bar\imath})h_{i,1} + (r_i - r_{\bar\imath})h_{i,2} &= 0, \\
\label{RRB} 2h'_{i,2} + t^\eta(4h_{i,1} h_{i,2}-2h_{i,3}) + (r_i + r_{\bar\imath})h_{i,2} + (r_i - r_{\bar\imath})h_{i,1} &= 0,\\
\label{RRC} 2h'_{i,3} + t^\eta\big( 2h_{i,1} h_{i,3} + h_{i,2}(h_{i,1}^2-h_{i,2}^2) \big) + 2(r_i + r_{\bar\imath})h_{i,3} &= 0,
\end{align}
\end{subequations}
whose coefficients are regular at $0$ since $\eta \in \Z_{\geq 0}$. The coupled system of equations \eqref{RegNonComm} admits a solution regular at $0$ for every set of initial conditions $h_{i,1}(0), h_{i,2}(0), h_{i,3}(0) \in \C$. It remains to check whether there are solutions for which the corresponding gauge transformation parameter $g$ defined in \eqref{TransfoNonCommOrbit} is $\hat\vs$-invariant so that $g \nabla g^{-1} \in \MOp^\Gamma_\g(\CP)$.

When phrased in terms of the regularised functions $h_{k, 1}$, $h_{k,2}$ and $h_{k,3}$, the conditions \eqref{f123 cycl relations} for $g \in N(\M)$ to be $\hat\vs$-invariant read
\begin{equation} \label{h123 cycl rel}
\begin{split}
h_{k,1}(t) &= \omega^{-(1+\eta) |\mathcal I| / 2} h_{k,1}(\omega^{-|\mathcal I| / 2}t), \qquad
h_{k,2}(t) = -\omega^{-(1+\eta) |\mathcal I| / 2} h_{k,2}(\omega^{-|\mathcal I| / 2}t),\\
&\qquad\qquad\qquad h_{k,3}(t) = -\omega^{-(1+\eta) |\mathcal I|} h_{k,3}(\omega^{-|\mathcal I|/2}t).
\end{split}
\end{equation}
Now by construction, both sides of the relations \eqref{h123 cycl rel} satisfy the same system of coupled equations \eqref{RegNonComm}. Thus \eqref{h123 cycl rel} hold for all $t \in \C$ if and only if they hold at the origin. It follows that $g \in N^{\hat\vs}(\M)$ if and only if
\begin{equation} \label{h1h2h3 origin}
\begin{split}
h_{k,1}(0) &= \omega^{-(1+\eta) |\mathcal I| / 2} h_{k,1}(0), \qquad
h_{k,2}(0) = -\omega^{-(1+\eta) |\mathcal I| / 2} h_{k,2}(0), \\
&\qquad\qquad h_{k,3}(0) = -\omega^{-(1+\eta) |\mathcal I|} h_{k,3}(0).
\end{split}
\end{equation}
Since the coupled equations \eqref{RegNonComm} are homogeneous, for the set of functions $h_{k,1}$, $h_{k,2}$ and $h_{k,3}$ to be non-trivial, and hence $g \neq \Id$, we should certainly have at least one of the initial conditions $h_{k,1}(0)$, $h_{k,2}(0)$ and $h_{k,3}(0)$ be non-zero. However, the above conditions \eqref{h1h2h3 origin} are consistent only if at most one of these initial conditions is non-zero. We therefore have three distinct possibilities:
\begin{itemize}
\item $\omega^{(1+\eta) |\mathcal I| / 2}=1$ \quad if \; $h_{k,2}(0)=h_{k,3}(0)=0$,
\item $\omega^{(1+\eta) |\mathcal I| / 2}=-1$ \quad if \; $h_{k,1}(0)=h_{k,3}(0)=0$,
\item $\omega^{(1+\eta) |\mathcal I|}=-1$ \quad if \; $h_{k,1}(0)=h_{k,2}(0)=0$.
\end{itemize}
In each case, the $\hat\vs$-invariance of $g$ reduces the number of free initial conditions to one. Therefore, we have once again a one-parameter family of cyclotomic gauge transformations. We note that the above three cases can be combined into the single condition $\omega^{2 |\mathcal I| (1+\eta)}=1$. Or since $\omega$ is a primitive $T^{\rm th}$-root of unity this is equivalent to $2 |\mathcal I| (1 + \eta) \equiv 0 \; \textup{mod} \; T$. Using the definition of $\eta$ we can also rewrite this as $|\mathcal I| \langle \alpha_k + \alpha_{\bar k}, \co\lambda_0 + \co\rho \rangle \equiv 0 \; \textup{mod} \; T$ and in turn, since $|\mathcal I|$ divides $T$, as
\begin{equation*}
\langle \alpha_k + \alpha_{\bar k}, \co\lambda_0 + \co\rho \rangle \equiv 0 \; \textup{mod} \; \frac{T}{|\mathcal I|}.
\end{equation*}

Since all the functions $f_{i,1}$ and $f_{i,2}$ are regular at $0$, it is clear that this cyclotomic reproduction does not change the residue at the origin, namely
\begin{equation*}
\res_0 \overline{g\nabla g^{-1}} = \res_0 \overline{\nabla} = - \co\lambda_0. \vspace{8pt}
\end{equation*}

\subsubsection{Singular reproduction procedure}

Suppose now that at least one of the functions $f_{k,1}$, $f_{k,2}$ and $f_{k,3}$ is singular at $0$. They are then described by case $(ii)$ of Lemma \ref{LemmaRegNonCom} and are of the form \eqref{SolSing}. Let us see when these solutions correspond to a cyclotomic reproduction, \emph{i.e.} when they verify the functional relation \eqref{f123 cycl relations}. It is clear that this relation is compatible with the form \eqref{SolSing} only if the coefficients $b$ and $c$ are zero. According to Lemma \ref{LemmaRegNonCom}, this only leaves case $(a)$ as a possibility, where $(a,b,c)=(2(\eta+1),0,0)$.

By construction, the sets of functions on the left and right hand sides of \eqref{f123 cycl relations} are solutions of the same equation $(R^k)$. Moreover, they are both of the form given in case $(ii)(a)$ of Lemma \ref{LemmaRegNonCom}. By uniqueness of the solutions of this form, these two sets of functions must be equal and thus the corresponding reproduction procedure is cyclotomic.

Finally, recall that the new $\g$-connection is given by
\begin{equation*}
g\nabla g^{-1} = \nabla + \sum_{i\in\mathcal{I}/2} \Bigl( f_{i,1} \left(\co\alpha_i+\co\alpha_{\bar\imath}\right) + f_{i,2}\left(\co\alpha_i-\co\alpha_{\bar\imath}\right) \Bigr)dt.
\end{equation*}
Using the fact that
\begin{equation*}
\res_0 f_{i,1} dt = 2(\eta+1) = \langle \alpha_i + \alpha_{\bar\imath}, \co\lambda_0 + \co\rho \rangle, \qquad
\res_0 f_{i,2} dt = 0 \quad \text{and} \quad \res_0 f_{i,3} dt = 0,
\end{equation*}
and Lemma \ref{lem: Weyl nu}, one finds
\begin{equation*}
-\res_0 \overline{g \nabla g^{-1}} = \left( \prod_{i \in \mathcal I / 2} s_i s_{\bar\imath} s_i \right) \cdot \co\lambda_0 = s^\nu_{\mathcal I} \cdot (-\res_0 \overline{\nabla}).\vspace{8pt}
\end{equation*}

\subsubsection{Residue at infinity}

The following discussion applies whether the reproduction procedure is regular or singular at the origin. In order to determine the residue of the new Miura $\g$-oper $g \nabla g^{-1}$ at infinity, recall that the gauge transformation parameter can also be written in the form of a product of exponentials of simple root generators as $g = \prod_{i \in \mathcal I/2} e^{k_{i, 3} E_i} e^{k_{i, 2} E_{\bar\imath}} e^{k_{i, 1} E_i}$, where $k_{i,1}, k_{i,2}, k_{i,3} \in \M$ can be expressed in terms of the functions $f_{i,1}, f_{i,2}, f_{i,3} \in \M$ constructed above. In particular, if we started with functions $f_{k,1}$, $f_{k,2}$ and $f_{k,3}$ not all zero and verifying the cyclotomy condition \eqref{f123 cycl relations}, then one checks that the functions $k_{i,1}$, $k_{i,2}$ and $k_{i,3}$ for $i\in\mathcal{I}/2$ are all non-zero.

Let $\co\mu = \res_\infty \overline{\nabla}$ and fix $i\in\mathcal{I}/2$. Suppose that $\kappa=\langle \alpha_i, \co\mu+\co\rho \rangle = \langle \alpha_{\bar\imath}, \co\mu+\co\rho \rangle$ is non-negative and consider the three successive reproductions with overall parameter $g_i=e^{k_{i, 3} E_i} e^{k_{i, 2} E_{\bar\imath}} e^{k_{i, 1} E_i}$. According to the last part of Lemma \ref{TransfSimple}, since $k_{i,1}$ is non-zero, it has a simple pole at infinity and the first reproduction then yields a new $\g$-connection with residue $s_i\cdot\co\mu$ at infinity. Recalling that $\langle \alpha_i, \co\alpha_{\bar\imath} \rangle = \langle \alpha_{\bar\imath}, \co\alpha_i \rangle = -1$, one then finds $\langle \alpha_{\bar\imath}, s_i\cdot\co\mu+\co\rho \rangle = 2\kappa$. Thus, the last part of Lemma \ref{TransfSimple} still applies to the second reproduction and the new $\g$-connection has residue $s_{\bar\imath}s_i\cdot\co\mu$ at infinity. Finally, one has $\langle \alpha_i, s_{\bar\imath}s_i\cdot\co\mu+\co\rho \rangle = \kappa$. Applying Lemma \ref{TransfSimple} one last time, one finds that after the gauge transformation by $g_i$, the residue at infinity of the new $\g$-connection $\overline{g_i\nabla g_i^{-1}}$ is $s_is_{\bar\imath}s_i\cdot\left(\res_\infty \overline{\nabla} \right)$.

If $\langle \alpha_k, \co\mu+\co\rho \rangle$ is non-negative for some $k\in\mathcal{I}/2$, then $\langle \alpha_i, \co\mu+\co\rho \rangle$ is non-negative for all $i\in\mathcal{I}/2$. As $\langle \alpha_i, \co\alpha_j \rangle=0$ for any $j\neq i,\bar\imath$, the group elements $g_i$, $i\in\mathcal{I}/2$ commute with one another and
\begin{equation}
\langle \alpha_i, \res_\infty \overline{g_j\nabla g_j^{-1}}+\co\rho \rangle = \langle \alpha_i, \res_\infty \overline{\nabla}+\co\rho \rangle.
\end{equation}
Thus, we can apply successively (and in any order) the reasoning of the previous paragraph to the gauge transformations with parameters $g_i$, for all $i\in\mathcal{I}/2$. We then deduce that
\begin{equation*}
\res_\infty \overline{g \nabla g^{-1}} = \left( \prod_{i \in \mathcal I / 2} s_i s_{\bar\imath} s_i \right) \cdot (\res_\infty \overline{\nabla}) = s^\nu_{\mathcal I} \cdot (\res_\infty \overline{\nabla}),
\end{equation*}
where we used Lemma \ref{lem: Weyl nu}. We have thus proved the following.

\begin{theorem}\label{NonCommutingOrbit}
Let $\mathcal I \in I/\nu$ be an orbit of type $A_2^{\times |\mathcal I|/2}$.
Fix a $k \in \mathcal I$ and let $f_{k,1}, f_{k,2}, f_{k,3} \in \M$ be solutions of the coupled differential equations \eqref{RiccatiNonCommuting}, not all zero. Define $f_{i,1}, f_{i,2}, f_{i,3} \in \M$ for all $i \in (\mathcal I/2) \setminus \{ k \}$ recursively by
\begin{equation*}
f_{i,1}(t) \coloneqq \omega^{-1}f_{\nu^{-1}(i),1}(\omega^{-1}t), \qquad
f_{i,2}(t) \coloneqq \omega^{-1}f_{\nu^{-1}(i),2}(\omega^{-1}t), \qquad
f_{i,3}(t) \coloneqq \omega^{-2}f_{\nu^{-1}(i),3}(\omega^{-1}t).
\end{equation*}
and let $g := \prod_{i \in \mathcal I / 2} e^{f_{i,1} (E_i + E_{\bar\imath}) + f_{i,2} (E_i - E_{\bar\imath}) + f_{i,3} [E_i, E_{\bar\imath}]}$.

If we consider functions $f_{k,1}, f_{k,2}, f_{k,3}$ regular at the origin, then they can be chosen so that $g$ belongs to $N^{\hat\vs}(\M)$ if and only if
\begin{equation*}
\langle \alpha_k + \alpha_{\bar k}, \co\lambda_0 + \co\rho \rangle \equiv 0 \;\;\textup{mod}\;\; \frac{T}{|\mathcal I|}.
\end{equation*}
Moreover, in this case we have $g \nabla g^{-1} \in \MOp^\Gamma_\g(\CP)$ with $\res_0 \overline{g \nabla g^{-1}} = \res_0 \overline{\nabla}$.

Suppose now that at least one of the functions $f_{k,1}, f_{k,2}, f_{k,3}$ is singular at the origin, so that they are described by case $(ii)$ of Lemma \ref{LemmaRegNonCom}. Then $g$ is in $N^{\hat{\vs}}(\M)$ if and only if the singular behaviour of $f_{k,1}, f_{k,2}, f_{k,3}$ is as in the subcase $(a)$, without any conditions on the coweight $\co\lambda_0$. Moreover, in this case we have $g \nabla g^{-1} \in \MOp^\Gamma_\g(\CP)$ with $-\res_0 \overline{g \nabla g^{-1}} = s^\nu_{\mathcal{I}}\cdot(-\res_0 \overline{\nabla})$.

Finally, if $\langle \alpha_k, \res_\infty \overline{\nabla} + \co\rho \rangle$ is non-negative, then we have $\res_\infty \overline{g \nabla g^{-1}} = s^\nu_{\mathcal I} \cdot (\res_\infty \overline{\nabla})$ (whether we started with functions regular or singular at the origin).
\qed
\end{theorem}

\begin{rem}
Note that the approach of this section may equally be applied to the case of an orbit of type $A_1^{\times |\mathcal I|}$ for even $|\mathcal I|$: in this case, there are no functions $f_{i,3}$ and we only get the first two conditions $\omega^{(1+\eta) |\mathcal I| / 2} = \pm 1$ above, which are indeed equivalent to the condition $\langle \alpha_k, \co\lambda_0 + \co\rho \rangle \equiv 0 \;\textup{mod}\; \frac{T}{|\mathcal I|}$ as given in Theorem \ref{CommutingOrbit}.
\end{rem}

\begin{exmp} \label{ReproSl3}
Let us illustrate Theorem \ref{NonCommutingOrbit} with an example. We consider the Lie algebra $\g=\sl_3$ and the cyclotomic Miura $\sl_3$-oper $\nabla$ defined in equation \eqref{ExpleSl3} of example \ref{Sl3}.

We consider the orbit of simple roots $\lbrace \alpha_1, \alpha_2 \rbrace$, which is of type $A_2$, and the unipotent element
\begin{equation}\label{gSl3}
g = e^{\tilde{f}_1 E_1 + \tilde{f}_2 E_2 + \tilde{f}_3 E_{12}},
\end{equation}
expressed in the basis $E_1, E_2, E_{12} \coloneqq [E_1,E_2]$. In Theorem \ref{NonCommutingOrbit} and the discussion above, $g$ is given by $g = e^{f_1 (E_1+E_2)+f_2 (E_1-E_2)+ f_3 E_{12}}$, in a basis of eigenvectors of $\vs$. The two expressions are simply related by $\tilde{f}_1=f_1+f_2$, $\tilde{f}_2=f_1-f_2$ and $\tilde{f}_3=f_3$. In the fundamental representation, we have
\begin{equation*}
g=\begin{pmatrix}
1 & \tilde{f}_1 & \tilde{f}_3 + \frac{1}{2}\tilde{f}_1\tilde{f}_2 \\ 0 & 1 & \tilde{f}_2 \\ 0 & 0 & 1
\end{pmatrix}.
\end{equation*}

The gauge transformation $g\nabla g^{-1}$ is a Miura $\sl_3$-oper if and only if the functions $\tilde{f}_i$ are solutions of the following system of coupled differential equations
\begin{eqnarray*}
&\tilde{f}'_1(t) = \dfrac{\eta}{t}\tilde{f}_1(t) - \tilde{f}_1(t)^2  - \frac{1}{2}\tilde{f}_1(t)\tilde{f}_2(t) - \tilde{f}_3(t),\\
&\tilde{f}'_2(t) = \dfrac{\eta}{t}\tilde{f}_2(t) - \tilde{f}_2(t)^2 - \frac{1}{2}\tilde{f}_1(t)\tilde{f}_2(t) + \tilde{f}_3(t), \\
&\tilde{f}'_3(t) = 2\dfrac{\eta}{t}\tilde{f}_3(t) - \Bigl(\tilde{f}_1(t)+\tilde{f}_2(t)\Bigr) \tilde{f}_3(t)   + \frac{1}{2} \Bigl(\tilde{f}_1(t)^2+\tilde{f}'_1(t)\Bigr)  \tilde{f}_2(t) - \frac{1}{2} \Bigl(\tilde{f}_2(t)^2+\tilde{f}'_2(t)\Bigr)  \tilde{f}_1(t).
\end{eqnarray*}
The regularized functions $\tilde{h}_1(t)=t^{-\eta}\tilde{f}_1(t)$, $\tilde{h}_2(t)=t^{-\eta}\tilde{f}_2(t)$ and $\tilde{h}_3(t)=t^{-2\eta}\tilde{f}_3(t)$ are solutions of differential equations regular at the origin. The solution of the above differential system with initial conditions $\tilde{h}_1(0)=a$, $\tilde{h}_2(0)=b$ and $\tilde{h}_3(0)=c$ is
\begin{subequations}\label{fSl3}
\begin{eqnarray}
&\tilde{f}_1(t) = 2\mu \, t^{\mu-1} \dfrac{(ab+2c)\,t^\mu+2\mu a}{(ab+2c)\,t^{2\mu}+4\mu a \,t^\mu + 4\mu^2}, \\
&\tilde{f}_2(t) = 2\mu \, t^{\mu-1} \dfrac{(ab-2c)\,t^\mu+2\mu b}{(ab-2c)\,t^{2\mu}+4\mu b \,t^\mu + 4\mu^2}, \\
&\tilde{f}_3(t) = 4\mu^3 \, t^{2\mu-2} \dfrac{\bigl( (ab+2c)b-(ab-2c)a \bigr) \, t^\mu + 4 \mu c}{\bigl((ab+2c)\,t^{2\mu}+4\mu a \,t^\mu + 4\mu^2\bigr)\bigl((ab-2c)\,t^{2\mu}+4\mu b \,t^\mu + 4\mu^2\bigr)}.
\end{eqnarray}
\end{subequations}
where we defined $\mu=\eta+1$.

To simplify the analysis of this solution, we define functions $A$ and $B$ such that $\tilde{f}_1(t)=A(t;a,b,c)$, $\tilde{f}_2(t)=A(t;b,a,-c)$ and $\tilde{f}_3(t)=B(t;a,b,c)=-B(t;b,a,-c)$. One checks easily that
\begin{align*}
\omega^{-1} A(\omega^{-1}t;a,b,c) &= A(t;\omega^{-\mu}a,\omega^{-\mu}b,\omega^{-2\mu}c), \\
\omega^{-2} B(\omega^{-1}t;a,b,c) &= B(t;\omega^{-\mu}a,\omega^{-\mu}b,\omega^{-2\mu}c).
\end{align*}
Noting that $g$ is $\hat{\vs}$-invariant if and only if we have $\omega^{-1}\tilde{f}_1(\omega^{-1}t)= \tilde f_2(t)$, $\omega^{-1}\tilde{f}_2(\omega^{-1}t)= \tilde f_1(t)$ and $\omega^{-2}\tilde{f}_3(\omega^{-1}t)=- \tilde f_3(t)$, the condition for the reproduction to be cyclotomic can be re-expressed as
\begin{align*}
A(t;\omega^{-\mu}a,\omega^{-\mu}b,\omega^{-2\mu}c) &= A(t;b,a,-c), \\
A(t;\omega^{-\mu}b,\omega^{-\mu}a,-\omega^{-2\mu}c) &= A(t;a,b,c), \\
B(t;\omega^{-\mu}a,\omega^{-\mu}b,\omega^{-2\mu}c) &= B(t;b,a,-c).
\end{align*}
These conditions are equivalent to $a=\omega^\mu b$, $b=\omega^\mu a$ and $c=-\omega^{2\mu}c$. As in the general discussion above, we distinguish three possible cases where these conditions have non-trivial solutions:
\begin{itemize}
\item $\omega^\mu = 1$, $a=b \in \C^\times$ and $c=0$, \emph{i.e.} $h_2(0)=h_3(0)=0$, in which case
\begin{equation*}
f_1(t) = \frac{2\mu a \, t^{\mu-1}}{a\, t^\mu+2\mu}, \; \; \; \; f_2(t)=0, \; \; \; \; f_3(t)=0,
\end{equation*}
\item $\omega^\mu = -1$, $a=-b \in \C^\times$ and $c=0$, \emph{i.e.} $h_1(0)=h_3(0)=0$, in which case
\begin{equation*}
f_1(t) = \frac{2\mu a^2 \, t^{2\mu-1}(a^2\,t^{2\mu}-12\mu^2)}{R(t)}, \; \; \; \; f_2(t)=\frac{4\mu^2 a \, t^{\mu-1}(a^2\,t^{2\mu}+4\mu^2)}{R(t)}, \; \; \; \; f_3(t)=\frac{8\mu^3 a^3 \, t^{3\mu-2}}{R(t)},
\end{equation*}
where $R(t) = a^4\, t^{4\mu}-24\mu^2a^2 \, t^{2\mu} + 16\mu^4$,
\item $\omega^{2\mu} = -1$, $a=b=0$ and $c\in\C^\times$, \emph{i.e.} $h_1(0)=h_2(0)=0$, in which case
\begin{equation}\label{Cas3}
f_1(t) = \frac{2\mu c^2\, t^{4\mu-1}}{c^2\, t^{4\mu}-4\mu^4}, \; \; \; \; f_2(t)=-\frac{4\mu^3 c\, t^{2\mu-1}}{c^2\, t^{4\mu}-4\mu^4}, \; \; \; \; f_3(t)=-\frac{4\mu^4 c\, t^{2\mu-2}}{c^2\, t^{4\mu}-4\mu^4},
\end{equation}
\end{itemize}
where the regularised functions $h_1=\frac{1}{2}(\tilde{h}_1+\tilde{h}_2)$, $h_2=\frac{1}{2}(\tilde{h}_1-\tilde{h}_2)$ and $h_3=\tilde{h}_3$ correspond to the basis of eigenvectors $\lbrace E_1+E_2, E_1-E_2, [E_1,E_2] \rbrace$. In each of these cases, the solution depends on an arbitrary parameter in $\C^\times$.

To end this example we discuss also the solutions of the differential system singular at the origin. Starting from any of the above three types of regular solutions, a singular solution can be obtained by considering the limit where the free parameter goes to infinity. One can check that each of the three cases yields the same singular solution
\begin{equation}\label{SingSl3}
g(t)=\exp\left(\frac{2(\eta+1)}{t}\left(E_1+E_2\right) \right).
\end{equation}
Note that this solution is always cyclotomic, without imposing any restriction on $\omega$ or $\eta$.
\end{exmp}

\section{Cyclotomic Miura $\g$-opers and flag varieties}
\label{ReproFlag}

Recall the setup of \S\ref{Repro}. Throughout this section we will also consider a cyclotomic Miura $\g$-oper of the form $\nabla=d+(p_{-1}+u)dt$ with $u \, dt \in \Omega^{\hat{\nu}}(\h)$ as in \eqref{cycl Miura lambda0}. In particular, we still assume for the time being that $-\res_0 \overline{\nabla}=\co\lambda_0$ is a $\nu$-invariant integral dominant coweight. We will further assume here that $\nabla$ is monodromy-free. In \S\ref{sec: non-integral}, however, we will show how to weaken the integrality assumption on the coweight $\co\lambda_0$ together with the related assumption on the trivial monodromy at the origin.

Consider the set
\begin{equation} \label{cycl MOp}
\MOp^\Gamma_\g(\CP)_{[\nabla]_\Gamma} \coloneqq \big\{ \widetilde{\nabla} \in \MOp^\Gamma_\g(\CP) \,\big|\, [\widetilde{\nabla}]_\Gamma = [\nabla]_\Gamma \big\}
\end{equation}
of cyclotomic Miura $\g$-opers with underlying cyclotomic $\g$-oper $[\nabla]_\Gamma$. By definition, any cyclotomic Miura $\g$-oper in $\MOp^\Gamma_\g(\CP)_{[\nabla]_\Gamma}$
is a representative of $[\nabla]_\Gamma$ and hence is of the form $\nabla^g \coloneqq g \nabla g^{-1}$ for some $g \in N^{\hat\vs}(\M)$. The goal of this section is to describe the space $\MOp^\Gamma_\g(\CP)_{[\nabla]_\Gamma}$.

A finite analog of \eqref{cycl MOp} is the set $\MOp_{\g^\nu, [\co\lambda_0]_{W^\nu}}^{\rm fin}$, introduced in \S\ref{finite opers}, consisting of all finite Miura $\g^\nu$-opers which are connected to the given finite Miura $\g^\nu$-oper $p_{-1} - \co\lambda_0 - \co\rho$ by the adjoint action of $N^\nu$. We saw in Theorem \ref{MOp fin iso W} that $\MOp_{\g^\nu, [\co\lambda_0]_{W^\nu}}^{\rm fin}$ is in bijection with the $\nu$-invariant subgroup $W^\nu$ of the Weyl group $W$. The affine counterpart of this statement is Theorem \ref{MiuraFlagCyc} which provides an isomorphism between $\MOp^\Gamma_\g(\CP)_{[\nabla]_\Gamma}$ and the subset $(G/B_-)^\vartheta$ of the flag variety $G/B_-$ consisting of points fixed by some automorphism $\vartheta$ which we introduce in \S\ref{ReproGen}. In \S\ref{sec: cell decomp} we show that $(G/B_-)^\vartheta$ admits a cell decomposition where the cells are labelled by elements of $W^\nu$.

\subsection{Generic cyclotomic Miura $\g$-opers}
\label{ReproGen}

If $\co\lambda_0 = 0$ then the origin is a regular point of the given cyclotomic Miura $\g$-oper $\nabla$. In this case we let $Y \in B_-(\M)$ be the solution of $\nabla Y = 0$ with $Y(0) = \Id$, which exists by Lemma \ref{BorelSol}.
Then for any $\nabla^g \in \MOp^\Gamma_\g(\CP)_{[\nabla]_\Gamma}$ such that $g \in N^{\hat\vs}(\M)$ is also regular at $0$, we will denote by $Y_g \in B_-(\M)$ the solution of $\nabla^g Y_g = 0$ with $Y_g(0)=\Id$.
If instead $\co\lambda_0 \neq 0$ then $\nabla$ is singular at $0$ so we cannot consider such a solution. In this case, we will work with a suitable regularisation of $\nabla$ defined as follows. Let
\begin{equation} \label{nabla reg}
\nabla_{\rm r} \coloneqq t^{-\co\lambda_0} \nabla t^{\co\lambda_0} = d + \sum_{k \in I} t^{\langle \alpha_k, \co\lambda_0 \rangle} F_k \, dt - \sum_{r=0}^{T-1} \left( \sum_{i=1}^N \frac{\nu^r (w_i\cdot\co \lambda_i)}{t - \omega^r z_i} + \sum_{j=1}^m \frac{\nu^r y_j \cdot 0}{t - \omega^r x_j} \right) dt.
\end{equation}
The assumption that the coweight $\co\lambda_0$ is integral dominant ensures this $\g$-connection is regular at the origin. We shall refer to $\nabla_{\rm r}$ as the \emph{regularisation} of $\nabla$ at the origin. By Lemma \ref{BorelSol}, we can then consider the solution $Y \in B_-(\M)$ of $\nabla_{\rm r} Y = 0$ such that $Y(0) = \Id$.
For any $\nabla^g \in \MOp^\Gamma_\g(\CP)_{[\nabla]_\Gamma}$ we consider its regularisation $(\nabla^g)_{\rm r} = (g \nabla g^{-1})_{\rm r} = g_{\rm r} \nabla_{\rm r} g_{\rm r}^{-1}$ where $g_{\rm r} \coloneqq t^{-\co\lambda_0} g t^{\co\lambda_0} \in N(\M)$. Provided the latter is regular at the origin, we can define the solution $Y_g$ of $(\nabla^g)_{\rm r} Y_g=0$ with $Y_g(0) = \Id$.

An important feature of the regularisation procedure is that all the regularised objects such as the $\g$-connection \eqref{nabla reg} are $\Gamma$-equivariant not with respect to the automorphism $\vs \in \Aut G$ but rather with respect to the automorphism
\begin{equation} \label{vartheta def}
\vartheta \coloneqq \Ad_{\omega^{-\co\lambda_0}} \circ \, \vs \in \Aut G.
\end{equation}
The following proposition describes the $\Gamma$-equivariance properties of the regularised $\g$-connection $\nabla_{\rm r}$, of the regularised gauge transformation parameter $g_{\rm r} \in N(\M)$ and of a solution $Y \in G(\M)$ of the equation $\nabla_{\rm r} Y = 0$.

\begin{proposition}\label{InvY}
Define $D_{\rm r} : G(\M) \to \Omega(\g)$, $Z \mapsto (\nabla_{\rm r} Z) Z^{-1}$ and ${\rm reg} : N(\M) \to N(\M)$, $g \mapsto g_{\rm r}$. Both are $\Gamma$-equivariant 
in the sense that the following diagrams are commutative
\begin{equation*}
\begin{tikzpicture}    
\matrix (m) [matrix of math nodes, row sep=3em, column sep=3em,text height=1.5ex, text depth=0.25ex]    
{
G(\M) & \Omega(\g)\\
G(\M) & \Omega(\g)\\
};
\path[->] (m-1-1) edge node[above]{$D_{\rm r}$} (m-1-2);
\path[->] (m-2-1) edge node[below]{$D_{\rm r}$} (m-2-2);
\path[->] (m-1-1) edge node[left]{$\hat\vartheta$} (m-2-1);
\path[->] (m-1-2) edge node[right]{$\hat\vartheta$} (m-2-2);
\end{tikzpicture}
\qquad
\begin{tikzpicture}    
\matrix (m) [matrix of math nodes, row sep=3em, column sep=3em,text height=1.5ex, text depth=0.25ex]    
{
N(\M) & N(\M)\\
N(\M) & N(\M)\\
};
\path[->] (m-1-1) edge node[above]{$\rm reg$} (m-1-2);
\path[->] (m-2-1) edge node[below]{$\rm reg$} (m-2-2);
\path[->] (m-1-1) edge node[left]{$\hat\vs$} (m-2-1);
\path[->] (m-1-2) edge node[right]{$\hat\vartheta$} (m-2-2);
\end{tikzpicture}
\end{equation*}
Let $Y \in G(\M)$ be any solution of $\nabla_{\rm r} Y = 0$ regular at $0$. Then $Y \in G^{\hat\vartheta}(\M)$ if and only if $Y(0) \in G^\vartheta$. In particular, $Y \in G^{\hat\vartheta}(\M)$ if $Y(0) = \Id$.
\begin{proof}
The $\Gamma$-equivariance of the map $D_{\rm r}$ is the statement that $\hat\vartheta \big( (\nabla_{\rm r} Z) Z^{-1} \big) = (\nabla_{\rm r} \hat\vartheta Z)(\hat\vartheta Z)^{-1}$ for every $Z \in G(\M)$. This follows from expression \eqref{nabla reg} for $\nabla_{\rm r}$ noting that $\hat\vartheta\big( dZ Z^{-1} \big) = (d \hat\vartheta Z)(\hat\vartheta Z)^{-1}$, the action of $\hat\vartheta$ on $\Omega^{\hat\vs}(\h)$ coincides with that of $\hat\vs$ and
\begin{equation*}
\hat\vartheta \Big( t^{- \co\lambda_0} p_{-1} t^{\co\lambda_0} dt \Big) = \Ad_{\omega^{-\co\lambda_0}} \big( (\omega^{-1} t)^{- \co\lambda_0} \vs p_{-1} (\omega^{-1} t)^{\co\lambda_0} d(\omega^{-1} t) \big) = t^{- \co\lambda_0} p_{-1} t^{\co\lambda_0} dt
\end{equation*}
using the $\vs$-invariance of $\co\lambda_0$ and Proposition \ref{ExpoSigma}.
Similarly, the $\Gamma$-equivariance of the map $g \mapsto g_{\rm r}$ is seen as follows
\begin{equation*}
\hat\vartheta g_{\rm r} = \Ad_{\omega^{- \co\lambda_0}} \circ \, \hat\vs\big( t^{- \co\lambda_0} g t^{\co\lambda_0} \big) = \Ad_{\omega^{- \co\lambda_0}} \big( (\omega^{-1} t)^{- \co\lambda_0} \hat\vs g (\omega^{-1} t)^{\co\lambda_0} \big) = (\hat\vs g)_{\rm r}.
\end{equation*}

For the last statement, suppose first that $Y \in G^{\hat\vartheta}(\M)$ is a solution to $\nabla_{\rm r} Y = 0$ which is regular at $0$. This implies $\hat\vartheta Y = Y$ or in other words $\vartheta \circ Y \circ \mu_{\omega^{-1}}^\ast = Y$. Evaluating this at the origin we obtain $\vartheta Y(0) = Y(0)$. Conversely, let $Y \in G(\M)$ be a solution of $\nabla_{\rm r} Y = 0$ such that $\vartheta Y(0) = Y(0)$. Using the first part of the proposition, or more precisely the $\Gamma$-equivariance of $D_{\rm r}$, we deduce that $\hat\vartheta Y$ also satisfies the same equation since $(\nabla_{\rm r} \hat\vartheta Y)(\hat\vartheta Y)^{-1} = \hat\vartheta \big( (\nabla_{\rm r} Y) Y^{-1} \big) = 0$. However, by assumption the two solutions $Y$ and $\hat\vartheta Y$ coincide at the origin and therefore are equal.
\end{proof}
\end{proposition}

Introduce the following subset of cyclotomic Miura $\g$-opers, which we call \emph{generic} at the origin,
\begin{equation} \label{generic Miura opers}
\MOp^\Gamma_\g(\CP)^{\rm gen}_{\nabla, 0} \coloneqq \{ \nabla^g \in \MOp^\Gamma_\g(\CP)_{[\nabla]_\Gamma} \,|\, g_{\rm r} \;\text{is regular at}\; 0 \}.
\end{equation}
The following theorem describes $\MOp^\Gamma_\g(\CP)^{\rm gen}_{[\nabla]_\Gamma, 0}$ as an $N^\vartheta$-orbit in $\MOp^\Gamma_\g(\CP)_{[\nabla]_\Gamma}$ through $\nabla$.

\begin{theorem}\label{CycReproReg}
Let $\nabla \in \MOp^\Gamma_\g(\CP)$. If $\nabla^g \in \MOp^\Gamma_\g(\CP)^{\rm gen}_{\nabla, 0}$ then $Y g_{\rm r}(0)^{-1} = g_{\rm r}^{-1} Y_g$ and $g_{\rm r}(0) \in N^\vartheta$.
Conversely, for every $g_0 \in N^{\vartheta}$ there exists unique $n \in N^{\hat\vartheta}(\M)$ and $\widetilde{Y} \in B^{\hat\vartheta}_-(\M)$ such that
\begin{equation*}
Y g^{-1}_0 = n^{-1} \widetilde{Y}.
\end{equation*}
Defining $g \coloneqq t^{\co\lambda_0} n t^{- \co\lambda_0}$ so that $n = g_{\rm r}$, we have $\nabla^g \in \MOp^\Gamma_\g(\CP)^{\rm gen}_{\nabla, 0}$, $g_{\rm r}(0) = g_0$ and $\widetilde{Y} = Y_g$.
\begin{proof}
Suppose that $\nabla^g \in \MOp^\Gamma_\g(\CP)^{\rm gen}_{\nabla, 0}$. By definition $g_{\rm r}$ is regular at the origin so $g_{\rm r}(0) \in N$ and clearly $\nabla_{\rm r} Y g_{\rm r}(0)^{-1} = 0$. On the other hand, we have
$\nabla_{\rm r} g_{\rm r}^{-1} Y_g = g_{\rm r}^{-1} (\nabla^g)_{\rm r} Y_g = 0$.
It therefore follows by the uniqueness of the solution to $\nabla_{\rm r} Z = 0$ with $Z(0) = g_{\rm r}(0)^{-1}$ that $Y g_{\rm r}(0)^{-1} = g_{\rm r}^{-1} Y_g$. It remains to show that $g_{\rm r}(0) \in N^\vartheta$. However, since $\nabla^g \in \MOp^\Gamma_\g(\CP)_{[\nabla]_\Gamma}$ we have $g \in N^{\hat\vs}(\M)$, \emph{i.e.} $\hat\vs g = g$, and therefore by Proposition \ref{InvY} we deduce that $\hat\vartheta g_{\rm r} = g_{\rm r}$. In particular, evaluating the latter at the origin gives $g_{\rm r}(0) \in N^\vartheta$.

Conversely, let $g_0 \in N$ and consider the element $Y g_0^{-1} \in G(\M)$. Let $\varphi \in \sO(G)$ be the coordinate function on $G$ whose nonzero set is the open cell $N B_- = \{ x \in G | \varphi(x) \neq 0 \}$ in $G$. Let $f \in \M$ denote the image of $\varphi \in \sO(G)$ under the homomorphism $Y g_0^{-1} : \sO(G) \to \M$.Note that $f$ is not identically zero since $f(0) =\varphi(Y(0) g_0^{-1}) \neq 0$, which follows from $Y(0) g_0^{-1} = g_0^{-1} \in N B_-$.
Since the number of zeroes of a non-zero meromorphic function on $\CP$ is finite, we have $\varphi(Y(t) g_0^{-1}) = f(t) \neq 0$ for all but finitely many $t \in \CP$. In other words, $Y(t) g_0^{-1} \in N B_-$ for all but finitely many $t \in \CP$. We may therefore write $Y g^{-1}_0 = n^{-1} \widetilde Y$, for some $n \in N(\M)$ and $\widetilde Y \in B_-(\M)$.

By definition of $g \in N(\M)$ we have $n = t^{-\co\lambda_0} g t^{\co\lambda_0} = g_{\rm r}$
so that $(\nabla^g)_{\rm r} \widetilde Y = g_{\rm r} \nabla_{\rm r} g_{\rm r}^{-1} \widetilde Y = n \nabla_{\rm r} n^{-1} \widetilde Y = n \nabla_{\rm r} Y g_0^{-1} = 0$. Since $\widetilde Y \in B_-(\M)$ it follows by the converse in Lemma \ref{BorelSol} that $(\nabla^g)_{\rm r} \in \Conn_{\b_-}(\CP)$ and hence $\nabla^g = t^{\co\lambda_0} (\nabla^g)_{\rm r} t^{-\co\lambda_0} \in \Conn_{\b_-}(\CP)$. But since $\nabla^g \in \op_\g(\CP)$ we conclude that $\nabla^g \in \MOp_\g(\CP)$ and hence $\nabla^g \in \MOp_\g(\CP)_{[\nabla]}$. Moreover, letting $t = 0$ we obtain $n(0)^{-1} \widetilde{Y}(0) = g_0^{-1}$, or equivalently $\widetilde{Y}(0) = g_{\rm r}(0) g_0^{-1}$. But $N \cap B_- = \{ \Id \}$ from which we deduce $g_{\rm r}(0) = g_0$ and $\widetilde{Y}(0) = \Id$. The first equality shows $g_{\rm r}$ is regular at $0$. The second equality implies that $\widetilde{Y} = Y_g$ by the uniqueness of the solution to $(\nabla^g)_{\rm r} Z = 0$ with $Z(0) = \Id$.
It just remains to show that $\hat\vartheta n = n$ and $\hat\vartheta \widetilde{Y} = \widetilde{Y}$. But from Proposition \ref{InvY} we know that $\hat\vartheta Y = Y$. So since $g_0 \in N^\vartheta$ it follows that $\hat\vartheta (Y g_0^{-1}) = Y g_0^{-1}$. Hence $n^{-1} \widetilde{Y} = \hat\vartheta(n^{-1} \widetilde{Y}) = (\hat\vartheta n)^{-1} \hat\vartheta \widetilde{Y}$, or equivalently $(\hat\vartheta n) n^{-1} = (\hat\vartheta \widetilde{Y}) \widetilde{Y}^{-1}$. Now $\vartheta$ stabilises $N$ and $B_-$ so that $\hat\vartheta n \in N(\M)$ and $\hat\vartheta \widetilde{Y} \in B_-(\M)$. Finally, since $N \cap B_- = \{ \Id \}$ we deduce that $\hat\vartheta n = n$ and $\hat\vartheta \widetilde{Y} = \widetilde{Y}$, as required.
\end{proof}
\end{theorem}

\begin{rem} \label{ReproRem}
In Theorems \ref{CommutingOrbit} and \ref{NonCommutingOrbit}, we found necessary and sufficient conditions on $\co\lambda_0$ for the existence of a cyclotomic gauge transformation along the $\Gamma$-orbit $\{ \alpha_i \}_{i \in \mathcal I}$ of some simple root $\alpha_k$, $k \in I$. These results can now be seen as a particular case of Theorem \ref{CycReproReg}: they correspond to the conditions on $\co\lambda_0$, and thus on $\vartheta$, for the existence of a $\vartheta$-invariant element in the subgroup of $N$ generated by the $G_{\alpha_i}$, $i \in \mathcal I$, cf. \S\ref{AdGroup}.

Consider the concrete example of a simple root $\alpha$ with an orbit $\{ \alpha, \beta \}$ of type $A_1 \times A_1$. Theorem \ref{CommutingOrbit} asserts that there exists a cyclotomic reproduction in the direction of this orbit if and only if $\langle \alpha, \co\lambda_0 + \co\rho \rangle \equiv 0 \; \text{mod} \; \frac{T}{2}$. In the above language, there exists such a cyclotomic reproduction if and only if a linear combination of $E_\alpha$ and $E_\beta$ lives in $\n^\vartheta$. Yet, we have $\vartheta(E_\alpha \pm E_\beta) = \pm \omega^{-\langle \alpha, \co\lambda_0 + \co\rho \rangle} (E_\alpha \pm E_\beta)$, so this condition is indeed equivalent to $\omega^{-\langle \alpha, \co\lambda_0 + \co\rho \rangle} = \pm 1$, and thus to $\langle \alpha, \co\lambda_0 + \co\rho \rangle \equiv 0 \; \text{mod} \; \frac{T}{2}$. See also example \ref{ReproGenSl3} below.
\end{rem}

\begin{exmp}\label{ReproGenSl3}
Recall the cyclotomic Miura $\sl_3$-oper $\nabla$ of example \ref{Sl3} given by \eqref{ExpleSl3}, whose cyclotomic reproduction we studied in example \ref{ReproSl3}. We now return to this example in light of Theorem \ref{CycReproReg}.

The regularised connection associated with $\nabla$ is simply
\begin{equation}\label{Sl3Reg}
\nabla_{\rm r} = t^{-\co\lambda_0} \nabla t^{\co\lambda_0} = d + \begin{pmatrix}
0 & 0 & 0 \\ t^\eta & 0 & 0 \\ 0 & t^\eta & 0
\end{pmatrix}dt, \;\;\;\;\; \text{with } \; \; t^{-\co\lambda_0} = \begin{pmatrix}
t^{-\eta} & 0 & 0 \\ 0 & 1 & 0 \\ 0 & 0 & t^\eta
\end{pmatrix}.
\end{equation}
The solution of equation $\nabla_{\rm r} Y = 0$ with initial condition $Y(0)=\Id$ is (recall that $\mu=\eta+1$)
\begin{equation}\label{Sl3SolReg}
Y(t) = \begin{pmatrix}
1 & 0 & 0 \\ -\dfrac{t^\mu}{\mu} & 1 & 0 \\ \dfrac{t^{2\mu}}{2\mu^2} & -\dfrac{t^\mu}{\mu} & 1
\end{pmatrix}.
\end{equation}
Let us now consider an arbitrary element $g_0$ of $N$ which we parametrise as
\begin{equation*}
g_0 = e^{aE_1 + bE_2 + cE_{12}} = \begin{pmatrix}
1 & a & c+\frac{1}{2}ab \\ 0 & 1 & b \\ 0 & 0 & 1
\end{pmatrix}.
\end{equation*}
Following Theorem \ref{CycReproReg}, we perform the $NB_-$ factorisation of $Y(t)g_0^{-1}$ and get $Y(t)g_0^{-1}=n(t)^{-1}\widetilde{Y}(t)$ with $n \in N(\M)$ and $\widetilde{Y} \in B_-(\M)$. In particular, we have
\begin{equation*}
n(t) = e^{\tilde{h}_1(t)E_1 + \tilde{h}_2(t)E_2 + \tilde{h}_3(t)E_{12}} = \begin{pmatrix}
1 & \tilde{h}_1(t) & \tilde{h}_3(t)+\frac{1}{2}\tilde{h}_1(t)\tilde{h}_2(t) \\ 0 & 1 & \tilde{h}_2(t) \\ 0 & 0 & 1
\end{pmatrix},
\end{equation*}
with
\begin{eqnarray*}
&\tilde{h}_1(t) = 2\mu \dfrac{(ab+2c)\,t^\mu+2\mu a}{(ab+2c)\,t^{2\mu}+4\mu a \,t^\mu + 4\mu^2}, \\
&\tilde{h}_2(t) = 2\mu \dfrac{(ab-2c)\,t^\mu+2\mu b}{(ab-2c)\,t^{2\mu}+4\mu b \,t^\mu + 4\mu^2}, \\
&\tilde{h}_3(t) = 4\mu^3 \dfrac{\bigl( (ab+2c)b-(ab-2c)a \bigr) \, t^\mu + 4 \mu c}{\bigl((ab+2c)\,t^{2\mu}+4\mu a \,t^\mu + 4\mu^2\bigr)\bigl((ab-2c)\,t^{2\mu}+4\mu b \,t^\mu + 4\mu^2\bigr)}.
\end{eqnarray*}
Defining $g=t^{\co\lambda_0} n t^{-\co\lambda_0}$ we recover the element \eqref{gSl3} with the functions $\tilde{f}_i$ given by \eqref{fSl3}. We have thus checked that the $NB_-$-factorisation of $Yg_0^{-1}$ allows one to find the most general reproduction (regular at the origin) of the Miura $\sl_3$-oper $\nabla$. Moreover, one also checks that $g_0=n(0)$ is indeed the initial condition of the regularised transformation $n= t^{-\co\lambda_0}g t^{\co\lambda_0}$.

We now discuss the condition under which this reproduction is cyclotomic, following Theorem \ref{CycReproReg}, and compare it to the results found in example \ref{ReproSl3}. A straightforward calculation gives
\begin{equation*}
\vartheta g_0 = \begin{pmatrix}
1 & \omega^{-\mu}b & \omega^{-2\mu}\left(-c+\frac{1}{2}ab\right) \\
0 & 1 & \omega^{-\mu}a \\ 0 & 0 & 1
\end{pmatrix}.
\end{equation*}
According to Theorem \ref{CycReproReg}, the reproduction is cyclotomic if the initial condition $g_0$ is $\vartheta$-invariant. Thus, we recover the conditions
\begin{equation*}
a=\omega^{-\mu}b, \;\;\;\;\; b=\omega^{-\mu}a \;\;\;\;\; \text{and} \;\;\;\;\; c=-\omega^{-2\mu}c,
\end{equation*}
which we found in example \ref{ReproSl3}. We verified on an explicit example the general statement of remark \ref{ReproRem}, namely that the condition $\langle \alpha_1+\alpha_2, \co\lambda_0+\co\rho \rangle \equiv 0 \; \text{mod} \; \frac{T}{2}$ given in Theorem \ref{NonCommutingOrbit} is simply the condition for the $\vartheta$-invariant subgroup $N^\vartheta$ to be non-trivial.
\end{exmp}

\subsection{Non-integral coweight $\co\lambda_0$} \label{sec: non-integral}

In \S\ref{Repro} and \S\ref{ReproGen}, we studied the reproduction of a cyclotomic Miura $\g$-oper with a dominant integral coweight $\co\lambda_0$ at the origin. By virtue of the $\Gamma$-equivariance, this procedure can be extended to any dominant coweight $\co\lambda_0$ such that for all $i \in I$,
\begin{equation} \label{non-integral}
\langle \alpha_i, \co\lambda_0 \rangle \in \mathbb{Q}.
\end{equation}

We start with a cyclotomic Miura $\g$-oper of the form
\begin{equation} \label{Miura non-int}
\nabla = d + p_{-1} dt - \frac{\co\lambda_0}{t} dt + \ms r,
\end{equation}
with $\co\lambda_0$ as above and $\ms r \in \h(\M)dt$ regular at the origin. One can always find a $q \in \Z_{\geq 1}$ such that $q \langle \alpha_i, \co\lambda_0 \rangle \in \Z$ for any $i \in I$, \emph{i.e.} such that $q \co\lambda_0$ is integral dominant. Let $\pi : \CP \to \CP$ be a $q$-sheeted cover branched over the origin and infinity, given by $u \mapsto t = u^q$ in terms of a global coordinate $u$ on $\C \subset \CP$.
Let $\widetilde{\M}$ denote the algebra of rational functions in $u$. We can regard $\M$ as a subalgebra of $\widetilde{\M}$ via the pullback $\pi^\ast : \M \hookrightarrow \widetilde{\M}$.

Let $\nabla = d + \ms A$ be a $\g$-connection meromorphic in $t$, \emph{i.e.} with $\ms A \in \g(\M) dt$, and consider the associated $\g$-connection $\widetilde{\nabla} \coloneqq d + \pi^\ast \ms A$ on the cover, noting that $\pi^\ast \ms A \in \g(\widetilde \M) du$. Define the primitive $(qT)^{\rm th}$-root of unity $\widetilde \omega \coloneqq \omega^{\frac{1}{q}}$ and let
\begin{equation*}
\widetilde \Gamma \coloneqq \langle \widetilde \omega \rangle \subset \C^\times
\end{equation*}
denote the corresponding copy of the cyclic group of order $qT$.
If the $\g$-connection $\nabla$ is cyclotomic then its $\Gamma$-equivariance property $\hat\vs(\ms A) = \ms A$ translates to the $\widetilde \Gamma$-equivariance property of $\widetilde \nabla$, namely $\widetilde\vs (\pi^\ast \ms A) = \pi^\ast \ms A$, where for any automorphism $\up \in \Aut \g$ and any $X \otimes \varpi \in \g(\widetilde\M) du$ we define
\begin{equation*}
\widetilde\up (X \otimes \varpi) \coloneqq \up(X) \otimes \mu^\ast_{\widetilde \omega^{-1}} \varpi.
\end{equation*}
Here we have used the fact that $\pi \circ \mu_{\widetilde \omega} = \mu_\omega \circ \pi$ and we note that $\vs^{qT} = \Id$.
The cyclotomic Miura $\g$-oper \eqref{Miura non-int} pulls back to the following $\g$-connection on the cover
\begin{equation*}
\widetilde \nabla = d + q u^{q-1} p_{-1} du - \frac{q \co\lambda_0}{u}du + \pi^\ast \ms r,  
\end{equation*}
where we note that $\pi^\ast \ms r \in \h(\widetilde \M) du$ is regular at the origin on the cover.

Although we could not define $t^{\co\lambda_0}$ in $G(\M)$, since $q \co\lambda_0$ is an integral coweight we can define $u^{q \co\lambda_0}$ in $G(\widetilde{\M})$ and compute the corresponding regularisation of the $\g$-connection $\widetilde \nabla$ over the cover, namely
\begin{equation*}
\widetilde{\nabla}_{\rm r} = u^{- q \co\lambda_0} \widetilde{\nabla} u^{q\co\lambda_0} = d + q \sum_{i\in I} u^{q-1+ q\langle \alpha_i, \co\lambda_0\rangle} F_i\, du + \pi^\ast \ms r.
\end{equation*}
This $\g$-connection is clearly regular at the origin, so we may apply the methods developed in \S\ref{ReproGen}. We consider the solution $Y \in B_-(\widetilde{\M})$ of $\widetilde{\nabla}_{\rm r} Y =0$ such that $Y(0)=\Id$, which will be $\widetilde \Gamma$-equivariant but with respect to the automorphism $\vartheta$. Now given any $g_0 \in N^\vartheta$, we perform the $NB_-$ factorisation $Y(u)g_0^{-1}=n(u)^{-1}\widetilde{Y}(u)$ as in Theorem \ref{CycReproReg}. Defining $\widetilde g \coloneqq u^{q\co\lambda_0}n u^{-q\co\lambda_0}$ we have that $\widetilde g\in N^{\widetilde\vs}(\widetilde{\M})$ and $\widetilde g \widetilde{\nabla} \widetilde g^{-1}$ is a cyclotomic Miura $\g$-oper. The $\widetilde\Gamma$-equivariance of $\widetilde g$ means that $\vs\big( \widetilde g(u) \big) = \widetilde g(\widetilde \omega u)$, which implies that $\widetilde g(u) = \widetilde g(\widetilde \omega^T u)$. And since $\widetilde \omega^T$ is a primitive $q^{\rm th}$-root of unity, it therefore follows that $\widetilde g(u) = g\big(\pi(u)\big) = g(t)$ for some $g \in N(\M)$. Moreover, $\vs\big( g(t) \big) = \vs\big( \widetilde g(u) \big) = \widetilde g(\widetilde \omega u) = g(\widetilde \omega^q t) = g(\omega t)$, so that $g \in N^{\hat\vs}(\M)$.
The cyclotomic Miura $\g$-oper $\widetilde g\widetilde{\nabla} \widetilde g^{-1}$ on the cover is therefore the pullback by $\pi$ of the cyclotomic Miura $\g$-oper $g \nabla g^{-1}$.

In summary, even when the dominant coweight $\co\lambda_0 \in \h^\nu$ is not integral but satisfies the weaker condition \eqref{non-integral}, we are still able to construct a new cyclotomic Miura $\g$-oper $g \nabla g^{-1}$ from any given $g_0 \in N^\vartheta$. In this case, $g_0$ can still be interpreted as the initial condition of the regularised gauge transformation parameter $g_{\rm r}$, but where $g_{\rm r}$ has to be considered as an element in $N(\widetilde{\M})$.

\begin{exmp}
Let us illustrate this construction with the help of examples \ref{ReproSl3} and \ref{ReproGenSl3}. We consider the cyclotomic Miura $\sl_3$-oper $\nabla$ defined in equation \eqref{ExpleSl3} but we relax here the condition that $\co\lambda_0$ is integral dominant, \textit{i.e.} that $\eta\in\Z_{\geq 0}$. Instead, we will consider $\eta$ positive and rational and thus write $\eta=\frac{p}{q}$, with $p,q \in \Z_{\geq 1}$ coprime.

We can consider the regularisation $\nabla_{\rm r}$ of $\nabla$ and its extended solution $Y$ as in equations \eqref{Sl3Reg} and \eqref{Sl3SolReg}, but over a $q$-sheeted cover of $\CP$, with coordinates $u=t^{\frac{1}{q}}$, so that the expression $t^\eta$ in these equations can be understood as $u^p$. Let us recall that the $\vartheta$-invariant subgroup $N^\vartheta$ is non-trivial if $\mu=\eta+1$ satisfies one of the following three conditions:
\begin{itemize}
\item $\omega^\mu=1$: In this case, $\mu$ must be integer (a multiple of $T$), hence $\eta$ is integer and we recover the usual setting. \vspace{5pt}

\item $\omega^\mu=-1$: In this case, $\mu$ is also integer (a multiple of $\frac{T}{2}$).\vspace{5pt}

\item $\omega^{2\mu}=-1$: This is the interesting case since $\mu$ can now be half-integer. Working with $\nabla_{\rm r}$ thus requires working in the variable $u=\sqrt{t}$. However, in this case the transformation parameter $g$ is described by the functions $f_i$ as given in equation \eqref{Cas3}. We see that they only depend on $t^{\mu}$ through the expression $t^{2\mu}$ and so are indeed meromorphic functions of $t$.
\end{itemize}
Finally, let us note that, whatever the value of $\eta$ we consider (possibly non-integer), the singular reproduction \eqref{SingSl3} is always meromorphic in $t$.
\end{exmp}

We end this subsection with a discussion of the monodromy around $0$ of the connection $\nabla$. As explained above, working with the cover $\pi : u \mapsto t=u^q$, one can consider a solution $Y \in B_-(\widetilde{\M})$ of the regularised equation $\widetilde{\nabla}_{\rm r}Y = 0$ with $Y(0)=\Id$. And since $\widetilde{\nabla}=u^{q\co\lambda_0}\widetilde{\nabla}_{\rm r}u^{-q\co\lambda_0}$, it follows that $X(u) \coloneqq u^{q\co\lambda_0}Y(u)$ is a solution of the equation $\widetilde{\nabla}X=0$.

Introducing the $q^{\rm th}$-root of unity $\zeta \coloneqq e^{2 \pi i/ q}$, we find that $\widetilde{X}(u) \coloneqq X(\zeta u)$ is also a solution of $\widetilde{\nabla} \widetilde{X}=0$. The monodromy of $\nabla$ around $0$ is then the unique $M_0\in G$ such that $X(\zeta u) = X(u) M_0$. To work out $M_0$ explicitly, we consider the corresponding regularised solution $\widetilde{Y}(u) \coloneqq u^{-q\co\lambda_0}\widetilde{X}(u)$ of $\widetilde{\nabla}_{\rm r} \widetilde{Y}=0$ which is related to $Y$ as $\widetilde{Y}(u)=u^{-q\co\lambda_0}X(\zeta u)=\zeta^{q\co\lambda_0}(\zeta u)^{-q\co\lambda_0}X(\zeta u)=\zeta^{q\co\lambda_0}Y(\zeta u)$. Thus $\widetilde{Y}(0)=\zeta^{q\co\lambda_0}$ and hence $\widetilde{Y} = Y \zeta^{q\co\lambda_0}$, or equivalently $\widetilde X = X \zeta^{q \co\lambda_0}$. Therefore the monodromy of $\nabla$ around the origin is
\begin{equation}
M_0 = \zeta^{q\co\lambda_0} = e^{2\pi i \co\lambda_0}.
\end{equation}
In this equation, the expression $e^{2\pi i \co\lambda_0}$ has to be considered as a formal notation for $\zeta^{q\co\lambda_0}$, which is a well defined element of $H$ as $q\co\lambda_0$ is an integral coweight. When $\co\lambda_0$ is an integral coweight, we find that $M_0=\Id$ as expected from remark \ref{apres proposition}.

\subsection{Isomorphism with $(G/B_-)^\vartheta$} \label{sec: iso flag}

Given any $\nabla \in \MOp_\g(\CP)$, consider the set
\begin{equation*}
\MOp_\g(\CP)_{[\nabla]} \coloneqq \big\{ \widetilde{\nabla} \in \MOp_\g(\CP) \,\big|\, [\widetilde{\nabla}] = [\nabla] \big\}
\end{equation*}
of all Miura $\g$-opers with underlying $\g$-oper $[\nabla]$.

Recall the definition \eqref{nabla reg} of the regularisation of a cyclotomic Miura $\g$-oper at the origin. If $\nabla$ is any Miura $\g$-oper with residue given by minus an integral dominant coweight $- \co\lambda \in \h$ at some point $x \in \C$, then we can similarly introduce its regularisation $\nabla_{\rm r}^x \coloneqq (t - x)^{-\co\lambda} \nabla (t - x)^{\co\lambda}$ at $x$ and for any $g \in N(\M)$ we let $g_{\rm r}^x \coloneqq (t - x)^{-\co\lambda} g (t - x)^{\co\lambda}$. For any $x \in \C$, we introduce the subset of Miura $\g$-opers which are generic at $x$ as, cf. \eqref{generic Miura opers} in the cyclotomic setting,
\begin{equation*}
\MOp_\g(\CP)^{\rm gen}_{\nabla, x} \coloneqq \{ \nabla^g \in \MOp_\g(\CP)_{[\nabla]} \,|\, g_{\rm r}^x \;\text{is regular at}\; x \}.
\end{equation*}
The following theorem is proved in \cite{MV} in the case when $\co\lambda = 0$, \emph{i.e.} when the Miura $\g$-oper $\nabla$ is regular at $x \in \C$, but similar arguments apply also when $\nabla$ has a regular singularity at the given point $x$ by working with its regularisation $\nabla_{\rm r}^x$. Recall the notation of \S\ref{WeylFlag}.

\begin{theorem}\label{MiuraFlag}
Let $\nabla \in \MOp_\g(\CP)$ be monodromy-free. For any $x \in \C$, the following map
\begin{align*}
\Phi: \MOp_\g(\CP)_{[\nabla]} & \longrightarrow G/B_- \\
       \nabla^g & \longmapsto  \lim_{t \rightarrow x} g_{\rm r}^x(t)^{-1} B_-
\end{align*}
is an isomorphism.
Moreover, $\Phi\big( \MOp_\g(\CP)^{\rm gen}_{[\nabla], x} \big) = C_\Id$ and $\overline{\MOp_\g(\CP)^{\rm gen}_{[\nabla], x}} = \MOp_\g(\CP)_{[\nabla]}$. \qed
\end{theorem}

In this section we derive a cyclotomic analog of Theorem \ref{MiuraFlag}. Specifically, from now on we restrict attention to the point $x = 0$, the fixed point of $\C$ under the action of $\Gamma$. Given a \emph{cyclotomic} Miura $\g$-oper $\nabla$ of the form specified at the start of this section, our aim is to describe the image of $\MOp^\Gamma_\g(\CP)_{[\nabla]_\Gamma} \subset \MOp_\g(\CP)_{[\nabla]}$ under the isomorphism $\Phi$.

Since the automorphism $\vs \in \Aut G$ stabilises the Borel subgroup $B_-$, so does the automorphism $\vartheta \in \Aut G$ defined in \eqref{vartheta def}. Hence we obtain a well defined induced action on the quotient $G/B_-$, which by abuse of notation we also denote $\vartheta$, given by
\begin{align*}
\vartheta : G/B_- &\longrightarrow G/B_- \\
              gB_- &\longmapsto  \vartheta(g)B_-.
\end{align*}
Note that the canonical projection $G \to G/B_-$ is $\Gamma$-equivariant with respect to this action, \emph{i.e.} the following diagram is commutative
\begin{equation*}
\begin{tikzpicture}
\matrix (m) [matrix of math nodes, row sep=3em, column sep=3em,text height=1.5ex, text depth=0.25ex]    
{
G & G/B_-\\
G & G/B_-\\
};
\path[->] (m-1-1) edge (m-1-2);
\path[->] (m-2-1) edge (m-2-2);
\path[->] (m-1-1) edge node[left]{$\vartheta$} (m-2-1);
\path[->] (m-1-2) edge node[right]{$\vartheta$} (m-2-2);
\end{tikzpicture}
\end{equation*}
We denote by $(G/B_-)^\vartheta$ the subset of fixed points of $G/B_-$ under $\vartheta$. 

\begin{lemma}\label{PhiEquiv}
Let $\nabla \in \MOp^\Gamma_\g(\CP)$ be monodromy-free. Define $\gamma : N(\M) \to \MOp_\g(\CP)_{[\nabla]}$, $g \mapsto \nabla^g$. The composition $\Phi \circ \gamma : N(\M) \to G/B_-$ is $\Gamma$-equivariant in the sense that the following diagram is commtuative
\begin{equation*}
\begin{tikzpicture}
\matrix (m) [matrix of math nodes, row sep=3em, column sep=3em,text height=1.5ex, text depth=0.25ex]    
{
N(\M) & G/B_-\\
N(\M) & G/B_-\\
};
\path[->] (m-1-1) edge node[above]{$\Phi \circ \gamma$} (m-1-2);
\path[->] (m-2-1) edge node[below]{$\Phi \circ \gamma$} (m-2-2);
\path[->] (m-1-1) edge node[left]{$\hat\vs$} (m-2-1);
\path[->] (m-1-2) edge node[right]{$\vartheta$} (m-2-2);
\end{tikzpicture}
\end{equation*}
In other words, for all $g \in N(\M)$ we have $\Phi(\nabla^{\hat\vs g}) = \vartheta\big( \Phi(\nabla^g) \big)$.
In particular,
\begin{equation*}
\Phi\big( \MOp^\Gamma_\g(\CP)_{[\nabla]_\Gamma} \big) \subset (G/B_-)^\vartheta
\end{equation*}
since $\MOp^\Gamma_\g(\CP)_{[\nabla]_\Gamma} = \gamma\big( N^{\hat\vs}(\M) \big)$.
\begin{proof} Using Proposition \ref{InvY} and the definition and continuity of $\vartheta$, we have
\begin{align*}
\Phi(\nabla^{\hat\vs g}) &= \lim_{t\rightarrow 0} (\hat\vs g)_{\rm r}(t)^{-1} B_-
= \lim_{t\rightarrow 0} \big(\hat\vartheta g_{\rm r} \big)(t)^{-1} B_-
= \lim_{t\rightarrow 0} \vartheta\left[g_{\rm r}(\omega^{-1}t)^{-1}\right] B_-\\
&= \lim_{t\rightarrow 0} \vartheta\left[g_{\rm r}(\omega^{-1}t)^{-1} B_-\right]
= \vartheta \left[ \lim_{t\rightarrow 0} g_{\rm r}(t)^{-1} B_- \right] = \vartheta\big( \Phi(\nabla^g) \big).
\qedhere
\end{align*}
\end{proof}
\end{lemma}

In general, $(G/B_-)^\vartheta$ is \emph{not} isomorphic to the flag variety $G^\vartheta/B_-^\vartheta$. However, we can always realise $G^\vartheta/B_-^\vartheta$ as a subset of $(G/B_-)^\vartheta$. Indeed, we have a well defined injection
\begin{align*}
G^\vartheta/B^\vartheta_- &\longhookrightarrow (G/B_-)^\vartheta\\
              gB^\vartheta_-  &\longmapsto  gB_-.
\end{align*}

The main result of this section is the following cyclotomic version of Theorem \ref{MiuraFlag}.

\begin{theorem}\label{MiuraFlagCyc}
Let $\nabla \in \MOp^\Gamma_\g(\CP)$ be monodromy-free and of the form $\nabla=d+(p_{-1}+u)dt$ with $u\, dt \in \Omega^{\hat{\nu}}(\h)$ as in \eqref{cycl Miura lambda0}. The following map
\begin{align*}
\Phi:  \MOp^\Gamma_\g(\CP)_{[\nabla]_\Gamma} & \longrightarrow (G/B_-)^\vartheta \\
          \nabla^g & \longmapsto  \lim_{t \rightarrow 0} g_{\rm r}(t)^{-1} B_-
\end{align*}
is an isomorphism.
Moreover, we have $\Phi\big( \MOp^\Gamma_\g(\CP)^{\rm gen}_{\nabla, 0} \big) = N^\vartheta B_-/B_-$.
\begin{proof}
This map is the restriction to $\MOp^\Gamma_\g(\CP)_{[\nabla]_\Gamma}$ of the one of Theorem \ref{MiuraFlag} with $x = 0$, and so is injective. Moreover, we know that it is valued in $(G/B_-)^\vartheta$ by Lemma \ref{PhiEquiv}. It therefore remains to prove that it is surjective.

Let $p \in (G/B_-)^\vartheta$. By Theorem \ref{MiuraFlag} there exists a $\nabla^g \in \MOp_\g(\CP)_{[\nabla]}$ such that $\Phi(\nabla^g) = p$. And according to Lemma \ref{PhiEquiv} we have $\Phi(\nabla^{\hat\vs g}) = \vartheta\big( \Phi(\nabla^g) \big) = \vartheta(p) = p = \Phi(\nabla^g)$. The injectivity of $\Phi$ implies that $\nabla^{\hat\vs g} = \nabla^g$. And since the action of $N(\M)$ on $\op_\g(\CP)$ is free by Theorem \ref{TwistedDS}, we deduce that $\hat\vs g = g$. In other words $g \in N^{\hat\vs}(\M)$ so that $\nabla^g \in \MOp^\Gamma_\g(\CP)_{[\nabla]_\Gamma}$, which proves the surjectivity.

The last claim follows at once from Theorem \ref{CycReproReg}.
\end{proof} 
\end{theorem}

\begin{rem}\label{Density}
In general, by contrast with the non-cyclotomic case, $\overline{\MOp^\Gamma_\g(\CP)^{\rm gen}_{\nabla, 0}} \neq \MOp^\Gamma_\g(\CP)_{[\nabla]_\Gamma}$. This feature of the cyclotomic case will be illustrated below in example \ref{FlagSl3}.
\end{rem}

\subsection{Cell decomposition of $(G/B_-)^\upsilon$} \label{sec: cell decomp}

Given any cyclotomic Miura $\g$-oper $\nabla \in \MOp^\Gamma_\g(\CP)$ of the form considered in this section, Theorem \ref{MiuraFlagCyc} establishes an isomorphism between $\MOp^\Gamma_\g(\CP)_{[\nabla]_\Gamma}$ and the $\vartheta$-invariant subvariety $(G/B_-)^\vartheta$ of the flag variety $G/B_-$. Recall the cell decomposition \eqref{Schubert} of the latter. Let $\upsilon \in \Aut \g$ be any automorphism of $\g$, such that $\upsilon^T = \Id$, whose diagram part is given by $\nu$, and which stabilises the Cartan decomposition $\g=\h\oplus\n\oplus\n_-$. An example of such an automorphism is provided by $\vartheta$. In this subsection we establish a similar cell decomposition for the $\upsilon$-invariant subvariety  $(G/B_-)^\upsilon$. We first need the following lemma.

\begin{lemma}
\label{EquivarianceWdot}
The map $\mathcal{K} : W \rightarrow G/B_-$ defined by $\mathcal{K}(w) \coloneqq \dot{w}B_-$ is $\Gamma$-equivariant, \emph{i.e.} the diagram
\begin{equation*}
\begin{tikzpicture}
\matrix (m) [matrix of math nodes, row sep=3em, column sep=3em,text height=1.5ex, text depth=0.25ex]    
{
W & G/B_-\\
W & G/B_-\\
};
\path[->] (m-1-1) edge node[above]{$\mathcal K$} (m-1-2);
\path[->] (m-2-1) edge node[below]{$\mathcal K$} (m-2-2);
\path[->] (m-1-1) edge node[left]{$\nu$} (m-2-1);
\path[->] (m-1-2) edge node[right]{$\upsilon$} (m-2-2);
\end{tikzpicture}
\end{equation*}
is commutative.
\begin{proof}
Recall from \S\ref{WeylFlag} that there is an isomorphism $\pi:W\SimTo N_G(\h)/H$, given by $w\mapsto \dot{w}H$, and that this isomorphism is $\upsilon$-equivariant with respect to the actions of $\upsilon$ on $W$ and $N_G(\h)/H$. Moreover, the map
\begin{align*}
\iota : N_G(\h)/H & \longmapsto G/B_- \\
        gH        & \longmapsto gB_-
\end{align*}
is well defined and $\upsilon$-equivariant. Therefore $\mathcal{K}=\iota\circ\pi$ is also $\upsilon$-equivariant. Finally, note that the action of $\upsilon$ on the Weyl group $W$ only depends on the diagram part $\nu:I\rightarrow I$ of the automorphism $\upsilon$, which coincide with that of $\nu \in \Aut \g$, hence the lemma.
\end{proof}
\end{lemma}

To obtain the cell decomposition of $(G/B_-)^\upsilon$ we begin by recalling some facts about the Schubert cells of $G/B_-$. Let $w$ be an element of the Weyl group $W$. We define the subset of positive roots $R(w) \coloneqq \{ \alpha\in\Phi^+ \, | \, w^{-1}\alpha \in \Phi^- \}$. Its cardinality $|R(w)|$ is equal to the length $l(w)$ of the Weyl group element $w$. Consider the subgroup $U_w \coloneqq \prod_{\alpha\in R(w)} G_\alpha$ of $N$ which does not depend on the chosen ordering on the roots $\alpha \in R(w)$.

Let $w_\circ$ be the longest element of $W$ (the unique element of $W$ of maximal length). Then the cell $C_w$ is isomorphic to $U_{w_\circ w}$. More precisely, every point of $C_w$ is of the form $n\dot{w}B_-=n\mathcal{K}(w)$, with $n\in U_{w_\circ w}$, and this form is unique. The following lemma describes the action of the automorphism $\upsilon$ on the subgroups $U_w$.

\begin{lemma}\label{EquiUw}
Let $w \in W$. We have $R(\nu w)= \nu\bigl(R(w)\bigr)$ and $\upsilon(U_w)=U_{\nu w}$. In particular, if $w \in W^\nu$ then $\upsilon$ stabilises $U_w$. Moreover, $\nu w_\circ = w_\circ$.
\begin{proof}
Let $\alpha\in R(w)$. By definition of $\nu w$, we have $(\nu w)^{-1} (\nu\alpha) = \nu \circ w^{-1} \circ \nu^{-1} (\nu \alpha)=\nu (w^{-1}\alpha)$. Yet, $\alpha\in R(w)$, hence $w^{-1}\alpha \in\Phi^-$. Thus $(\nu w)^{-1} (\nu\alpha)=\nu (w^{-1}\alpha)$ belongs to $\Phi^-$, \emph{i.e.} $\nu\alpha \in R(\nu w)$. Therefore, $\nu\bigl(R(w)\bigr) \subset R(\nu w)$. Applying this result to $\nu w \in W$ and the automorphism $\nu^{-1}$, we get $\nu^{-1}\bigl(R(\nu w)\bigr) \subset R(w)$, hence (applying $\nu$ to this inclusion) $R(\nu w) \subset \nu\bigl( R(w) \bigr)$. By double inclusion, we conclude that $R(\nu w)=\nu\bigl(R(w)\bigr)$. The claim that $\upsilon(U_w)=U_{\nu w}$ follows directly, noting that, for any root $\alpha\in\Phi$, $G_{\nu\alpha}=\nu(G_\alpha)=\upsilon(G_\alpha)$ (see equation \eqref{AdaptedBasis}).

As $R(\nu w)=\nu\bigl(R(w)\bigr)$, we have $l(\nu w) = |R(\nu w)| = |R(w)| = l(w)$, \emph{i.e.} the action of $\nu$ on the Weyl group $W$ preserves the length. In particular, $l(\nu w_\circ)=l(w_\circ)$. As $w_\circ$ is the unique element of $W$ of maximal length, we deduce that $\nu w_\circ = w_\circ$. 
\end{proof}
\end{lemma}

We are now in a position to state and prove the sought after cell decomposition of $(G/B_-)^\upsilon$.

\begin{theorem}\label{CellUpsilon}
The space $(G/B_-)^\upsilon$ admits the following cell decompositon
\begin{equation}\label{SchubertTheta}
(G/B_-)^\upsilon = \displaystyle \bigsqcup_{w \in W^\nu} N^\upsilon \dot{w} B_-/B_- =: \bigsqcup_{w \in W^\nu} C^\upsilon_w.
\end{equation}
Moreover, each element of the cell $C^\upsilon_w$ can be written uniquely in the form $n\dot{w}B_-$ with $n \in U_{w_\circ w}^\upsilon$.
\begin{proof}
Let $w\in W^\nu$ and $n\in N^\upsilon$. Using Lemma \ref{EquivarianceWdot}, we have $\upsilon\bigl(\mathcal{K}(w)\bigr)=\mathcal{K}(\nu w)= \mathcal{K}(w)$ and hence $\upsilon(n\dot{w}B_-)=\upsilon(n) \upsilon\bigl(\mathcal{K}(w)\bigr) = n \mathcal{K}(w) = n\dot{w}B_-$. In other words, $n\dot{w}B_-$ belongs to $(G/B_-)^\upsilon$. This proves the inclusion $\displaystyle \bigsqcup_{w \in W^\nu} N^\upsilon \dot{w} B_-/B_- \subset (G/B_-)^\upsilon$.

Conversely, let us consider an element $p \in (G/B_-)^\upsilon$. As an element of $G/B_-$, it belongs to a cell $C_w$, for some $w \in W$. Thus, we can write $p=n\dot{w}B_-=n\mathcal{K}(w)$ for some unique $n \in U_{w_\circ w}$. Using Lemma \ref{EquivarianceWdot} we get $p=\upsilon(p)=\upsilon(n) \upsilon\bigl( \mathcal{K}(w) \bigr) = \upsilon(n) \mathcal{K}(\nu w)$. As $\upsilon$ stabilises $N$, $\upsilon(n) \in N$ and so $p \in C_{\nu w}$. Yet, the cells $C_{w'}$ are disjoint in $G/B_-$, hence $w=\nu w$, \emph{i.e.} $w\in W^\nu$. Lemma \ref{EquiUw} implies that $w_\circ w \in W^\nu$ and that $U_{w_\circ w}$ is stabilised by $\upsilon$. Thus, we have $p=n\mathcal{K}(w)=\upsilon(n)\mathcal{K}(w)$, with $n$ and $\upsilon(n)$ in $U_{w_\circ w}$. By unicity of this form, we deduce that $n \in U_{w_\circ w}^\upsilon$. This provides the opposite inclusion $(G/B_-)^\upsilon \subset \displaystyle \bigsqcup_{w \in W^\nu} N^\upsilon \dot{w} B_-/B_-$ and proves also the last statement.
\end{proof}
\end{theorem}

\begin{exmp}\label{FlagSl3}
Consider the cyclotomic Miura $\sl_3$-oper $\nabla$ studied in examples \ref{Sl3}, \ref{ReproSl3} and \ref{ReproGenSl3}. We will illustrate Theorem \ref{MiuraFlagCyc} for $\nabla$, with the help of Theorem \ref{CellUpsilon}. The $\nu$-invariant Weyl group $W^\nu$ is composed of two elements: the identity and the simple reflection $s^\nu=s_1s_2s_1$ associated with the simple folded root $\alpha^\nu=\frac{1}{2}(\alpha_1+\alpha_2)$. By Theorem \ref{CellUpsilon}, the variety $(G/B_-)^\vartheta$ must be composed of two cells $C^\vartheta_\Id$ and $C^\vartheta_{s^\nu}$.

In the examples \ref{ReproSl3} and \ref{ReproGenSl3} we described all possible cyclotomic reproductions of $\nabla$. They are of the following two types:
\begin{itemize}
\item the regular ones, characterised by a choice of initial condition $g_0 \in N^\vartheta$,
\item the (unique) singular one, given by equation \eqref{SingSl3}.
\end{itemize}

As explained in \S\ref{ReproGen}, the regular reproductions generate the subset of generic cyclotomic Miura $\sl_3$-opers $\MOp^\Gamma_\g(\CP)^{\rm gen}_{\nabla, 0}$. Moreover, it follows from Theorem \ref{MiuraFlagCyc} that this subset corresponds to the big cell $C^\vartheta_\Id=N^\vartheta B_-/B_-$ of $(G/B_-)^\vartheta$. More precisely, if $g \in N^{\hat\vs}(\M)$ is such that $g_{\rm r}$ is regular with initial condition $g_{\rm r}(0)=g_0 \in N^\vartheta$ then $\Phi(g\nabla g^{-1})=g_0 B_-$.

The singular reproduction, with parameter $g$ given by \eqref{SingSl3}, gives rise to a non-generic cyclotomic Miura $\g$-oper $g\nabla g^{-1}$, corresponding to the point $\Phi(g\nabla g^{-1})=\dot{s}^\nu B_-$ in the variety $(G/B_-)^\vartheta$ forming the one-point cell $C^\vartheta_{s^\nu}=\lbrace \dot{s}^\nu B_- \rbrace$.

As explained in example \ref{ReproGenSl3}, the $\vartheta$-invariant subgroup $N^\vartheta$ is non-trivial if and only if $\omega^{4(\eta+1)}=1$. In the case where it is non-trivial, the singular reproduction \eqref{SingSl3} can be obtained as a limit of the regular ones (cf. example \ref{ReproSl3}): in the flag variety, this corresponds to the density of the big cell $N^\vartheta B_-/B_-$ in the variety $(G/B_-)^\vartheta$. However, when $N^\vartheta$ is trivial, the variety $(G/B_-)^\vartheta$ is composed of only two points: $B_-$ and $\dot{s}B_-$. In this case, the big cell $\lbrace B_- \rbrace$ is not dense in $(G/B_-)^\vartheta$, which illustrates the point made in remark \ref{Density} that in general $\overline{\MOp^\Gamma_\g(\CP)^{\rm gen}_{\nabla, 0}} \neq \MOp^\Gamma_\g(\CP)_{[\nabla]_\Gamma}$.
\end{exmp}

\section{Cyclotomic Gaudin models}
\label{GaudinModels}

In this section we use the notion of cyclotomic $\g$-oper from \S\ref{DS} to formulate the conjecture relating the (maximal) spectrum of the cyclotomic Gaudin algebra $\mathscr Z_{(z_i)}^\Gamma(\g)$, introduced in \cite{VY1} and whose definition was briefly recalled in the introduction, to a certain space of cyclotomic $\Lg$-opers for the Langlands dual Lie algebra $\Lg$ of $\g$.

We begin recalling the Bethe ansatz solution of the cyclotomic Gaudin model obtained in \cite{VY1} when the spin chain is a tensor product of Verma modules. The joint eigenvectors of the cyclotomic Gaudin algebra $\mathscr Z_{(z_i)}^\Gamma(\g)$ are characterised in this case by certain $\Gamma$-equivariant rational functions $\CP \to \h^\ast$ valued in weight space $\h^\ast$ with poles at the points $z_i$, $i =1, \ldots, N$ of the cyclotomic Gaudin model (and their $\Gamma$-orbits), at a further collection of \emph{Bethe roots} $x_j$, $j = 1, \ldots, m$ (and their $\Gamma$-orbits) satisfying the cyclotomic Bethe ansatz equations, as well as an additional pole at the origin with residue given by a special weight depending on the cyclotomic Gaudin model datum.

Noting that $\h^\ast$ is canonically identified with the Cartan subalgebra of the Langlands dual Lie algebra $\Lg$, such a rational function $\CP \to \h^\ast$ can be used to define a cyclotomic Miura $\Lg$-oper. We conjecture in this case that the eigenvalues of the cyclotomic Gaudin Hamiltonians can be read off from the canonical form of the underlying cyclotomic $\Lg$-oper.

\subsection{Bethe ansatz} \label{sec: Bethe}

Let $\g$ be a finite-dimensional complex semisimple Lie algebra. We follow the notation and conventions of \S\ref{SSLieAlgebras}. As in \S\ref{DS} we let $T \in \Z_{\geq 1}$, we pick a primitive $T^{\rm th}$-root of unity $\omega$ and let $\Gamma \coloneqq \langle \omega \rangle \cong \Z/ T \Z$ be a copy of the cyclic group of order $T$ acting on $\CP$ by multiplication. Let $\sigma \in \Aut \g$ be such that $\sigma^T = \Id$ with diagram part $\nu : I \to I$. We can always choose a Cartan subalgebra $\h$ adapted to $\sigma$ with corresponding Chevalley-Serre generators $E_i$, $\co\alpha_i$, $F_i$ for $i \in I$ such that, see \emph{e.g.} \cite{Kac},
\begin{equation*}
\sigma(E_i) = \tau_{\alpha_i} E_{\nu(i)}, \qquad
\sigma(\co\alpha_i) = \co\alpha_{\nu(i)}, \qquad
\sigma(F_i) = \tau_{\alpha_i}^{-1} F_{\nu(i)},
\end{equation*}
for some $T^{\rm th}$-roots of unity $\tau_{\alpha_i} \in \Gamma$.

Let $N \in \Z_{\geq 1}$ and fix a set of $N$ distinct finite non-zero points $z_i \in \C^\times$, $i = 1, \ldots, N$ with disjoint $\Gamma$-orbits. Let $u \in \C^\times \setminus \{ z_1, \ldots, z_N \}$ be any other point whose $\Gamma$-orbit is also disjoint from those of the $z_i$. The \emph{cyclotomic Gaudin algebra} $\mathscr Z_{(z_i)}^\Gamma(\g) \subset U(\g)^{\otimes N}$ is defined in \cite{VY1} as the image of a homomorphism of commutative $\C$-algebras
\begin{equation} \label{hom V0 to ZG 7}
\Psi^\Gamma_{(z_i), u} : \mathfrak z \big( \mathbb V_{0,u}^{\rm crit}(\g) \big) \longrightarrow U(\g)^{\otimes N},
\end{equation}
where $\mathfrak z \big( \mathbb V_{0,u}^{\rm crit}(\g) \big)$ denotes the space of singular vectors in the vacuum Verma module $\mathbb V_{0,u}^{\rm crit}(\g)$ over the affine Kac-Moody algebra $\widehat \g$ at the \emph{critical level}.
The construction of \eqref{hom V0 to ZG 7} in \cite{VY1}, which generalises arguments of \cite{FFR} to the cyclotomic setting, makes essential use of the notion of \emph{cyclotomic coinvariants} of a tensor product of $\widehat \g$-modules introduced in \cite{VY2} (see also \cite{FreSzcz}).

Pick a collection of weights $\lambda_i \in \h^\ast$ for $i = 1,\ldots, N$. Let $M_\lambda \coloneqq U(\g) \otimes_{U(\b)} \C \, v_\lambda$ be the Verma module of highest weight $\lambda \in \h^\ast$. The problem of diagonalising the cyclotomic Gaudin algebra $\mathscr Z_{(z_i)}^\Gamma(\g)$ on the tensor product $\bigotimes_{i=1}^N M_{\lambda_i}$ was addressed in \cite{VY1} by using a generalisation of the Bethe ansatz construction of \cite{FFR} to the case $\Gamma \neq \{ 1 \}$. The central idea of this approach is to utilise the rich structure of modules over the affine Kac-Moody algebra $\widehat\g$ at the critical level. In particular, one has access to the so called \emph{Wakimoto modules} which are parametrised by $\h$-valued formal Laurent series. The eigenvectors of the cyclotomic Gaudin algebra $\mathscr Z_{(z_i)}^\Gamma(\g)$ in $\bigotimes_{i=1}^N M_{\lambda_i}$ are then obtained from cyclotomic coinvariants of tensor products of such modules. The main novelty in the cyclotomic setting compared to the usual case is the need to also introduce a Wakimoto module over the twisted affine Kac-Moody algebra $\widehat \g^\sigma$ at the origin, \emph{i.e.} the fixed point of the action of $\Gamma$ in $\C$. Below we will only recall the end result of the construction, referring the interested reader to \cite{VY1, VY2} and \cite{FFR} for further details.

Let $m \in \Z_{\geq 0}$ and suppose we are given any map $c : \{ 1, \ldots, m \} \to I$, the so called \emph{colour function}. The corresponding \emph{cyclotomic Bethe ansatz equations} are a collection of $m$ algebraic equations on a subset of $m$ points $\{ x_j \}_{j=1}^m \subset \C^\times$ given by
\begin{equation}\label{Bethe}
0 = \sum_{r=0}^{T-1} \sum_{i=1}^N \frac{(\alpha_{c(j)} | \nu^r \lambda_i )}{x_j-\omega^r z_i} - \underset{(r, k)\neq (0, j)}{\sum_{r=0}^{T-1} \sum_{k=1}^m} \frac{( \alpha_{c(j)} | \nu^r \alpha_{c(k)} )}{x_j-\omega^r x_k} + \frac{1}{x_j} (\alpha_{c(j)}| \lambda_0), \qquad j = 1, \ldots, m.
\end{equation}
Here the \emph{weight at the origin} $\lambda_0 \in \h^\ast$ is defined in terms of the automorphism $\sigma \in \Aut \g$ and the root of unity $\omega$ as
\begin{equation} \label{lambda0 def}
\lambda_0(h) \coloneqq \sum_{r=1}^{T-1} \frac{\tr_\n (\sigma^{-r} \circ \ad_h)}{1 - \omega^r}.
\end{equation}
One checks that this weight is $\nu$-invariant using the equality $\ad_{\nu h} = \sigma \circ \ad_h \circ \sigma^{-1}$ as endomorphisms of the nilpotent Lie algebra $\n$, together with the cyclicity of the trace $\tr_\n$ over $\n$.

To any solution of \eqref{Bethe} we can associate an eigenvector $\psi_{(x_j)} \in \bigotimes_{i=1}^N M_{\lambda_i}$ of the cyclotomic Gaudin algebra $\mathscr Z^\Gamma_{(z_i)}(\g)$. Such a \emph{Bethe vector} can be expressed as a linear combination of vectors of the form
$\bigotimes_{i=1}^N F_{c(p^i_1)} \ldots F_{c(p^i_{n_i})} v_{\lambda_i}$, where $n_1 + \ldots + n_N = m$ is a composition of the integer $m$ and $(p^1_1, \ldots, p^1_{n_1}; \ldots; p^N_1, \ldots, p^N_{n_N})$ is a permutation of $(1, 2, \ldots, m)$, with the coefficients depending rationally on all the Bethe roots $x_j$, $j = 1, \ldots, m$ and on the points $z_i$, $i = 1, \ldots, N$. See \cite[Proposition 4.6]{VY1} for the explicit expression of $\psi_{(x_j)}$.
The eigenvalues of the family of $N$ quadratic cyclotomic Gaudin Hamiltonians, cf. \eqref{CycloH}, on this eigenvector read
\begin{equation}\label{Energy}
\mathcal E_i = \underset{(r, k)\neq (0, j)}{\sum_{r=0}^{T-1} \sum_{j=1}^N} \frac{(\lambda_i | \nu^r \lambda_j )}{z_i - \omega^r z_j} - \sum_{r=0}^{T-1} \sum_{j=1}^m \frac{(\lambda_i| \nu^r \alpha_{c(j)} )}{z_i - \omega^r x_j} + \frac{1}{z_i} (\lambda_i | \lambda_0 ), \qquad i = 1, \ldots, N.
\end{equation}
In other words, we have $\mathcal H_i \psi_{(x_j)} = \mathcal E_i \psi_{(x_j)}$ for all $i =1, \ldots, N$ provided the cyclotomic Bethe ansatz equations \eqref{Bethe} hold.
Note that the question of whether the Bethe vectors $\psi_{(x_j)}$ so constructed are non-zero was not addressed in \cite{VY1}. This problem was considered recently in \cite{CharlesVar2} where it was proved, at least in the case when $\sigma \in \Aut \g$ is a diagram automorphism, that the Bethe vectors $\psi_{(x_j)}$ are indeed non-zero (see also \cite{SchVar, MV2, ChSch} in the non-cyclotomic case).

\subsection{Cyclotomic Miura $\Lg$-opers} \label{sec: Miura Gaudin}

In the non-cyclotomic case the eigenvectors $\psi_{(x_j)}$ obtained by the Bethe ansatz, or Wakimoto construction, are naturally parameterised by an $\h^\ast$-valued connection on $\CP$ with poles at each $z_i$, $i = 1, \ldots, N$ with residue $- \lambda_i$ and at each Bethe root $x_j$, $j = 1, \ldots, m$ with residue $\alpha_{c(j)}$ \cite{Frenkel2}. Likewise, in the cyclotomic analog of the Wakimoto construction \cite{VY1}, each Bethe vector $\psi_{(x_j)}$ corresponds to a $\Gamma$-equivariant $\h^\ast$-valued meromorphic differential on $\CP$,
\begin{equation}\label{lambda}
\lambda(t) dt \coloneqq \frac{\lambda_0}{t} dt + \sum_{r=0}^{T-1} \left( \sum_{i=1}^N \frac{\nu^r \lambda_i}{t-\omega^rz_i} - \sum_{j=1}^m \frac{\nu^r \alpha_{c(j)}}{t-\omega^r x_j} \right) dt \in \Omega^{\hat \nu}(\h^\ast).
\end{equation}
In particular, taking the Laurent expansion of $\lambda(t)$ at $x_j$ we find
\begin{equation*}
\lambda(t) = - \frac{\alpha_{c(j)}}{t - x_j} + r(t)
\end{equation*}
where $r \in \h(\M)$ is regular at $x_j$, and the cyclotomic Bethe ansatz equations \eqref{Bethe} then simply read $(\alpha_{c(j)} | r(x_j)) = 0$. Similarly, the energies $\mathcal E_i$ of the quadratic cyclotomic Gaudin Hamiltonians $\mathcal H_i$ are given in terms of \eqref{lambda} simply by
\begin{equation}\label{lambda squared}
\mathcal E_i = \res_{z_i} \big( \ha ( \lambda(t) | \lambda(t) ) - ( \lambda'(t) | \rho ) \big) dt.
\end{equation}
Both of these statements can be formulated as properties of a cyclotomic Miura $\Lg$-oper built from \eqref{lambda}, where $\Lg$ denotes the \emph{Langlands dual} of $\g$.

The set $\Phi^\vee \coloneqq \{ \co \alpha \in \h \,|\, \alpha \in \Phi \}$ of all coroots forms a root system in $\h$.
The Langlands dual of $\g$ is defined as the finite-dimensional complex semisimple Lie algebra $\Lg$ with this root system $\Phi^\vee$. Its Cartan matrix is then the transpose of that of $\g$. The Cartan subalgebra of $\Lg$ is canonically identified with the dual $\h^*$ of $\h$. Denote by $\co E_i$, $\alpha_i$, $\co F_i$ for $i \in I$ the Chevalley-Serre generators of $\Lg$. We let $\omega_i$, $i \in I$ denote the basis of fundamental weights of $\g$, \emph{i.e.} fundamental coweights of $\Lg$.

As in \S\ref{Principal}, the principal $\sl_2$-triple is defined by its regular nilpotent and semisimple elements
\begin{equation*}
\co p_{-1} \coloneqq \sum_{i \in I} \co F_i, \qquad
\rho \coloneqq \sum_{i \in I} \omega_i.
\end{equation*}
We denote this $\sl_2$-triple by $\{ \co p_{-1}, \rho, \co p_1 \}$.
To define the notion of cyclotomic $\Lg$-oper we introduce a specific automorphism $\vs \in \Aut \Lg$ of $\Lg$, as in \S\ref{sec: cyclo op}, by letting
\begin{equation*}
\vs(\co E_i) = \omega^{-1} \co E_{\nu(i)}, \qquad
\vs(\alpha_i) = \alpha_{\nu(i)}, \qquad
\vs(\co F_i) = \omega \co F_{\nu(i)},
\end{equation*}
for all $i \in I$. We may use \eqref{lambda} to define a cyclotomic Miura $\Lg$-oper as
\begin{equation}\label{lambda Miura}
\nabla \coloneqq d + \co p_{-1} dt - \lambda(t) dt \in \MOp_{\Lg}^\Gamma(\CP).
\end{equation}
It follows from Proposition \ref{RegBetheSimple} that the cyclotomic Bethe ansatz equations \eqref{Bethe} are precisely the conditions for the underlying cyclotomic $\Lg$-oper $[\nabla]_\Gamma$ to be regular at all the Bethe roots $x_j$, $j = 1, \ldots, m$, see also remark \ref{apres proposition 0}. We define the \emph{weight at infinity} as the unique $\nu$-invariant weight $\lambda_\infty \in \h^{\ast, \nu}$ in the shifted $W^\nu$-orbit of $- \res_\infty \lambda(t) dt$ such that $\lambda_\infty + \rho$ is dominant. We then have
\begin{equation}\label{lambda inf}
\lambda_0 + \sum_{r=0}^{T-1} \left( \sum_{i=1}^N \nu^r \lambda_i - \sum_{j=1}^m \nu^r \alpha_{c(j)} \right) = w_\infty \cdot \lambda_\infty
\end{equation}
for some $w_\infty \in W^\nu$.
We deduce at once from Theorem \ref{GeneralForm} that the cyclotomic $\Lg$-oper $[\nabla]_\Gamma$ lives in the space of $\Lg$-opers $\Op_{\Lg}^\Gamma(\CP)^{\text{RS}}_{\bm z; \bm \lambda}$ introduced at the end of \S\ref{DS}, where $\bm z \coloneqq \{ 0, z_1, \ldots, z_N, \infty \} \subset \CP$ and $\bm \lambda \coloneqq \{ \lambda_0, \lambda_1, \ldots, \lambda_N, \lambda_\infty \} \subset \h^\ast$. If, moreover, the weights $\lambda_i \in \h^\ast$, $i = 1,\ldots, N$ and $\nu$-invariant weight $\lambda_0 \in \h^{\ast, \nu}$ are all integral dominant, then so are all the weights $\nu^r \lambda_i \in \h^\ast$ for $r=0, \ldots, T-1$ and $i = 1,\ldots, N$ from which it follows that $\nabla$ is monodromy-free, cf. remark \ref{apres proposition}. In this case $[\nabla]_\Gamma$ defines an element of the subset $\Op_{\Lg}^\Gamma(\CP)_{\bm z; \bm \lambda} \subset \Op_{\Lg}^\Gamma(\CP)^{\text{RS}}_{\bm z; \bm \lambda}$ introduced at the end of \S\ref{DS}. In fact, it follows from \S\ref{sec: non-integral} that the integrality assumption on $\lambda_0$ can be relaxed to $\langle \co\alpha_i, \lambda_0 \rangle \in \mathbb{Q}$ for all $i \in I$, in which case $\nabla$ has monodromy $e^{2 \pi i \lambda_0} \in \null^L G$ at the origin, where $\null^L G$ denotes the adjoint group of the Langlands dual Lie algebra $\Lg$. The meromorphic $\Lg$-connection $\nabla$ becomes a monodromy-free when pulled back to a cover of $\CP$.

\subsection{Spectrum of $\mathscr Z_{(z_i)}^\Gamma(\g)$}

In this final section we state two conjectures relating the spectrum of the cyclotomic Gaudin model to cyclotomic $\Lg$-opers. The first gives a description of the common eigenvalues of elements of the cyclotomic Gaudin algebra $\mathscr Z_{(z_i)}^\Gamma(\g)$ on a joint eigenvector $\psi_{(x_j)}$ in the tensor product of Verma modules $\bigotimes_{i=1}^N M_{\lambda_i}$, given the corresponding cyclotomic Miura $\Lg$-oper.
The second describes the maximal spectrum of the cyclotomic Gaudin algebra $\mathscr Z_{(z_i)}^\Gamma(\g)$ in terms of cyclotomic $\Lg$-opers with regular singularities at the $z_i$, $i = 1, \ldots, N$, the origin and infinity.

Recall that the center $Z(\g) \subset U(\g)$ of the enveloping algebra $U(\g)$ is isomorphic to the polynomial algebra $\C[P_k]_{k \in E}$ in $\rk \g$ generators, where $P_k$ has degree $k+1$ in the canonical filtration on $U(\g)$. Recall here that $E$ denotes the \emph{multiset} of exponents of the semisimple Lie algebra $\g$. Similarly, in the affine setting it follows from the results of Feigin and Frenkel \cite{FF1} that, as a commutative algebra, the subspace of singular vectors $\mathfrak z\big( \mathbb V^{\rm crit}_{0, u}(\g) \big) \subset \mathbb V^{\rm crit}_{0, u}(\g)$ in the vacuum Verma module $\mathbb V^{\rm crit}_{0, u}(\g)$ at critical level over $\widehat \g_u$ is isomorphic to the polynomial algebra in infinitely many variables
\begin{equation*}
\mathfrak z\big( \mathbb V^{\rm crit}_{0, u}(\g) \big) \simeq \C[\mathcal D^r S_k]_{r \in \Z_{\geq 0}; k \in E},
\end{equation*}
where $\mathcal D \in \End \mathbb V^{\rm crit}_{0, u}(\g)$ is defined by the property $[\mathcal D, a(n)] = - n\, a(n-1)$ for all $a \in \g$ and $n \in \Z$, where $a(n) \coloneqq a \otimes (t-u)^n \in \widehat \g_u$, and $\mathcal D v_0 = 0$. The element $S_k \in \mathbb V^{\rm crit}_{0, u}(\g)$ is of degree $k + 1$ in the natural filtration on the $\widehat \g_u$-module $\mathbb V^{\rm crit}_{0, u}(\g)$. 

\begin{conj} \label{conj: Miura oper}
Let $\nabla \in \MOp^\Gamma_{\Lg}(\CP)$ correspond to an eigenvector $\psi_{(x_j)}$ of the cyclotomic Gaudin algebra $\mathscr Z^\Gamma_{(z_i)}(\g)$ on the tensor product of Verma modules $\bigotimes_{i=1}^N M_{\lambda_i}$, as considered in \S\ref{sec: Bethe} and \S\ref{sec: Miura Gaudin}. Let $d + \co p_{-1} dt + \sum_{k \in E} \ms c_k$ be the canonical representative of the underlying
cyclotomic $\Lg$-oper $[\nabla]_\Gamma$.
There exists a basis\footnote{Note that bases of the centraliser of $\co p_1$ in $\Lg$ are naturally labelled by elements of the multiset of exponents $\co E$ of $\Lg$, cf. \S\ref{Principal}. Here we implicitly make use of the fact that $\co E$ coincides with the multiset $E$ of exponents of $\g$.} $\{ \co p_i \}_{i \in E}$ of the centraliser of $\co p_1$ in $\Lg$, which is independent of $\nabla$, such that for every $r \in \Z_{\geq 0}$ and $k \in E$, the eigenvalue of $\Psi^\Gamma_{(z_i), t}(\mathcal D^r S_k)$ on $\psi_{(x_j)}$ is given by $\partial_t^r u_k(t)$, where $\ms c_k = u_k(t) \, \co p_k \, dt$.
\end{conj}

To motivate the conjecture, consider the quadratic singular vector $S \coloneqq S_1 = \ha I_a(-1) I^a(-1) v_0$ in $\mathbb V^{\rm crit}_{0, u}(\g)$ corresponding to the exponent $1$. The quadratic cyclotomic Gaudin Hamiltonians are obtained from it as $\mathcal H_j = \res_{z_j} \Psi^\Gamma_{(z_i), t}(S) dt$ for each $j = 1, \ldots, N$. According to Proposition \ref{u1}, the component $\ms c_1 \in \Omega^{\hat\vs}(\a)$ of the canonical representative of $\nabla = d + \co p_{-1} dt - \lambda(t) dt \in \MOp^\Gamma_{\Lg}(\CP)$ is given by $\ms c_1 = u_1 \, \co p_1 dt$ where
\begin{equation*}
u_1 = \frac{1}{2 (\co\rho | \co\rho)} \big( \ha ( \lambda(t) | \lambda(t)) - (\lambda'(t) | \rho )\big).
\end{equation*}
By suitably normalising the basis element $\co p_1$ we can remove the overall factor of $1/ 2 (\co\rho | \co\rho)$ from $u_1$. It then follows from \eqref{lambda squared} that the eigenvalue of $\mathcal H_i$ on $\psi_{(x_j)}$ is given by the residue $\res_{z_i} u_1(t) dt$, in agreement with Conjecture \ref{conj: Miura oper} .

At the end of \S\ref{DS} we defined the set $\Op_{\Lg}^\Gamma(\CP)^{\text{RS}}_{\bm z}$ of cyclotomic $\Lg$-opers with regular singularities at the points in the set $\bm z = \{ 0, z_1, \ldots, z_N, \infty \}$. Let $\lambda_0 \in \h^{\ast, \nu}$ be the $\nu$-invariant weight defined in \eqref{lambda0 def}. We shall assume that $\lambda_0$ is dominant and such that $\langle \co\alpha_i, \lambda_0 \rangle \in \mathbb Q$ for all $i \in I$, as in the setting of \S\ref{sec: non-integral}. Let $\Op_{\Lg}^\Gamma(\CP)^{\text{RS}}_{\bm z; \lambda_0}$ be the subset of cyclotomic $\Lg$-opers in $\Op_{\Lg}^\Gamma(\CP)^{\text{RS}}_{\bm z}$ whose residue at the origin is given by the finite $\Lg^\nu$-oper $[\lambda_0]_{W^\nu}$ and with monodromy $e^{2 \pi i \lambda_0} \in \null^L G$ at the origin.

\begin{conj} \label{conj: cyclo opers}
The cyclotomic Gaudin algebra $\mathscr Z^\Gamma_{(z_i)}(\g)$ is isomorphic to the algebra of functions on the space $\Op_{\Lg}^\Gamma(\CP)^{\textup{RS}}_{\bm z; \lambda_0}$. In other words,
\begin{equation*}
\Spec \mathscr Z_{(z_i)}^\Gamma(\g) \simeq \Op_{\Lg}^\Gamma(\CP)^{\textup{RS}}_{\bm z; \lambda_0}.
\end{equation*}
\end{conj}

In the non-cyclotomic case where $\Gamma = \{ 1 \}$, the cyclotomic Gaudin algebra $\mathscr Z^\Gamma_{(z_i)}(\g)$ reduces to the ordinary Gaudin algebra $\mathscr Z_{(z_i)}(\g)$. On the other hand, since $T = 1$ the weight at the origin becomes trivial, \emph{i.e.} $\lambda_0 = 0$. Hence $\Op_{\Lg}^\Gamma(\CP)^{\textup{RS}}_{\bm z; \lambda_0}$ consists of $\Lg$-opers with regular singularity at the points of $\bm z$, with trivial monodromy at the origin and residue there equal to $[0]_{W^\nu}$. It follows from \cite[Lemma 2.4]{Frenkel} that such $\Lg$-opers are regular at the origin. Thus Conjecture \ref{conj: cyclo opers} reduces in the case $\Gamma = \{ 1 \}$ to \cite[Theorem 2.7$(1)$]{Frenkel2} describing the spectrum of the usual Gaudin algebra $\mathscr Z_{(z_i)}(\g)$.

\providecommand{\bysame}{\leavevmode\hbox to3em{\hrulefill}\thinspace}
\providecommand{\MR}{\relax\ifhmode\unskip\space\fi MR }
\providecommand{\MRhref}[2]{%
  \href{http://www.ams.org/mathscinet-getitem?mr=#1}{#2}
}
\providecommand{\href}[2]{#2}

\end{document}